\documentclass{article}
\usepackage[utf8]{inputenc}
\usepackage{amsmath}
\usepackage{amsfonts}
\usepackage{tikz}
\usepackage{xcolor}
\usepackage{newtxtext,newtxmath}
\usepackage{caption}
\usepackage{subcaption}
\usepackage{amsmath}
\usepackage{hyperref}
\usepackage{authblk}
\usepackage{adjustbox}
\usepackage{float}

\usepackage{amsthm}

\title{An optimisation--based domain--decomposition reduced order model for the incompressible Navier-Stokes equations}

\date{}
\author[a]{Ivan Prusak\footnote{\href{mailto:iprusak@sissa.it}{iprusak@sissa.it}}}

\author[b]{Monica Nonino\footnote{\href{mailto:monica.nonino@univie.ac.at}{monica.nonino@univie.ac.at}}}

\author[a]{Davide Torlo\footnote{\href{mailto:dtorlo@sissa.it}{dtorlo@sissa.it}}}
\author[c]{Francesco Ballarin\footnote{\href{mailto:francesco.ballarin@unicatt.it}{francesco.ballarin@unicatt.it}}}
\author[a]{Gianluigi Rozza \footnote{\href{mailto:grozza@sissa.it}{grozza@sissa.it}}} 

\affil[a]{Mathematics Area, mathLab, SISSA, 34136 Trieste, Italy}

\affil[b]{Fakultät für Mathematik,
Universität Wien, 1090 Wien, Austria}
            
\affil[c]{Department of Mathematics and Physics, Università Cattolica del Sacro Cuore, 25133 Brescia, Italy, }

\providecommand{\keywords}[1]
{
  \small	
  \textbf{\textit{Keywords:}} #1
}

\newcommand{\reviewerA}[1]{{#1}}
\newcommand{\reviewerB}[1]{{#1}}
\newcommand{\reviewerAB}[1]{{#1}}

\theoremstyle{definition}

\theoremstyle{remark}
\newtheorem*{remark}{Remark}

\begin{document}

\maketitle

\begin{abstract}
\label{abstact}

The aim of this work is to present a model reduction technique in the framework of optimal control problems for partial differential equations. We combine two approaches used for reducing the computational cost of the mathematical numerical models: domain--decomposition (DD) methods and reduced--order modelling (ROM). In particular, we consider an optimisation--based domain--decomposition algorithm for the parameter--dependent stationary incompressible Navier--Stokes equations. Firstly, the problem is described on the subdomains coupled at the interface and solved through an optimal control problem, which leads to the complete separation of the subdomain problems in the DD method. On top of that, a reduced model for the obtained optimal--control problem is built; the  procedure is based on the Proper Orthogonal Decomposition technique and a further Galerkin projection. The presented methodology is tested on two fluid dynamics benchmarks: the stationary backward--facing step and lid-driven cavity flow. The numerical tests show a significant reduction of the computational costs in terms of both the problem dimensions and the number of optimisation iterations in the domain--decomposition algorithm.  
\end{abstract}


\section{Introduction}
\label{introduction}

In the last decades, there has been a growing interest in approximation techniques for partial differential equations (PDEs) that exploit high--performance computing within different fields of applications: industrial applications, naval engineering, aeronautics engineering, medical engineering, etc. Very often these problems have prohibitively high computational costs, and there is always the need of much more effective algorithms in order to alleviate the complexities of numerical models. 

Two of the most investigated and most important topics for rendering low computational costs are the reduced--order modelling for parameter--dependent PDEs \cite{Rozza_book} and domain--decomposition algorithms \cite{QuarteroniValiDD}. In the former case, equations of interest usually depend on a given set of parameters; these parameters can describe either the physical properties of the sought quantities or the geometrical configuration of the physical domain over which the problem is posed. Model--order reduction is a technique based on the effective decoupling of the computationally expensive offline and usually computationally cheap online phase which provides a solution for any parameter value: for details we refer to \cite{Rozza_book}. Model order reduction has been successfully employed in different fields such as fluid dynamics \cite{ALI20202399,carere2021weighted,crisovan2019model, DeparisRozza2009, LassilaManzoniQuarteroniRozza2014, Rozza2009ReducedBM,stabile2019reduced,stabile2018finite,strazzullo2020pod,strazzullo2022consistency,tezzele2020enhancing,torlo2021model} and structural mechanics \cite{ballarin2016pod,ballarin2017reduced, Haasdonk, RozzaHuynhPatera2007,torlo2018stabilized,venturi2019weighted}. Among the aforementioned applications a significant type of problems often emerges, namely saddle-point problems \cite{benzi_golub_liesen_2005, GernerVeroy2012}, for which special care has to be taken in order to construct stable pairs of the reduced spaces; in particular, in fluid dynamics problems this is achieved by introducing so-called velocity supremisers, see, for instance, \cite{BallarinManzoniQuarteroniRozza2015, GernerVeroy2012,LassilaManzoniQuarteroniRozza2014,tsiolakis2020nonintrusive}.

Another very efficient way for reducing the computational complexity of numerical modes is Domain Decomposition (DD) method. Any domain decomposition method is based on the assumption that a given physical domain of interest is partitioned into subdomains; the original problem is then recast upon each subdomain yielding a family of subproblems of reduced size that are coupled to one another through the values and fluxes of the unknown solution at the subdomain interfaces \cite{QuarteroniValiNumerics, QuarteroniValiDD}. Very often the interface coupling is relaxed at the expense of providing an iterative process among subdomains, allowing a split of each of the subdomain solvers and making it computationally feasible. Domain decomposition methods can be extremely advantageous in the case of very complex geometries as well as in the case of multi-physics problems. 
The latter is even more attractive if we consider that there are often available state-of-the-art codes for a subcomponent model of a multi-physics problem which can be effectively exploited by decoupling algorithms; see, for instance, \cite{ErvinJenkinsLee2014, GosseletChiaruttiniReyFeyel2012, HoangLee2021, lagnese2004domain}.

In this paper, we bring our attention to domain--decomposition methods using an optimisation approach to ensure the coupling of the interface conditions between subdomains as it is presented, for example, in \cite{GunzburgerLee2000, GUNZBURGER2000177}. In particular, we exploit both aforementioned techniques: optimisation--based domain decomposition algorithm in combination with projection--based reduced--order models. This paper is the first step towards the development of an efficient reduced--order model for an optimisation--based domain--decomposition algorithm for Fluid--Structure Interaction (FSI) problems \cite{ballarin2016pod}. It is even more attractive in the view of the articles \cite{KUBERRY2013594, Kuberry2015} where the authors are suggesting that this approach leads to a stable segregated model for FSI problems in the case of added-mass effect \cite{CAUSIN20054506}; we also mention here some already successful ROM results in developing stable semi-implicit partitioned approaches, e.g., \cite{AstorinoChoulyFernandez2010, Nonino2020, Nonino2021}. 

Very recently, authors of the paper \cite{Tezaur2022} have introduced a novel partitioned approach for ROMs, where they couple either two different reduced--order models on each subdomain or a reduced--order model on one subdomain and a full--order (Finite Element) model on the other for the case of nonstationary diffusion--advection problems. In this context, the construction followed in this paper could be also applicable to the coupling presented in \cite{Tezaur2022}, as long as there is a way of casting functions defined on the subdomain interface onto the approximation spaces used on the corresponding subdomains; this will be subject of future work. 

As mentioned before, the use of optimisation--based domain--decomposition methods for PDEs goes as back as the end of the 1990s, e.g., \cite{GUNZBURGER199977}. It had been successfully studied in the case of Navier--Stokes equations as well, see \cite{GunzburgerLee2000}. As for the novelty of this work, to the best of the authors’ knowledge, this is the first attempt of combining the aforementioned technique with projection–based Reduced Order Models in order to provide a computationally efficient algorithm for parametrised PDEs.
Other works deal with model order reduction and domain decomposition but basing their work on other algorithms, e.g. on Schwarz domain overlapping methods \cite{del2023boundary, iollo2023one}. 



A possible extension of current work could be the application of the technique described in this paper to optimal--control problems; in this case, as, for example, in \cite{GUNZBURGER199977}, we have to deal with multi--objective optimisation problems - one for the optimal control and another one for the domain--decomposition part. 

This work is outlined as follows. In Section \ref{problem_formulation} we introduce the monolithic and the optimisation--based domain-decomposition formulation of the incompressible Navier--Stokes equations in both strong and weak forms. Furthermore, we derive the optimality condition for the resulting optimal control problem and compute the expression for the gradient of the objective functional. Section \ref{gradiend_optimisation} lists a gradient--based optimisation algorithm for the problem derived in the previous section. In section \ref{high_fidelity} we describe the Finite Element discretisation of the problem of interest and provide a finite--dimensional high--fidelity optimisation problem. Section \ref{ROM} deals with the reduced--order model which is based on a reduced basis generation by Proper Orthogonal Decomposition methodology and the Galerkin projection of the high--fidelity problem onto the lower-dimensional reduced spaces. In Section \ref{results} we show some numerical results for two toy problems: the backward--facing step and the lid--driven cavity flows. Conclusions will follow in Section \ref{conclusions}.

\section{Problem formulation}
\label{problem_formulation}

In this section, starting with a monolithic formulation of the incompressible Navier-Stokes equations we introduce a two--domain optimisation-based domain-decomposition formulation in both strong and weak forms. Then, the optimality conditions of the resulting optimal control problem are derived followed by the expression of the gradient of the objective functional obtained by sensitivity analysis.  

\subsection{Monolithic formulation}
\label{monolithic_formulation}

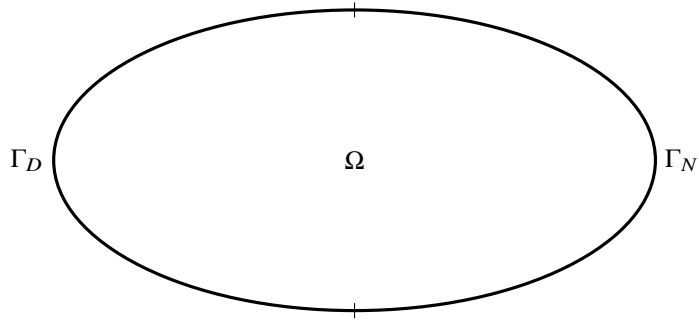
\begin{figure}
    \centering
    \begin{tikzpicture}
 \draw[very thick] (0,0) ellipse (4cm and 2cm);
\node at (0,0) {$\Omega$};
\draw (0,1.9) -- (0,2.1);
\draw (0,-1.9) -- (0,-2.1);
\node[anchor=east] at (-4,0) {$\Gamma_{D}$};
\node[anchor=west] at (4,0) {$\Gamma_{N}$};
\end{tikzpicture}
    \caption{Physical domain}
    \label{fig:mono_domain}
\end{figure}

Let $\Omega$ be a physical domain of interest: we assume $\Omega$ to be an open subset of $\mathbb{R}^2$ and $\Gamma$ to be the boundary of $\Omega$.
Let $f: \Omega \rightarrow \mathbb{R}^2$ be the forcing term, $\nu$ the kinematic viscosity, $u_{D}$ a given Dirichlet datum.  The problem reads as follows: find the velocity field $u: \Omega \rightarrow \mathbb{R}^2 $ and the pressure $p: \Omega \rightarrow \mathbb{R}$ s.t.

\begin{subequations}   \label{eq:mono_equation}
\begin{eqnarray}
    -\nu \Delta u + \left( u \cdot \nabla \right) u + \nabla p = f & \text{in} & \Omega, \label{eq:mono1}\\
    -\text{div} u = 0& \text{in} & \Omega, \label{eq:mono2}  \\
    u = u_{D}  & \text{on} & \Gamma_{D}, \label{eq:mono3} \\
   \nu \frac{\partial u}{\partial n} - p n = 0 & \text{on} & \Gamma_{N} \label{eq:mono4},
\end{eqnarray}
\end{subequations}
where $\Gamma_{D}$ and $\Gamma_{N}$ are  disjoint subsets of $\Gamma$ (as it is shown in Figure \ref{fig:mono_domain}) and $n$ is an outward unit normal vector to $\Gamma_{N}$.

\subsection{Domain Decomposition (DD) formulation} 
\label{DD}
\reviewerA{For the sake of simplicity of exposition, we restricted ourselves to the two--domain decomposition method, but the multi--domain splitting case is a straightforward extension of the two--domain case and should also bring more computational efficiency.} 

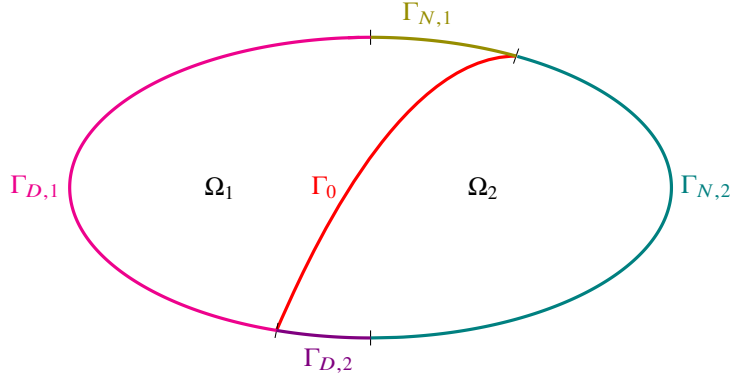
\begin{figure}
    \centering
    \begin{tikzpicture}

\draw[red, very thick]  (1.9365,1.75) parabola (-1.25,-1.9);
\node[red, very thick, anchor=east] at (-0.3, 0) {\textbf{$\Gamma_0$}};

\node at (-2,0) {$\Omega_1$};
\node at (1.5,0) {$\Omega_2$};
\draw (0,1.9) -- (0,2.1);
\draw (0,-1.9) -- (0,-2.1);
\draw (1.9,1.65) -- (1.97,1.85);
\draw (-1.27,-2) -- (-1.23,-1.8);

\node[anchor=east, color=magenta] at (-4,0) {$\Gamma_{D,1}$};
\node[anchor=west, color=teal] at (4,0) {$\Gamma_{N,2}$};
\node[anchor=west, color=olive] at (0.3,2.3) {$\Gamma_{N,1}$};
\node[anchor=east, color=violet] at (-0.1,-2.3) {$\Gamma_{D,2}$};

\draw[domain=0:1.9365, smooth, very thick, variable=\x, samples=51, olive] plot ({\x}, {2*sqrt(1.0-\x*\x/16.0)});

\draw[domain=-2:1.75, smooth, very thick, variable=\y, samples=501, teal] plot ({sqrt(16-4*\y*\y)}, {\y});

\draw[domain=-1.25:0, smooth, very thick, variable=\x, samples=51, violet] plot ({\x}, {-2*sqrt(1.0-\x*\x/16.0)});

\draw[domain=-1.9:2, smooth, very thick, variable=\y, samples=501, magenta] plot ({-4*sqrt(1-\y*\y/4)}, {\y});
\draw[domain=-0.3:0, smooth, very thick, variable=\x, samples=51, magenta] plot ({\x}, {2*sqrt(1.0-\x*\x/16.0)});

\end{tikzpicture}
    \caption{Domain Decomposition of the fluid domain}
    \label{fig:dd_domain}
\end{figure}

Let $\Omega_i, \ i=1,2$ be open subsets of $\Omega$, such that  $\overline{\Omega} = \overline{\Omega_1 \cup \Omega_2}$, \  $\Omega_1 \cap \Omega_2 = \emptyset$. Denote $\Gamma_i := \partial \Omega_i \cap \Gamma, \ i=1,2$ and $\Gamma_0 := \overline{\Omega_1} \cap \overline{\Omega_2}$. In the same way we define the corresponding boundary subsets $\Gamma_{i,D}$ and $\Gamma_{i, N}$, $i=1,2$; see Figure \ref{fig:dd_domain}. \\
Then the DD formulation reads as follows: for $i=1,2$, given $f_i: \Omega_i \rightarrow \mathbb{R}^2$ and $u_{i,D}:\Gamma_{i,D} \rightarrow \mathbb{R}^2$, find $u_i: \Omega_i \rightarrow \mathbb{R}^2$, $p_i: \Omega_i \rightarrow \mathbb{R}$  s.t.

\begin{subequations}\label{eq:dd}
\begin{eqnarray}
    -\nu \Delta u_i + \left( u_i \cdot \nabla \right) u_i + \nabla p_i = f_i & \text{in} & \Omega_i, \label{eq:dd1} \\
    -\text{div} u_i = 0& \text{in} & \Omega_i,  \label{eq:dd2} \\
    u_i = u_{i, D}  & \text{on} & \Gamma_{i, D},  \label{eq:dd3}\\
       \nu \frac{\partial u_i}{\partial n_i} - p_i n_i = 0 & \text{on} & \Gamma_{i, N},   \label{eq:dd4}\\
       \nu \frac{\partial u_i}{\partial n_i} - p_i n_i = (-1)^{i+1}g & \text{on} & \Gamma_{0}  \label{eq:dd5},
\end{eqnarray}
\end{subequations}
for some $g: \Gamma_0 \rightarrow \mathbb{R}^2$ \reviewerB{ such that the functions defined in the following way
\begin{eqnarray*}
     u :=  \begin{cases}
          u_1, \quad \text{   in } \Omega_1 \cup \Gamma_0,  \\
          u_2, \quad  \text{   in } \Omega_2 \cup \Gamma_0,
     \end{cases}
         p :=  \begin{cases}
          p_1, \quad \text{   in } \Omega_1 \cup \Gamma_0,  \\
          p_2, \quad  \text{   in } \Omega_2 \cup \Gamma_0,
     \end{cases}
 \end{eqnarray*}
satisfy the monolithic equations \eqref{eq:mono_equation}.}

Even though in the numerical simulations we will focus on the cases where $f_i = f|_{\Omega_i}$, $u_{i, D} = u_{D}|_{\Gamma_{i, D}}$ for $i=1,2$, the whole theoretical exposition in this paper works just as well for more general functions $f_1, f_2, u_{1, D}$ and $u_{2, D}$.

For any $g$ the solution to the problem \eqref{eq:dd} is not the same as the solution to the problem \eqref{eq:mono_equation}, that is $u_1 \neq u|_{\Omega_1}$,  $p_1 \neq p|_{\Omega_1}$ , $u_2 \neq u|_{\Omega_2}$ and $p_2 \neq p|_{\Omega_2}$. On the other hand, there exists a choice for $g$,  $g = \left(    \nu \frac{\partial u_1}{\partial n_1} - p_1 n_1 \right)|_{\Gamma_0}  = - \left(    \nu \frac{\partial u_2}{\partial n_2}- p_2 n_2 \right)|_{\Gamma_0} $, such that the solutions to  \eqref{eq:dd} coincide with the solution to \eqref{eq:mono_equation} on the corresponding subdomains. Therefore, we must find such a $g$, so that $u_1$ is as close as possible to $u_2$ at the interface $\Gamma_0$. One way to accomplish this is to minimise the functional 
\begin{equation}
    \mathcal J(u_1, u_2) =: \frac{1}{2} \int_{\Gamma_0} \left| u_1 - u_2\right|^2 d\Gamma. \label{eq:unregularised_func}
\end{equation}
Instead of \eqref{eq:unregularised_func} we can also consider the penalised or regularised functional 
\begin{equation}
  \label{eq:functional}  \mathcal J_\gamma(u_1, u_2; g) =: \frac{1}{2} \int_{\Gamma_0} \left| u_1 - u_2\right|^2 d\Gamma + \frac{\gamma}{2}\int_{\Gamma_0} \left|g\right|^2 d\Gamma,
\end{equation}
where $\gamma$ is a constant that can be chosen to change the relative importance of the terms in \eqref{eq:functional}.
Thus we face an optimisation problem under PDE constraints: minimise the functional  \eqref{eq:unregularised_func} (or \eqref{eq:functional}) over a suitable function $g$, subject to \eqref{eq:dd}.

\subsection{Variational Formulation of the PDE constraints}
\label{variational_formulation}

For $i=1,2$ define the following spaces \reviewerB{and the norms with which each of them is endowed}:
\begin{itemize}
    \item $V_i := \left\{ u \in H^1(\Omega_i; \mathbb{R}^2)\right\}$,  \quad \reviewerB{$|| \cdot||_{V_i} = || \cdot ||_{H^1(\Omega_i)}$},
    \item $V_{i,0} := \left\{ u \in H^1(\Omega_i; \mathbb{R}^2): \ u|_{\Gamma_{i, D}  } =0\right\}$, \quad \reviewerB{$|| \cdot||_{V_{i,0}} = || \cdot ||_{H^1_0(\Omega_i)}$},
    \item $Q_i := \left\{ p \in L^2(\Omega_i; \mathbb{R}) \right\}$, \quad \reviewerB{$|| \cdot||_{Q_i} = || \cdot ||_{L^2(\Omega_i)}$}.
\end{itemize} 

Then, we define the following bilinear and trilinear forms: for i=1,2
\begin{itemize}
    \item  $
        a_i: V_i \times V_{i,0} \rightarrow \mathbb{R}, \quad a_i(u_i,v_i) = \nu (\nabla u_i, \nabla v_i)_{\Omega_i}    
    $,
     \item  $
        b_i: V_i \times Q_{i} \rightarrow \mathbb{R}, \quad b_i(v_i,q_i) = - (\text{div} v_i, q_i)_{\Omega_i}    
    $,
     \item  $
        c_i: V_i \times V_i \times V_{i,0} \rightarrow \mathbb{R}, \quad c_i(u_i,w_i, v_i) = \left( (u_i \cdot \nabla)w_i, v_i\right)_{\Omega_i}   
    $,
     
\end{itemize}
where $(\cdot, \cdot)_{\omega}$ indicates the $L^2(\omega)$ inner product. 

Consequently, the variational counterpart of \eqref{eq:dd} reads as follows: for $i=1,2$, find $u_i \in V_i$ and $p_i \in Q_i$ s.t.
\begin{subequations}\label{eq:state}
\begin{align}
 \begin{split}
      a_i(u_i, v_i) + c_i(u_i, u_i, v_i) + b_i(v_i, p_i)    & =   (f_i, v_i)_{\Omega_i}
    \\ \label{eq:state1} & \quad  + \left( (-1)^{i+1} g, v_i \right)_{\Gamma_0} 
 \end{split}  &    \forall v_i \in V_{i,0},      \\
  \label{eq:state2}  b_i(u_i, q_i)  &  = 0 &  \forall q _i \in Q_{i},  \\ \label{eq:state3}
  u_i  &  =   u_{i, D}  &  \text{on}  \ \Gamma_{i, D}.
\end{align}
\end{subequations}

\reviewerB{
\begin{remark}
 In general, the fluxes through an interface $\Gamma_0$ for the weak formulation of Navier--Stokes equation lives in the space $H^{-\frac{1}{2}}(\Gamma_0)$ so that, in theory, the definition \eqref{eq:functional} of functional $\mathcal J_\gamma$ is not justified as it includes the $L^2(\Gamma_0)$--norm of the function $g$. Although, as it will be evident in Section \ref{gradiend_optimisation}, the family of optimisation algorithms which are used to tackle the optimal--control problem in hand, in fact, define the respective approximation of $g$ that belongs to the space $H^{\frac{1}{2}}(\Gamma_0)$.
\end{remark}

}

\subsection{Optimality system}
\label{optimality_system}
One of the ways to address the constrained optimisation problem is to reformulate the initial problem in terms of a Lagrangian functional by introducing the so--called adjoint variables. In this way, the optimal solution to the original problem is sought among the stationary points of the Lagrangian, see, for instance, \cite{Gunzburger_book, Hinze2009}.

We define the Lagrangian functional as follows:

\begin{eqnarray}
    \label{eq:lagrangian} & & \mathcal{L}(u_1, p_1, u_2, p_2, \xi_1, \xi_2, \lambda_1, \lambda_2; g)  :=  \mathcal J_\gamma(u_1, u_2; g)   
    - \sum\limits_{i=1}^2 \left[ a_i(u_i, \xi_i) \quad   \right. \\ \nonumber & &  \quad +c_i(u_i, u_i, \xi_i)    \left.
 + b_i(\xi_i, p_i)  + b_i( u_i, \lambda_i) \right]  + \sum\limits_{i=1}^2(f_i, \xi_i)_{\Omega_i} + \sum\limits_{i=1}^2 ((-1)^{i+1}g, \xi_i)_{\Gamma_0}.     
\end{eqnarray}
Notice that technically we should have also included Lagrange multipliers corresponding to the non--homogeneous Dirichlet boundary conditions \eqref{eq:state3} in the definition of the functional $\mathcal L$, but since the functional $\mathcal J_\gamma$ does not explicitly depend on $u_{1,D}$ and $u_{2, D}$ the corresponding Dirichlet boundary conditions for the adjoint equation that we are going to derive below will be homogeneous on these parts of the boundaries.

We now apply the necessary conditions for finding stationary points of $\mathcal L$. Setting to zero the first variations  w.r.t. $\xi_i$ and $\lambda_i$, $i=1,2$ yields the state equations \eqref{eq:state1}-\eqref{eq:state2}. Setting to zero the first variations  w.r.t. $u_1$, $p_1$, $u_2$ and $p_2$ yields the adjoint equations:
\begin{subequations}\label{eq:adjoint}
\begin{align}
\begin{split}
a_i(\eta_i, \xi_i )+c_i \left(  \eta_i , u_i, \xi_i \right) & + c_i\left( u_i , \eta_i, \xi_i \right)  
+  b_i ( \eta_i, \lambda_i)  \\
\label{eq:adjoint1}&  =  ((-1)^{i+1}\eta_i, u_1 - u_2)_{\Gamma_0},
\end{split}
 & \forall \eta_i \in V_{i,0}, \\
 \label{eq:adjoint2}  b_i (\xi_i, \mu_i) & =0  , & \forall \mu_i\in Q_i.
\end{align}
\end{subequations}

Finally, setting to zero the first variations  w.r.t. $g$ yields the optimality condition:

\begin{equation}
\label{eq:optimality_condition}    \gamma(h, g)_{\Gamma_0} + (h, \xi_1 - \xi_2)_{\Gamma_0} = 0 , \quad \forall h \in L^2(\Gamma_0).
\end{equation}

\subsection{Sensitivity derivatives}
\label{sensivity_derivatives}

In order to obtain the expression for the gradient of the optimisation problem in hand, we will resort to the sensitivity approach, see, for instance, \cite{Gunzburger_book, Hinze2009}. The approach consists of finding equations for direction derivatives of the state variable with respect to control, called sensitivities.

The first derivative $\frac{d\mathcal J_\gamma}{d g}$ of $\mathcal J_{\gamma}$ is defined through its action on variation $\tilde g$ as follows:
\begin{equation}
\label{eq:sensitivity_derivative}    \left<\frac{d  \mathcal{J}_\gamma}{d g}, \tilde g \right> = (u_1 - u_2, \tilde u_1 - \tilde u_2)_{\Gamma_0} + \gamma(g, \tilde g)_{\Gamma_0}, 
\end{equation}
where $\tilde u_1 \in V_{1,0}$, $\tilde u_2 \in V_{2,0}$ are the solutions to:

\begin{subequations}\label{eq:sensitivity_eq}
\begin{align}
\begin{split}
    a_i(\tilde u_i,  v_i) + c_i(\tilde u_i, u_i, v_i)   + c_i(u_i, \tilde u_i, v_i) &
    \\  \label{eq:sensitivity_eq1}    \ + b_i(v_i, \tilde p_i) &   =  ((-1)^{i+1}\tilde g, v_i )_{\Gamma_0}  
\end{split}    
     &  \forall v_i \in V_{i,0},  \\ \label{eq:sensitivity_eq2}
  b_i(\tilde u_i, q_i) & =  0  &\forall q  _i \in Q_{i}. 
\end{align}
\end{subequations}

We can make use of the adjoint equations \eqref{eq:adjoint} in order to find the representation of the gradient of the functional $\mathcal J_\gamma$. Let $\xi_1$ and $\xi_2$ be the solutions to \eqref{eq:adjoint}, $\tilde u_1$ and $\tilde u_2$ be the solutions to \eqref{eq:sensitivity_eq}. By setting $\eta_i = \tilde u_i$ in \eqref{eq:adjoint1}, $\mu_i = \tilde p_i$ in (\ref{eq:adjoint2}), $v_i = \xi_i$ in \eqref{eq:sensitivity_eq1} and $q_i = \lambda_i$ in \eqref{eq:sensitivity_eq2} we obtain:
\begin{equation}
 \nonumber   (u_1-u_2, \tilde u_1 - \tilde u_2)_{\Gamma_0}  = (\tilde g, \xi_1 - \xi_2)_{\Gamma_0},
\end{equation}
so that it yields the explicit formula for the gradient of $\mathcal J_\gamma$:
\begin{equation}
  \label{eq:gradient}  \frac{d\mathcal{J}_\gamma}{dg}(u_1, u_2; g) = \gamma g + (\xi_1 - \xi_2)|_{\Gamma_0},
\end{equation}
where $\xi_1$ and $\xi_2$ are determined from $g$ through \eqref{eq:adjoint}.
Notice that the gradient expression \eqref{eq:gradient} is consistent with the optimality condition \eqref{eq:optimality_condition} derived in the previous section.

\section{Gradient--based algorithm for PDE--constraint optimisation problem}
\label{gradiend_optimisation}

In view of being able to provide a closed--form formula for the gradient for the objective functional $\mathcal J_\gamma$, the natural way to proceed is to resort to a gradient--based iterative optimisation algorithm.

\reviewerB{In order to keep the exposition simple, we consider the following simple gradient method with a constant step size $\alpha>0$}: given a starting guess $g^{(0)}$, let
    \begin{equation}
   \label{eq:step_update}     g^{(n+1)} = g^{(n)} - \alpha \frac{d \mathcal J_\gamma}{dg}\left(u_1^{(n)}, u_2^{(n)}; g^{(n)}\right).
    \end{equation}
    Combining this with \eqref{eq:gradient} we obtain:
        \begin{equation}
     \label{eq:step_update_explicit1}   g^{(n+1)} = g^{(n)} - \alpha \left( \gamma g^{(n)} + (\xi_1^{(n)} - \xi_2^{(n)})|_{\Gamma_0} \right),
    \end{equation}
    or
    \begin{equation}
     \label{eq:step_update_explicit2}   g^{(n+1)} = \left( 1 - \alpha\gamma\right) g^{(n)} - \alpha  (\xi_1^{(n)} - \xi_2^{(n)})|_{\Gamma_0},
    \end{equation}
    where $\xi_1^{(n)}$ and $\xi_2^{(n)}$ are determined from \eqref{eq:adjoint} with $g$ replaced by $g^{(n)}$.
    \\
    In summary, we have the following algorithm:
\textbf{Algorithm 1.}
\begin{enumerate} 
\item Choose $g^{(0)}$, \reviewerB{ $\alpha>0$}.
    \item For n=0,1,2,... until convergence
    \begin{enumerate}
        \item Determine $u_1^{(n)} \in V_1$, $u_2^{(n)} \in V_2$ by solving \eqref{eq:state1}--\eqref{eq:state2} with $g = g^{(n)}$.
        \item Determine $\xi_1^{(n)} \in V_{1,0}$, $\xi_2^{(n)} \in V_{2,0}$ by solving \eqref{eq:adjoint} with
        $u_1 = u_1^{(n)}$, $u_2 = u_2^{(n)}$.
        \item Update $g^{(n+1)}$ by setting
        \begin{equation*}
            g^{(n+1)} := \left(1-\alpha \gamma\right) g^{(n)} -\alpha\left( \xi_1^{(n)} - \xi_2^{(n)}\right)|_{\Gamma_0}.
        \end{equation*}
    \end{enumerate}
    
\end{enumerate}
 
In practice, the typical methods used to solve problems like the one considered in this paper are Broyden–Fletcher–Goldfarb–Shanno (BFGS) and Newton Conjugate Gradient (CG) algorithms which tend to show much faster convergence and higher efficiency with respect to the steepest-decent algorithm.

\section{Finite Element Discretisation}
\label{high_fidelity}
In this section, we present the Finite Element spatial discretisation for the optimal control problem previously introduced. \reviewerB{In order to be able to apply FE discretisation, the domains $\Omega_i,i=1,2$ and the interface $\Gamma_0$ are assumed to be polygonal}. We consider two well-defined triangulations $\mathcal T_1$ and $\mathcal T_2$ over the domains $\Omega_1$ and $\Omega_2$ respectively, and an extra lower-dimensional triangulation $\mathcal T_0$ of the interface $\Gamma_0$; additionally, we assume that $\mathcal T_1$, $\mathcal T_2$ and $\mathcal T_0$ share the same degrees on freedom relative to the interface $\Gamma_0$. We can then define usual Lagrangian FE spaces $V_{i,h} \subset V_{i}$, $V_{i,0, h} \subset V_{i,0}$, $Q_{i,h} \subset Q_i$, $i=1,2$ and $X_h \subset L^2(\Gamma_0)$ endowed with $L^2(\Gamma_0)$-norm; \reviewerB{the spaces $V_{i,h}$, $V_{i,0,h}$ and $Q_{i,h}$ for $i=1,2$ are endowed the same norms as their continuous counterparts}. Since the problems at hand have a saddle-point structure, in order to guarantee the well-posedness of the discretised problem, we require the FE spaces to satisfy the following inf-sup conditions: there exist positive constants $c_1,  c_2, c_3$ and $c_4$ s.t. 

\begin{equation}
\label{eq:infsup}    \inf\limits_{q_{i,h} \in Q_{i,h} \backslash \{0\}} \sup\limits_{v_{i,h} \in V_{i,h} \backslash \{0\}} \frac{b_i(v_{i,h}, q_{i,h})}{||v_{i,h} ||_{V_{i,h}} ||q_{i,h} ||_{Q_{i,h}} } \geq c_i, \quad i=1,2,
\end{equation}
\begin{equation}
  \label{eq:infsup0}  \inf\limits_{q_{i,h} \in Q_{i,h} \backslash \{0\}} \sup\limits_{v_{i,h} \in V_{i,0,h} \backslash \{0\}} \frac{b_i(v_{i,h}, q_{i,h})}{||v_{i,h} ||_{V_{i,0,h}} ||q_{i,h} ||_{Q_{i,h}} } \geq c_{i+2}, \quad i=1,2 .
\end{equation}
A very common choice in this framework is to use the so-called Taylor--Hood finite element spaces, namely the Lagrange polynomial approximation of the second-order for velocity and of the first-order for pressure. We point out that the order of the polynomial space $X_h$ will not lead to big computational efforts as it is defined on the 1-dimensional curve $\Gamma_0$.

Using the Galerkin projection we can derive the following discretised optimisation problem: 
minimise over $g_h \in X_h$ the functional:
\begin{equation}
\label{eq:functional_fem}   \mathcal J_{\gamma,h} (u_{1,h}, u_{2,h}; g_h) := \frac{1}{2}\int_{\Gamma_0} \left| u_{1,h} - u_{2,h}\right|^2 d\Gamma + \frac{\gamma}{2}\int_{\Gamma_0} \left| g_h \right|^2 d\Gamma
\end{equation}
under the constraints that $u_{i, h} \in V_{i,h}$, $p_{i,h} \in Q_{i,h}$ satisfy the following variational equations for $i=1,2$:
\begin{subequations}\label{eq:state_fem}
\begin{align}
 \begin{split}
     a_i(u_{i,h}, v_{i,h}) + c_i(u_{i,h}, u_{i,h}, v_{i,h}) & + b_i(v_{i,h}, p_{i,h})     \\  \label{eq:state_fem1}   =   (f_i, v_{i,h})_{\Omega_i} 
       & + ((-1)^{i+1}g_h, v_{i,h} )_{\Gamma_0}
 \end{split}      &    \forall v_{i,h} \in V_{i,0,h},    \\ \label{eq:state_fem2}
    b_i(u_{i,h}, q_{i,h})  &  = 0 &  \forall q_{i,h} \in Q_{i,h},  \\ \label{eq:state_fem3}
  u_{i,h}  &  =   u_{i,D,h}  &  \text{on}  \ \Gamma_{i, D},
\end{align}
\end{subequations}
where $u_{i,D,h}$ is the Galerkin projection of $u_{i,D}$ onto the trace--space $V_{i,h}|_{\Gamma_{i,D}}$.

Notice that the structure of the equations \eqref{eq:state_fem} and of the functional \eqref{eq:functional_fem} is the same as the one of the continuous case so that it enables us to provide the following expression of the gradient of the discretised functional \eqref{eq:functional_fem}:
\begin{equation}
 \label{eq:gradient_fem}   \frac{d\mathcal{J}_{\gamma, h}}{dg_h}(u_{1,h}, u_{2,h}; g_h) = \gamma g_h + (\xi_{1,h} - \xi_{2,h})|_{\Gamma_0},
\end{equation}
where $\xi_{1,h}$ and $\xi_{2,h}$ are the solutions to the discretised adjoint problem: for $i=1,2$ find $\xi_{i,h} \in V_{i,0,h}$ and $\lambda_{i,h} \in Q_{i,h}$ that satisfy 
\begin{subequations}\label{eq:adjoint_fem}
\begin{align}
\begin{split}
 a_i(\eta_{i,h}, \xi_{i,h} ) & +  c_i \left(  \eta_{i,h} , u_{i,h}, \xi_{i,h} \right) 
 + c_i\left( u_{i,h} , \eta_{i,h}, \xi_{i,h} \right)  
  \\ \label{eq:adjoint_fem1}  & 
+    b_i ( \eta_{i,h}, \lambda_{i,h})  
 =   ((-1)^{i+1}\eta_{i,h}, u_{1,h} - u_{2,h})_{\Gamma_0},
\end{split}  & \forall \eta_{i,h} \in V_{i,0,h}, \\ \label{eq:adjoint_fem2}
 &  b_i (\xi_{i,h}, \mu_{i,h})   = 0  , & \forall \mu_{i,h}\in Q_{i,h}.
\end{align}
\end{subequations}

We would also like to stress that at the algebraic level the discretised minimisation problem can be recast in the setting of the finite--dimensional space $\mathbb{R}^p$, where $p$ is the number of Finite Element degrees of freedom which belong to the interface $\Gamma_0$.

\section{Reduced-Order Model}
\label{ROM}

As it was highlighted in section \ref{introduction}, Reduced--Order methods are efficient tools for significant reduction of the parameter--dependent PDEs.
This section deals with the reduced--order model for the problem obtained in the previous section, where the state equations, namely Navier--Stokes equations, are assumed to be dependent on a set of physical parameters. First, we introduce two practical ingredients we will be using in the course of the reduced--basis generation, namely a lifting function and the pressure supremiser enrichment. Then, we describe the offline phase based on the Proper Orthogonal Decomposition technique, which is followed by the online phase based on a Galerkin projection onto the reduced spaces.

\subsection{Lifting Function and Velocity Supremiser Enrichment}
\label{lifting_supremiser}

In the following, we are going to discuss a snapshot compression technique for the generation of reduced basis functions. In order to do so we need to introduce two important ingredients in this context, namely the lifting function technique and the supremiser enrichment of the velocity space.

The use of lifting functions is quite common in the reduced basis method (RBM) framework; see, for example, \cite{Rozza_book, BallarinManzoniQuarteroniRozza2015}. It is motivated by the fact that in the chosen model we are supposed to tackle the non-homogeneous Dirichlet boundary condition on the parts of the boundaries $\Gamma_{i,D}, \,i=1,2$. From the implementation point of view, this does not present any problem when dealing with the high-fidelity model since there are several well-known techniques for non-homogeneous essential conditions, in particular at the algebraic level. However, these boundary conditions  create some problems when dealing with the reduced basis methods. Indeed, we seek to generate a linear vector space which is obtained by the compression of the set of snapshots, and this clearly cannot be achieved by using snapshots which satisfy different Dirichlet conditions -- the resulting space would not be linear. This problem is solved by introducing a lifting function $l_{i,h} \in V_{i,h},\, i =1,2$ during the offline stage, such that $l_{i,h} = u_{i, D, h}$ on $\Gamma_{i,D}$. We define two new variables $u_{i,0,h} \in V_{i,0,h }, \,i=1,2$ by setting $u_{i,0,h}:=u_{i,h}-l_{i,h}$. Clearly, the variables $u_{i,0,h}, \, i=1,2$ satisfy the homogeneous condition $u_{i,0,h} = 0$ on $\Gamma_{i,D}$ and so they can be used to generate the reduced basis linear space. We remark that the lifting function is needed only in the domain where the Dirichlet boundary is non--empty, i.e. where $\Gamma_{i,D} \neq \emptyset$ for $i=1,2$. It is important to point out that the choice of lifting functions is not unique; in our work, we chose to use the solution of the Stokes problem in one of the domains $\Omega$, $\Omega_1$ or $\Omega_2$ (depending on the particular model we are investigating) with the velocity equal to $u_{D}$ on the corresponding parts of the boundaries and the homogeneous Neumann conditions analogous to the original problem setting.   

The other ingredient we will use in the following exposition is the so-called velocity supremiser. This is necessary to obtain a stable approximation of the saddle-point problem at the reduced level discussed in the following subsections. The well--posedness of the problem is again assured by satisfying the inf-sup conditions like \eqref{eq:infsup0}. The supremiser variables $s_{i,h}$, $i=1,2$ are defined as the solution to the following problem: find $s_{i,h} \in V_{i,0,h}$ such that 

\begin{equation}
    \left( \nabla v_{i,h},\nabla s_{i,h}\right) = b_{i,h}\left( v_{i,h}, p_{i,h} \right) \quad \forall v_{i,h} \in V_{i,0,h}, \label{eq:supremiser}
\end{equation}
where $p_{i,h}, \ i=1,2$ are the finite-element pressure solutions of the Navier-Stokes problem and the left-hand side is the scalar product which defines a norm on the space $V_{i,0,h}$. For more details, we refer to \cite{BallarinManzoniQuarteroniRozza2015, GernerVeroy2012}. \reviewerB{Another way to apply the supremiser is to apply it directly to the reduced basis of the velocity spaces, but this might lead to parameter--dependent reduced spaces \cite{BallarinManzoniQuarteroniRozza2015}. Other simplifications may work in a similar fashion, we might compare them in future works.} 

\subsection{Reduced Basis Generation}
\label{POD}

Once we obtain the homogenised snapshots $u_{i,0,h}$ and the pressure supremisers $s_{i,h}$ for $i=1,2$, we are ready to construct a set of reduced basis functions. A very common choice when dealing with Navier-Stokes equations is to use the Proper Orthogonal Decomposition (POD) technique, which is based on the Singular Value Decomposition of the snapshot matrices; see, for instance, \cite{Rozza_book}. In order to implement this technique we will need two main ingredients: the matrices of the inner products and the snapshot matrices. First, we define the basis functions for the FE element spaces used in the weak formulation \eqref{eq:functional_fem}, \eqref{eq:state_fem} and \eqref{eq:adjoint_fem} as follows:
\begin{eqnarray*}
\mathcal U_{i,0,h} = \left\{ \phi_1^{u_i}, ..., \phi_{\mathcal N_h^{u_i}}^{u_i} \right\} - \text{the FE basis of the space $V_{i,0,h}, i=1,2$}, \\
\mathcal P_{i,h}=\left\{ \phi_1^{p_i}, ..., \phi_{\mathcal N_h^{p_i}}^{p_i} \right\} - \text{the FE basis of the space $Q_{i,h}, i=1,2$},\\
\reviewerB{\Xi_{i,0,h} :=\mathcal U_{i,0,h}, \quad \mathcal N_h^{\xi_i} := \mathcal N_h^{u_i},  i=1,2},\quad \quad \quad \quad \quad \quad   \\
\mathcal G_{i,h} =\left\{ \phi_1^{g}, ..., \phi_{\mathcal N_h^{g}}^{g} \right\} - \text{the FE basis of the space $X_h$}, \quad \quad  \quad \quad 
\end{eqnarray*}
where $\mathcal N_h^{*}, * \in \left\{ u_1, p_1, u_2, p_2, g\right\}$ denotes the dimension of the corresponding FE space. 

We proceed by building the snapshot matrices. In doing so we sample a parameter space and draw a discrete set of $M$ parameter values; there are various sampling techniques, among which we point out the uniform sampling. Then, the snapshots are taken as a high--fidelity, i.e. Finite Element, solutions at each parameter value in the sampling set. 

We proceed by building the snapshot matrices $\mathcal S_{u_i} \in \mathbb{R}^{\mathcal N_h^s\times 4M}$,  $\mathcal S_{s_i} \in \mathbb{R}^{\mathcal N_h^s\times 4M}$, $\mathcal S_{p_i} \in \mathbb{R}^{\mathcal N_h^s\times 4M}$, $\mathcal S_{\xi_i} \in \mathbb{R}^{\mathcal N_h^a\times 2M}$ for $i=1,2$ and $\mathcal S_{g} \in \mathbb{R}^{\mathcal N_h^g\times M}$ defined as follows: 
\begin{eqnarray*}
\mathcal S_{u_1} & = & [u_{1,0,h}^{1},..., u_{1,0,h}^{M}, 0,...,0 , 0,...,0, 0,...,0], \\
\mathcal S_{s_1} & = & [s_{1,h}^{1},..., s_{1,h}^{M}, 0,...,0 , 0,...,0, 0,...,0], \\
\mathcal S_{p_1} & = & [0,..., 0, p_{1,h}^{1},..., p_{1,h}^{M}, 0,...,0 , 0,...,0], \\
\mathcal S_{u_2} & = & [0,..., 0, 0,..., 0, u_{2,0,h}^{1},..., u_{2,0,h}^{M}, 0,...,0], \\
\mathcal S_{s_2} & = & [0,...,0, 0, ..., 0,s_{2,h}^{1},..., s_{2,h}^{M}, 0,...,0], \\
\mathcal S_{p_2} & = & [0,...,0,0,..., 0,0,..., 0, p_{2,h}^{1},..., p_{2,h}^{M}], 
\\
\mathcal S_{\xi_1}  & = &  [\xi_{1,h}^{1},..., \xi_{1,h}^{M}, 0,..., 0], \quad 
\mathcal S_{\xi_2}  =  [0, ..., 0, \xi_{2,h}^{1},..., \xi_{2,h}^{M}], \\
\mathcal S_{g} & = & [g_{h}^{1},..., g_{h}^{M}],
\end{eqnarray*}
where $\mathcal N_h^{s} =\mathcal N_h^{u_1}+\mathcal N_h^{p_1}+\mathcal N_h^{u_2}+\mathcal N_h^{p_2}$, $\mathcal N_h^{a} = \mathcal   N_h^{\xi_1}+\mathcal N_h^{\xi_2}$ and $M$ is the number of snapshots.

Notice that since all the snapshots of the variables $\xi_{1,h}$ and $\xi_{2,h}$ are divergence-free on the domain of definition, the reduced spaces constructed for those variables will already contain this information, so that it allows us not to store the snapshots of the variables $\lambda_{1,h}$ and $\lambda_{2,h}$, which are playing the role of the Lagrange multipliers relative to the divergence free-conditions, as they do not contain any important information.   

The next step is to define the inner-product matrices $X_{u_i}$, $X_{p_i}$, $X_{\xi_i}$ for $i=1,2$ and $X_g$. These matrices have the block diagonal structure as follows:
\begin{eqnarray*}
X_{u_1} & = & \text{diag}\left(x_{u_1}, 0_{p_1}, 0_{u_2}, 0_{p_2}\right), \\
X_{p_1} & = & \text{diag}\left(0_{u_1}, x_{p_1}, 0_{u_2}, 0_{p_2} \right), \\
X_{u_2} & = & \text{diag}\left(0_{u_1}, 0_{p_1}, x_{u_2}, 0_{p_2}\right), \\
X_{p_2} & = & \text{diag}\left(0_{u_1}, 0_{p_1}, 0_{u_2}, x_{p_2} \right), \\
X_{\xi_1} & = & \text{diag}\left( \reviewerB{x_{u_1}}, 0_{\xi_2}\right), \\
X_{\xi_2} & = & \text{diag}\left( 0_{\xi_1}, \reviewerB{x_{u_2}} \right), \\
X_{g} & = & x_{g}.
\end{eqnarray*}

Above, we used the following notations: $0_{*} \in \mathbb R^{\mathcal N_h^{*}\times \mathcal N_h^{*}}$ is a zero square matrix of dimension $N_h^{*}\times \mathcal N_h^{*}$, where $* \in \left\{ u_1, p_1, u_2, p_2, \xi_1, \xi_2, g \right\}$ and 

\begin{eqnarray*}
 (x_{u_i})_{jk} & = & \left( \nabla \phi_{k}^{u_i}, \nabla \phi_{j}^{u_i} \right)_{\Omega_i}, \quad \text{for } j, k = 1,..., \mathcal{N}_h^{u_i}, \ i=1,2, \\
 (x_{p_i})_{jk} & = & \left(  \phi_{k}^{p_i},  \phi_{j}^{p_i} \right)_{\Omega_i}, \quad \text{for } j, k = 1,..., \mathcal{N}_h^{p_i}, \ i=1,2, \\
 (x_{g})_{jk} & = & \left(  \phi_{k}^{g},  \phi_{j}^{g} \right)_{\Gamma_0}, \quad \text{for } j, k = 1,..., \mathcal{N}_h^{g}.
\end{eqnarray*}

We are now ready to introduce the correlation matrices $\mathcal{C}_{u_i}$, $\mathcal{C}_{s_i}$, $\mathcal{C}_{p_i}$, $\mathcal{C}_{\xi_i}$ for $i=1,2$ and $\mathcal{C}_g$, all of dimension $M\times M$, as:
\begin{eqnarray*}
\mathcal{C}_{*} := \mathcal S_{*}^T X_{*} S_{*}
\end{eqnarray*}
for every $* \in \{u_1, p_1, u_2, p_2, \xi_1, \xi_2, g\}$
and
\begin{eqnarray*}
\mathcal{C}_{s_i} := \mathcal S_{s_i}^T X_{u_i} S_{s_i}, \ i=1,2.
\end{eqnarray*}

Once we have built the correlation matrices, we are able to carry out a POD compression on the sets of snapshots. This can be achieved by solving the following eigenvalue problems:
\begin{eqnarray}
\mathcal{C}_{*}\mathcal{Q}_{*} = \mathcal{Q}_{*} \Lambda_{*}  \label{eq:eigenvalue_problem}
\end{eqnarray}
where $* \in \{u_1, s_1, p_1, u_2, s_2, p_2, \xi_1, \xi_2, g\}$, $\mathcal Q_{*}$ is the eigenvectors matrix and $\Lambda_{*}$ is the diagonal eigenvalues matrix with eigenvalues ordered by decreasing order of their magnitude. The $k$-th reduced basis function for the component ${*}$ is then obtained by applying the matrix $\mathcal S_{*}$ to $\underline v_k^{*}$ -- the $k$-th column vector of the matrix $\mathcal{Q}_{*}$:
\begin{equation*}
    \Phi_k^{*}:=\frac{1}{\sqrt{\lambda_k^{*}}} \mathcal{S}_{*} \underline v_k^{*},  
\end{equation*}
where $\lambda_k^{*}$ is the $k$-th eigenvalue from \eqref{eq:eigenvalue_problem}. 
Therefore, we are able to form the set of reduced basis as
\begin{equation*}
\mathcal A^s := \bigcup\limits_{ * \in \{u_1, s_1, p_1, u_2, s_2, p_2 \} }\left\{ \Psi_1^{*}, ..., \Psi_{N_{*}}^{*}   \right\},
\end{equation*}
\begin{equation*}
\mathcal A^a := \bigcup\limits_{ * \in \{\xi_1, \xi_2\} }\left\{ \Psi_1^{*}, ..., \Psi_{N_{*}}^{*}   \right\},
\end{equation*}
\begin{equation*}
\mathcal A^g := \left\{ \Phi_1^{g}, ..., \Phi_{N_{g}}^{g}   \right\},
\end{equation*}
where the integer numbers $N_{*}$ indicate the number of the basis functions used for each component and 
\begin{equation*}
    \Psi_k^{u_1} = \left( 
    \begin{array}{c}
         \Phi_k^{u_1}  \\
         0 \\
         0 \\
         0 
    \end{array}
    \right), \
    \Psi_k^{s_1} = \left( 
    \begin{array}{c}
         \Phi_k^{s_1}  \\
         0 \\
         0 \\
         0 
    \end{array}
    \right),\
    \Psi_k^{p_1} = \left( 
    \begin{array}{c}
         0 \\
         \Phi_k^{p_1}  \\
         0 \\
         0 
    \end{array}
    \right),\
        \Psi_k^{u_2} = \left( 
    \begin{array}{c}
      0 \\
         0 \\
         \Phi_k^{u_2}  \\
         0 
    \end{array}
    \right), 
\end{equation*}
\begin{equation*}
\Psi_k^{s_2} = \left( 
    \begin{array}{c}
       0 \\
         0 \\
         \Phi_k^{s_2}  \\
         0 
    \end{array}
    \right),\
    \Psi_k^{p_2} = \left( 
    \begin{array}{c}
         0 \\
         0 \\
         0 \\
         \Phi_k^{p_2} 
    \end{array}
    \right),\
    \Psi_k^{\xi_1} = \left( 
    \begin{array}{c}
         \Phi_k^{\xi_1}  \\
         0 
    \end{array}
    \right),\
    \Psi_k^{\xi_2} = \left( 
    \begin{array}{c}
         0 \\
         \Phi_k^{\xi_2}  
    \end{array}
    \right).
\end{equation*}
We note that the first and the third blocks include both the $u_1$, $s_1$ and the $u_2$, $s_2$ basis functions - it is here that we use the pressure supremiser enrichment of the velocities spaces discussed at the beginning of this section. We provide the following renumbering of the functions for further simplicity: 
\begin{equation*}
    \Phi_{N_{u_i}+k}^{u_i}:=\Phi_{k}^{s_i}, \ \Psi_{N_{u_i}+k}^{u_i}:=\Psi_{k}^{s_i}, \quad \text{for} \ k=1,...,N_{s_i}, \ i=1,2,
\end{equation*}
and we redefine $N_{u_i}:=N_{u_i}+N_{s_i}, \ i=1,2$.

Finally, we introduce three separate reduced basis spaces - for the state, the adjoint and the control variables, respectively:
\begin{eqnarray*}
    V_{N}^s = \text{span} \mathcal{A}^s, & & \text{dim} V_N^s = N_{u_1}+N_{p_1}+N_{u_2}+N_{p_2},
\\
V_{N}^a = \text{span} \mathcal{A}^a, & & \text{dim} V_N^s = N_{\xi_1}+N_{\xi_2},
\\
V_N^g = \text{span} \mathcal{A}^g, & & \text{dim} V_N^s = N_{g}.
\end{eqnarray*}

\subsection{Online Phase}
\label{online}

Once we have introduced the reduced basis spaces we can define the reduced function expansions 
\begin{equation*}
U_N = (u_{1,0, N}, p_{1, N}, u_{2,0, N}, p_{2, N}) \in V_{N}^s, \Xi_N = (\xi_{1, N}, \xi_{2, N}) \in V_N^a, g_{ N} \in V_N^g
\end{equation*}
as
\begin{align*}
u_{i, 0, N}   :=  \sum\limits_{k=1}^{N_{u_i}} \underline u_{i,0,k}\Phi_k^{u_i}, \ i=1,2, & &  \xi_{i, N}   :=   \sum\limits_{k=1}^{N_{\xi_i}} \underline \xi_{i,k}\Phi_k^{\xi_i}, \ i=1,2,
\\
p_{i, N}  :=   \sum\limits_{k=1}^{N_{p_i}} \underline p_{i,k}\Phi_k^{p_i}, \ i=1,2, & &
g_{ N}   :=   \sum\limits_{k=1}^{N_{g}} \underline g_{k}\Phi_k^{g} \quad \quad\quad.
\end{align*}

In the previous equations, the underscore indicates the coefficients of the basis expansion of the reduced solution. Then the online reduced problem reads as follows: minimise over $g_N \in V_N^g$ the functional
\begin{equation}
    \mathcal J_{\gamma, N} (u_{1,N}, u_{2,N}; g_N) := \frac{1}{2}\int_{\Gamma_0} \left| u_{1,N} - u_{2,N}\right|^2 d\Gamma + \frac{\gamma}{2}\int_{\Gamma_0} \left| g_N \right|^2 d\Gamma  \label{eq:functional_rom}
\end{equation}
where $u_{1,N} = u_{1,0,N}+l_{1,N}$, $u_{2,N} = u_{2,0,N}+l_{2,N}$ for $(u_{1,0,N}, p_{1,N}, u_{2,0,N}, p_{2,N}) \in V_N^s$  satisfy the following reduced equations $\forall v_N  = ( v_{1,N}, q_{1,N}, v_{2,N}, q_{2,N})  \in V_N^s $:

\begin{subequations}\label{eq:state_rom}
\begin{eqnarray}
  \nonumber   a_i(u_{i,0,N}, v_{i,N}) &  + &   c_i(u_{i,0,N}, u_{i,0,N}, v_{i,N}) 
 +  c_i(u_{i,0,N}, l_{i,N}, v_{i,N}) 
    \\\nonumber 
 &  + &   c_i(l_{i,N}, u_{i,0,N}, v_{i,N}) 
     +   b_i(v_{i,N}, p_{i,N})
        \\ \label{eq:state_rom1} & = &   (f_i, v_{i,N})_{\Omega_i}
      +((-1)^{i+1}g_N, v_{i,N} )_{\Gamma_0}          
       \\\nonumber   &  & - a_i(l_{i,N}, v_{i,N}) 
      -c_i(l_{i,N},l_{i,N},v_{i,N}) 
       \\  \label{eq:state_rom2}
    b_i(u_{i,0,N}, q_{i,N})   & = & -b_i(l_{i,N}, q_{i,N}),
\end{eqnarray}
\end{subequations}
where $l_{i,N}$ is the Galerkin projection of the lifting function $l_{i,h}$ to the finite dimensional vector space spanned by the $i$-th velocity basis functions and $i=1,2$.

Similarly to the offline phase, we notice that the structure of the equations \eqref{eq:state_rom} and the functional \eqref{eq:functional_rom} are the same as the ones of the continuous case, so this enables us to provide the following expression of the gradient of the reduced functional \eqref{eq:functional_rom}:
\begin{equation}
 \label{eq:gradient_rom}   \frac{d\mathcal{J}_{\gamma, N}}{dg_N}(u_{1,N}, u_{2,N}; g_N) = \gamma g_{N} + (\xi_{1,N} - \xi_{2,N})|_{\Gamma_0},
\end{equation}
where $(\xi_{1,N},\xi_{2,N}) \in V_N^a$ are the solutions to the reduced adjoint problem: find $(\xi_{1,N}, \xi_{2,N}) \in V_N^a $ such that it satisfies, for each pair of test functions $(\eta_{1,N}, \eta_{2,N}) \in V_N^a$ and $i=1,2$, 
\begin{align}
\begin{split}
 a_i(\eta_{i,N}, \xi_{i,N} )  +  c_i \left(  \eta_{i,N} , u_{i,N}, \xi_{i,N} \right) &
 + c_i\left( u_{i,N} , \eta_{i,N}, \xi_{i,N} \right) \\  \label{eq:adjoint_rom}   &  =   ((-1)^{i+1}\eta_{i,N}, u_{1,N} - u_{2,N})_{\Gamma_0}. 
\end{split}
\end{align}
Notice that the reduced adjoint equations no longer contain any terms corresponding to the bilinear forms $b_i(\cdot, \cdot), \, i=1,2$. Indeed, as was previously mentioned, all the functions belonging to the reduced space $V_N^a$ are already divergence-free by construction, so the aforementioned terms are automatically satisfied.  

We would also like to stress that from the numerical implementation point of view the reduced minimisation problem can be recast in the setting of the finite-dimensional space $\mathbb{R}^{p}$, where $p$ is the number of reduced basis function used for the control variable $g_N$ in the online phase, that is $p=N_{g}$.

\section{Numerical Results}
\label{results}
We now present some numerical results obtained by applying the two-domain decomposition optimisation algorithm to the backward--facing step and the lid-driven cavity flow benchmarks. 

All the numerical simulations for the offline phase were obtained using the software multiphenics \cite{multiphenics}, whereas the online phase simulations were carried out using RBniCS \cite{rbnics}.

\subsection{Backward-facing step test case}
\label{backward_facing_step}

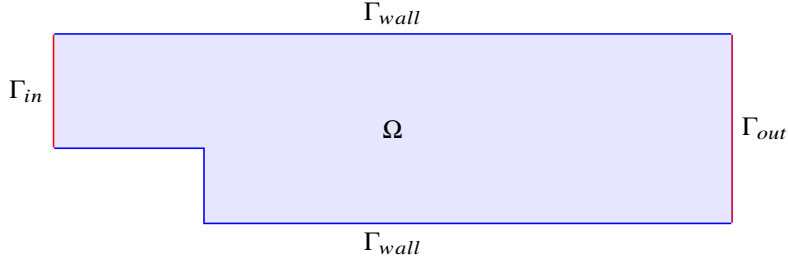
\begin{figure}
    \centering
    
\begin{tikzpicture}

    \draw[draw=red,very thick] (0,1) -- (0, 1.75) node [anchor=east]{$\Gamma_{in}$}  -- (0,2.5);
    \draw[draw=blue,very thick] (0,2.5) --  (4.5, 2.5) node [anchor=south] {$\Gamma_{wall}$}-- (9,2.5);
    \draw[draw=purple,very thick](9,2.5) -- (9, 1.25) node [anchor=west] {$\Gamma_{out}$}-- (9,0);
    \draw[draw=blue,very thick] (9,0) -- (4.5, 0) node [anchor=north]{$\Gamma_{wall}$} -- (2,0)--(2,1)--(0,1);

    \fill[fill=blue!10] (0,1)  -- (0,2.5) -- (9,2.5) -- (9, 1.25) -- (9,0)  -- (2,0)--(2,1)--(0,1);
    
    \draw (4.5, 1.25) node {$\Omega$};
\end{tikzpicture}
\caption{Physical domain for the backward-facing step problem\label{fig:Mono_domain_bfs}
}

\end{figure}
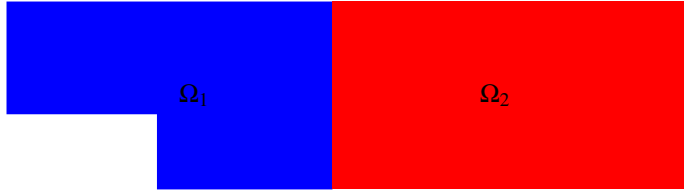
\begin{figure}
    \centering
    
\begin{tikzpicture}

    \fill[fill=blue] (0,1)  -- (0,2.5) -- (4.3333,2.5) -- (4.3333,0)  -- (2,0)--(2,1)--(0,1);
    
    \draw (2.5, 1.25) node {$\Omega_1$};
    
    \fill[fill=red] (4.3333,0) -- (4.3333,2.5) -- (9,2.5) -- (9,0)-- (4.3333,0) ;
    
    \draw (6.5, 1.25) node {$\Omega_2$};
\end{tikzpicture}
\caption{Domain decomposition for the backward-facing step problem domain\label{fig:dd_domain_bfs}
}

\end{figure}

\begin{table}[ht]
\begin{center}
\begin{tabular}{ c c }

\hline
    &   \\
    Physical parameters & $2: \nu, \bar U$ \\
    Range $\nu$ & [0.5, 2] \\
    Range $\bar U$ & [0.5, 6.5] \\
    \reviewerAB{Resulting $Re$ number} & \reviewerAB{[0.75, 40]} \\
    &   \\
    FE velocity order & 2 \\
    FE pressure order & 1 \\
    Total number of FE dofs & 27,890 \\
    Number of FE dofs at the interface & 130 \\
    & \\
    Optimisation algorithm & L-BFGS-B \\
    $It_{max}$ & 40 \\
    $Tol_{opt}$ & $10^{-5}$ \\
    & \\
    $M$ & 900 \\
    $N_{max}$ & 50 \\
    \hline
\end{tabular}
\caption{Computational details of the offline stage. \label{table:offline_bfs}}
\end{center}
\end{table}

We start with introducing the backward--facing step flow test case. Figure~\ref{fig:Mono_domain_bfs} represents the  physical domain of interest. The upper part of the channel has a length of 18 cm, the lower part 14 cm; the height of the left chamber is 3 cm, and the height of the right one is 5 cm. 
The splitting into two domains is performed by dissecting the domain by a vertical segment at the distance $\frac{26}{3}$ cm from the beginning of the channel as shown in Figure \ref{fig:dd_domain_bfs}.

We impose homogeneous Dirichlet boundary conditions on the top and the bottom walls of the boundary $\Gamma_{wall}$ for the fluid velocity, and homogeneous Neumann conditions on the outlet $\Gamma_{out}$, meaning that we assume free outflow on this portion of the boundary. 

We impose a parabolic profile $u_{in}$ on the inlet boundary $\Gamma_{in}$, where
\begin{equation}
    u_{in}(x,y) = 
    \left( 
    \begin{array}{c}
         w(y)  \\
         0
    \end{array}
    \right)
\end{equation}
with $w(y) = \bar U \times \frac{4}{9} (y-2)(5-y), \ y \in [2,5]$; values of $\bar U$ are reported in Table~\ref{table:offline_bfs}.
Two physical parameters are considered: the viscosity $\nu$ and the maximal magnitude $\bar U$ of the inlet velocity profile $u_{in}$. \reviewerB{Both parameters concur to the definition of the only physically relevant parameter, the Reynolds number $Re=L \frac{\bar U}{\nu}$, where $L$ is the characteristic length. Hence, we indicate for all tests also the corresponding $Re$.} Details of the offline stage and the finite-element discretisation are summarised in Table \ref{table:offline_bfs}. High-fidelity solutions are obtained by carrying out the minimisation in the space of dimension equal to the number of degrees of freedom at the interface, which is 130 in our test case. The best performance has been achieved by using the limited-memory Broyden–Fletcher–Goldfarb–Shanno (L-BFGS-B) optimisation algorithm, and two stopping criteria were applied: either the maximal number of iteration $\ It_{max}$ is reached or the gradient norm of the target functional is less than the given tolerance $Tol_{opt}$.

\begin{figure}
    \centering
    \begin{subfigure}[b]{0.49\textwidth}
        \includegraphics[width=\textwidth]{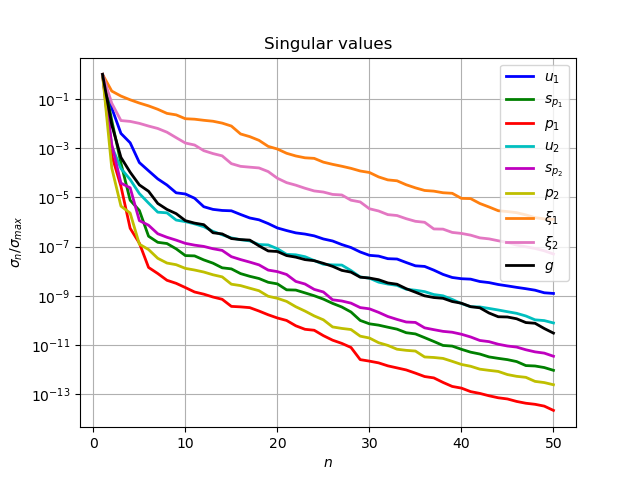}
        \caption{POD singular values as a function of number $n$ of POD modes (log scaling in $y$-direction)}
         \label{fig:singlular_values_bfs}
    \end{subfigure}
    \hfill
    \begin{subfigure}[b]{0.49\textwidth}
        \includegraphics[width=\textwidth]{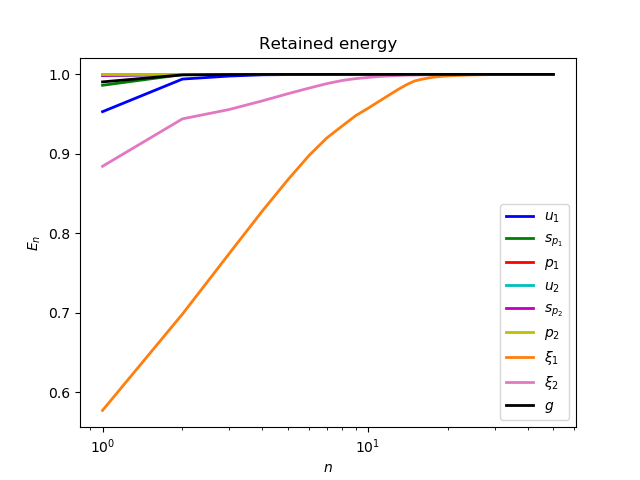}
        \caption{Energy retained by the first $N_{max}$ POD modes (log scaling in $x$-direction)}
         \label{fig:retained_energy_bfs}
    \end{subfigure}      
    \caption{Results of the offline stage: POD singular eigenvalue decay (a) and retained energy (b) of the first $N_{max}$ POD modes}
    \label{fig:pod_modes_bfs}
\end{figure}

Snapshots are sampled from a training set of $M$ parameters uniformly distributed in the 2-dimensional parameter space, and the first $N_{max}$ POD modes have been retained. Figure~\ref{fig:singlular_values_bfs} shows the POD singular values for all the state, the adjoint and the control variables. As it can be seen, the POD singular values corresponding to the adjoint velocities $\xi_1$ and $\xi_2$ feature a slower decay compared to the one for the other variables. In Figure~\ref{fig:retained_energy_bfs}, we can see the behaviour of the energy $E_n$ retained by the first $N$  modes for different components of the solution. Here, the retained energy for the component $* \in \{u_1, s_1, p_1, u_2, s_2, p_2, \xi_1, \xi_2, g\}$ is defined as
\begin{equation*}
    E_n^{*} := \frac{\sum_{k=1}^n|\lambda_k^*|}{\sum_{k=1}^{N_*}|\lambda_k^{*}|}.
\end{equation*}

The retained energy gives us an idea on the number of modes we would need to choose to preserve all the necessary physical information in the reduced model. In particular, we can see that a higher number of modes is needed to correctly represent the adjoint variables $\xi_1$ and $\xi_2$.   

\begin{figure}[H]
    \centering
    \begin{subfigure}[b]{0.49\textwidth}
        \includegraphics[width=\textwidth]{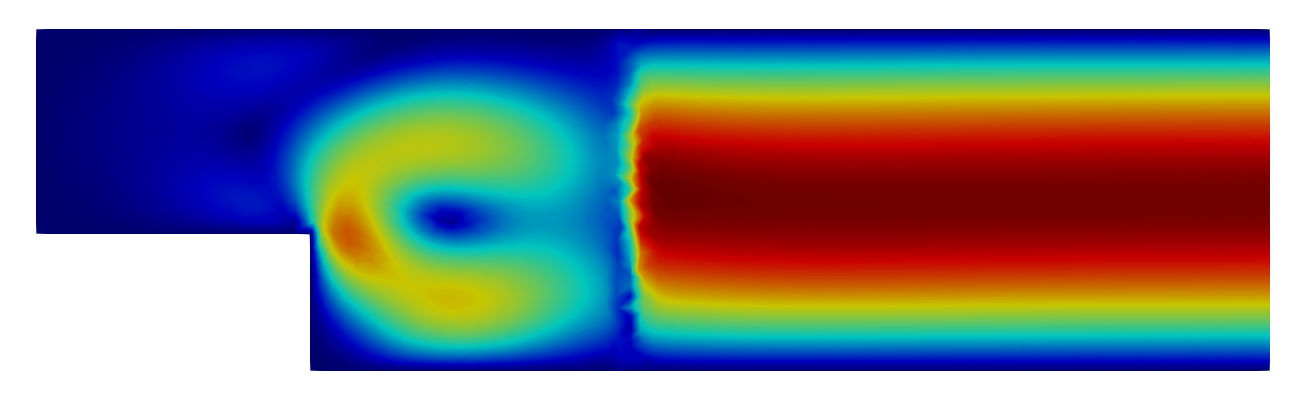}
         \label{fig:mode_u_1_bfs}
         \caption{\reviewerA{The first mode}}
    \end{subfigure}
    \hfill
    \begin{subfigure}[b]{0.49\textwidth}
        \includegraphics[width=\textwidth]{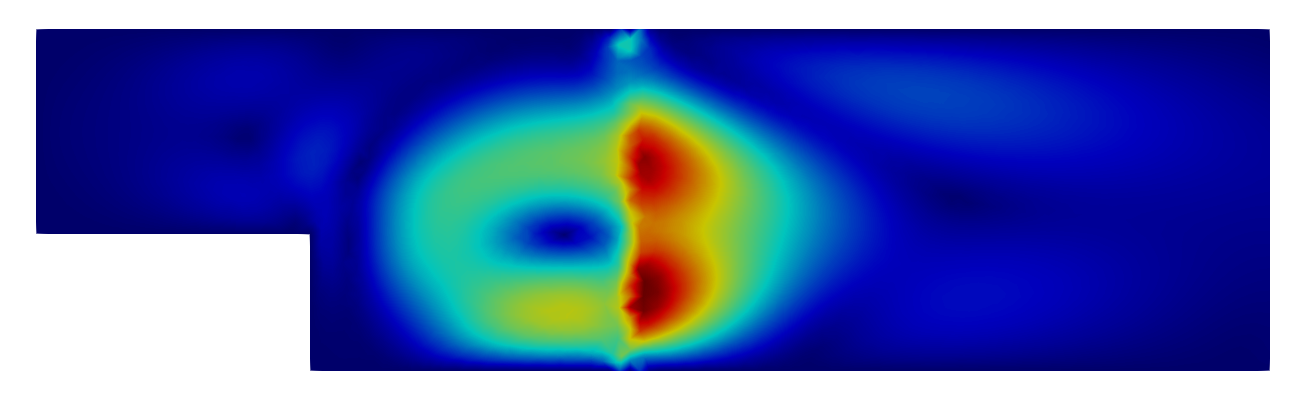}
         \label{fig:mode_u_2_bfs}
         \caption{\reviewerA{The second mode}}
    \end{subfigure}
    \begin{subfigure}[b]{0.49\textwidth}
        \includegraphics[width=\textwidth]{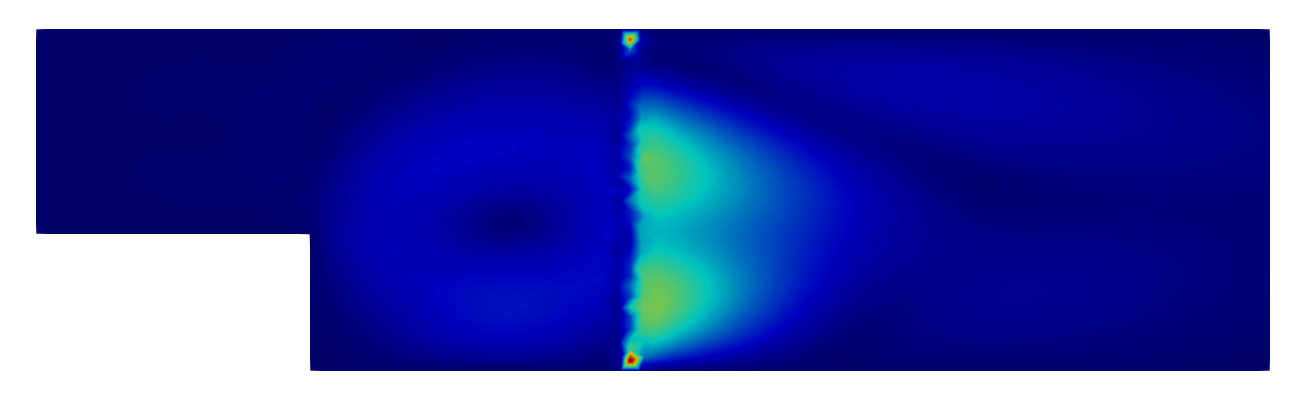}
        \caption{\reviewerA{The third mode}}
         \label{fig:mode_u_3_bfs}
    \end{subfigure}
    \hfill
    \begin{subfigure}[b]{0.49\textwidth}
        \includegraphics[width=\textwidth]{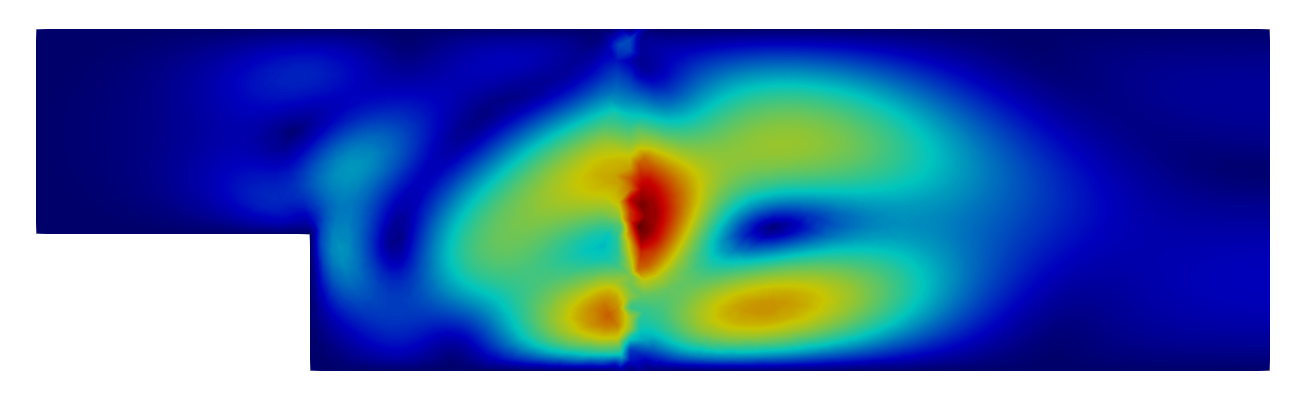}
         \caption{\reviewerA{The fourth mode}}
         \label{fig:mode_u_4_bfs}
    \end{subfigure}
    \caption{The first POD modes for the velocities $u_1$ and $u_2$ (subdomain functions are glued together for visualisation purposes).}
    \label{fig:pod_modes_u_bfs}
\end{figure}

\begin{figure}[H]
    \centering
    \begin{subfigure}[b]{0.49\textwidth}
        \includegraphics[width=\textwidth]{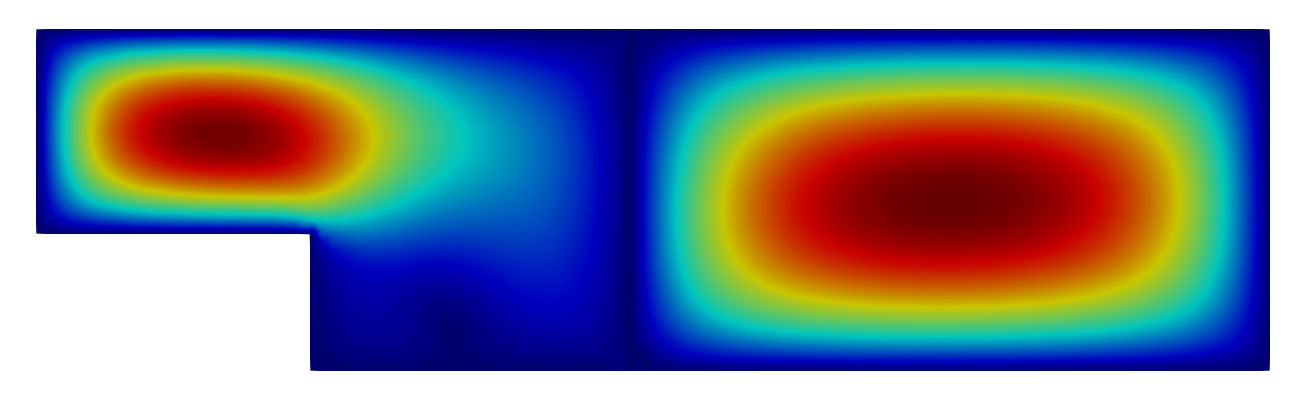}
         \label{fig:mode_s_1_bfs}
         \caption{\reviewerA{The first mode}}
    \end{subfigure}
    \hfill
    \begin{subfigure}[b]{0.49\textwidth}
        \includegraphics[width=\textwidth]{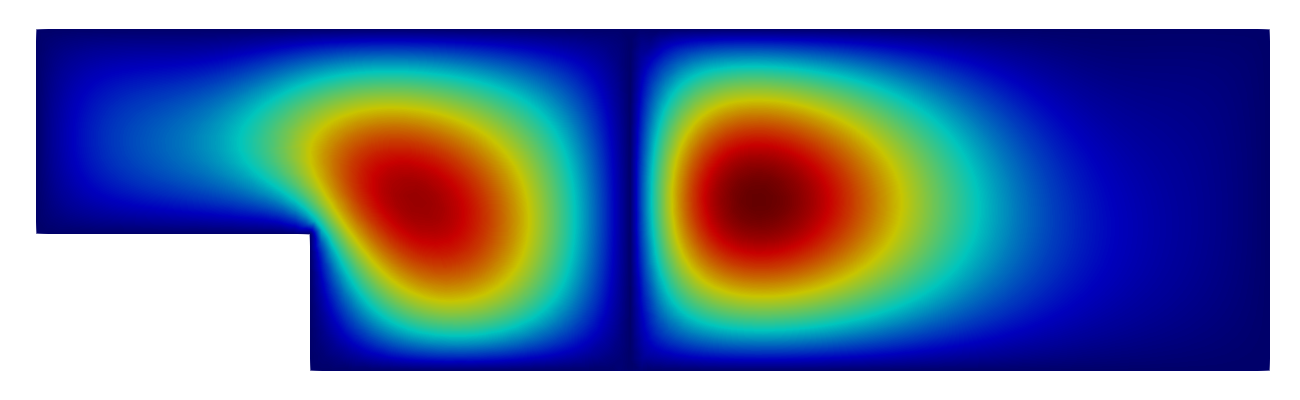}
         \label{fig:mode_s_2_bfs}
         \caption{\reviewerA{The second mode}}
    \end{subfigure}
    \begin{subfigure}[b]{0.49\textwidth}
        \includegraphics[width=\textwidth]{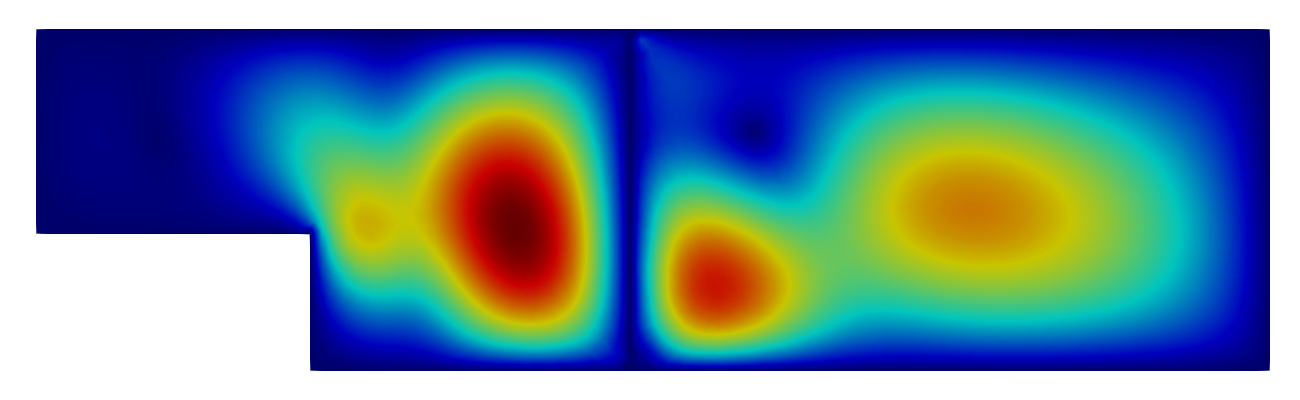}
         \caption{\reviewerA{The third mode}}
         \label{fig:mode_s_3_bfs}
    \end{subfigure}
    \hfill
    \begin{subfigure}[b]{0.49\textwidth}
        \includegraphics[width=\textwidth]{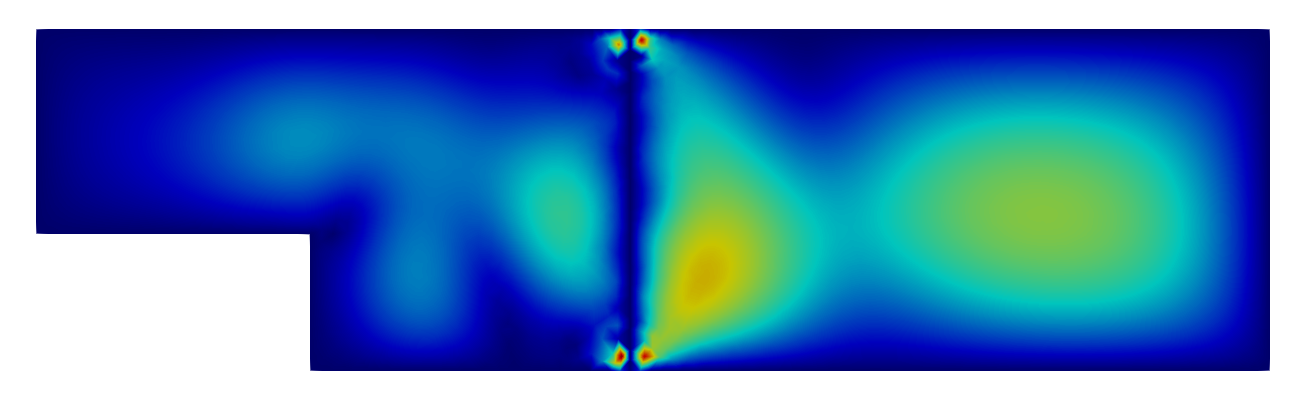}
         \caption{\reviewerA{The fourth mode}}
         \label{fig:mode_s_4_bfs}
    \end{subfigure}
    \caption{The first POD modes for the pressure supremisers $s_1$ and $s_2$ (subdomain functions are glued together for visualisation purposes).}
    \label{fig:pod_modes_s_bfs}
\end{figure}

\begin{figure}[H]
    \centering
    \begin{subfigure}[b]{0.49\textwidth}
        \includegraphics[width=\textwidth]{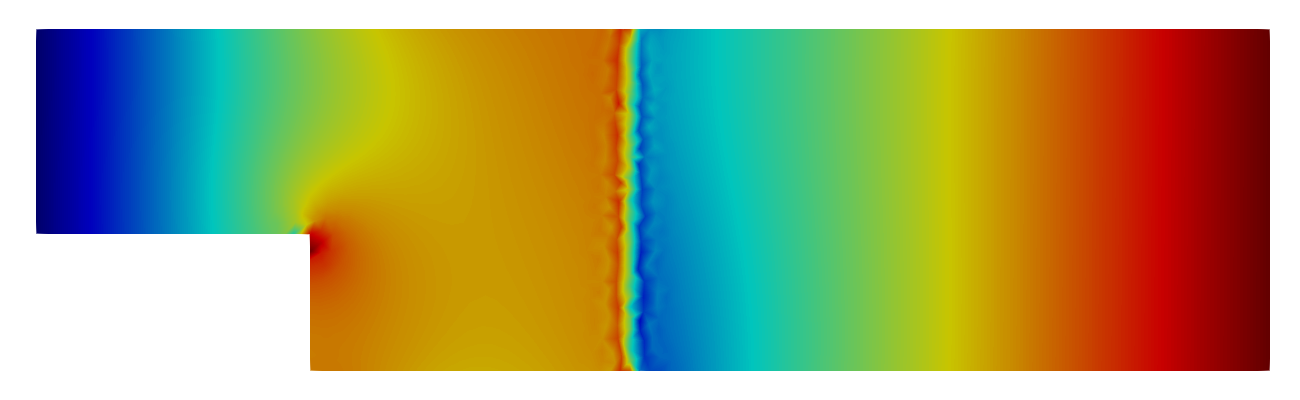}
         \caption{\reviewerA{The first mode}}
         \label{fig:mode_p_1_bfs}
    \end{subfigure}
    \hfill
    \begin{subfigure}[b]{0.49\textwidth}
        \includegraphics[width=\textwidth]{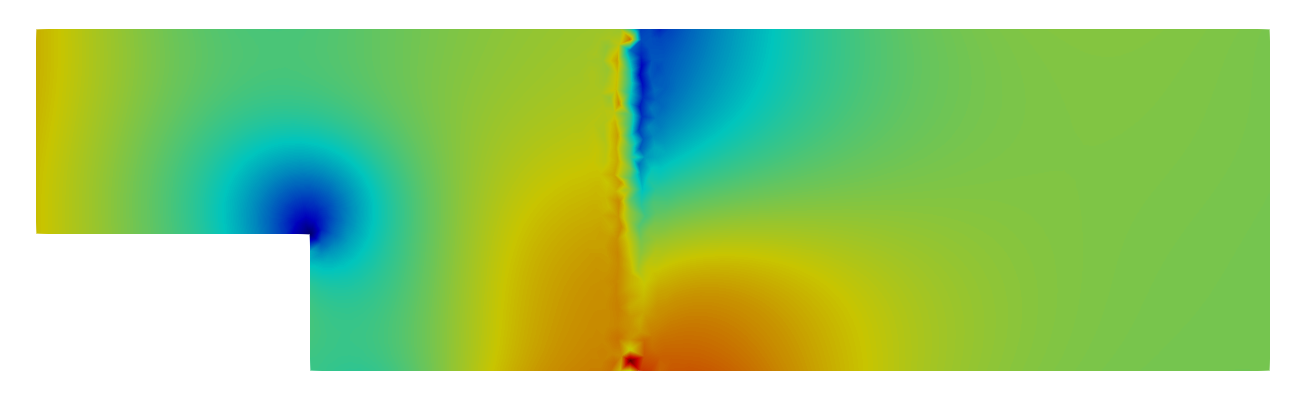}
         \caption{\reviewerA{The second mode}}
         \label{fig:mode_p_2_bfs}
    \end{subfigure}
    \begin{subfigure}[b]{0.49\textwidth}
        \includegraphics[width=\textwidth]{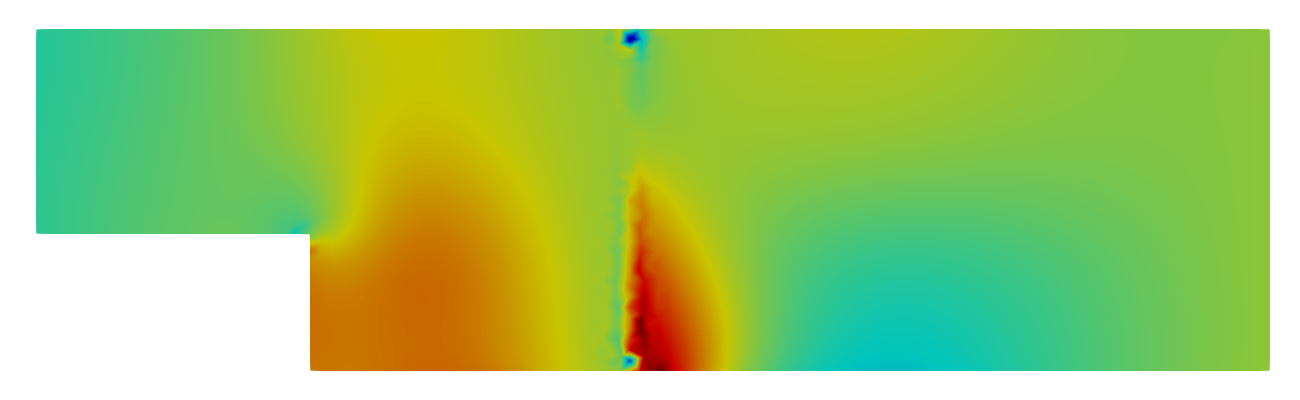}
         \caption{\reviewerA{The third mode}}
         \label{fig:mode_p_3_bfs}
    \end{subfigure}
    \hfill
    \begin{subfigure}[b]{0.49\textwidth}
        \includegraphics[width=\textwidth]{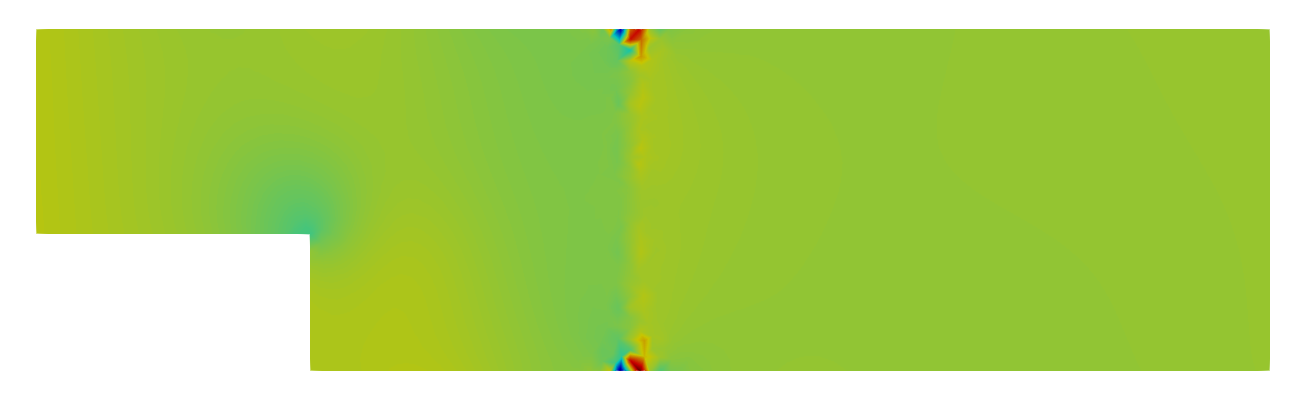}
         \caption{\reviewerA{The fourth mode}}
         \label{fig:mode_p_4_bfs}
    \end{subfigure} 
    \caption{The first POD modes for the pressures $p_1$ and $p_2$ (subdomain functions are glued together for visualisation purposes).}
    \label{fig:pod_modes_p_bfs}
\end{figure}

\begin{figure}[H]
    \centering
    \begin{subfigure}[b]{0.49\textwidth}
        \includegraphics[width=\textwidth]{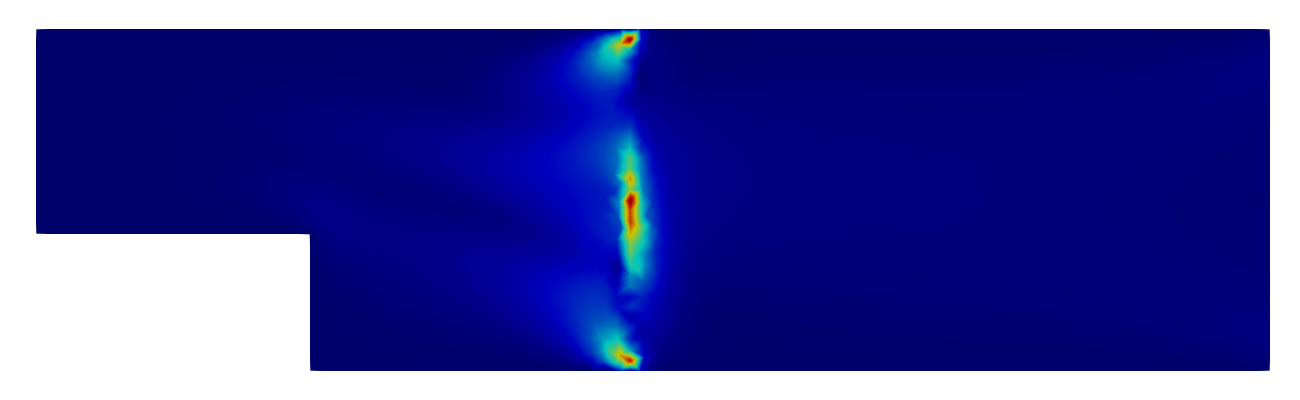}
         \caption{\reviewerA{The first mode}}
         \label{fig:mode_x_1_bfs}
    \end{subfigure}
    \hfill
    \begin{subfigure}[b]{0.49\textwidth}
        \includegraphics[width=\textwidth]{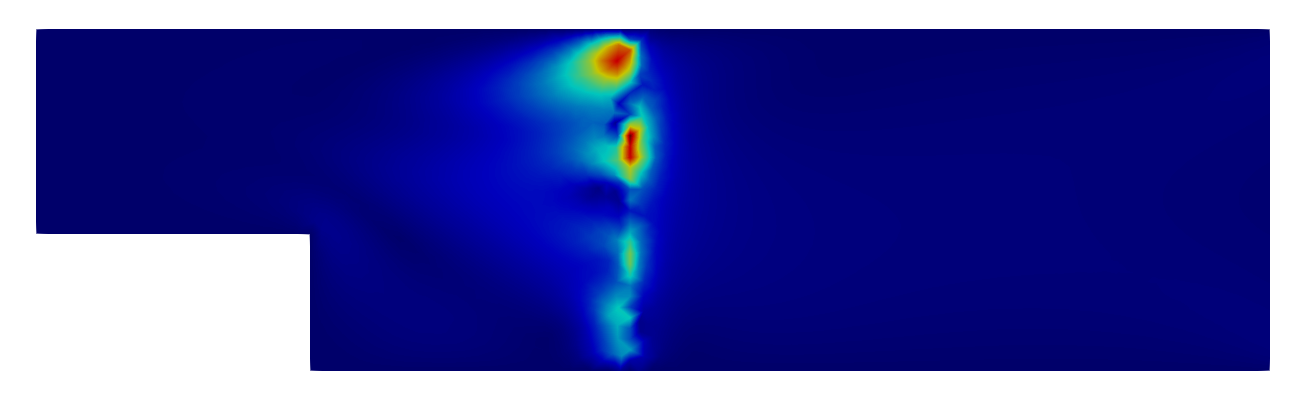}
         \caption{\reviewerA{The second mode}}
         \label{fig:mode_x_2_bfs}
    \end{subfigure}
    \begin{subfigure}[b]{0.49\textwidth}
        \includegraphics[width=\textwidth]{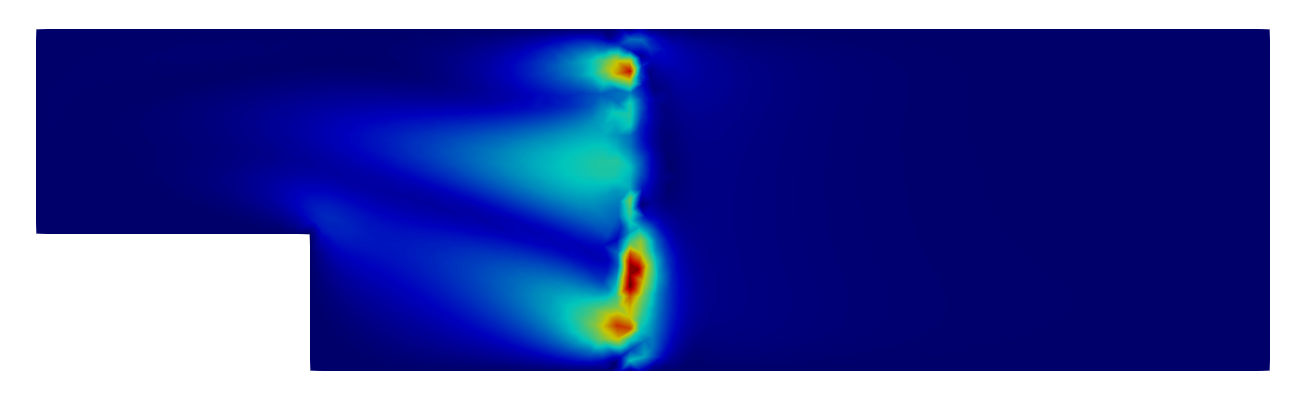}
         \caption{\reviewerA{The third mode}}
         \label{fig:mode_x_3_bfs}
    \end{subfigure}
    \hfill
    \begin{subfigure}[b]{0.49\textwidth}
        \includegraphics[width=\textwidth]{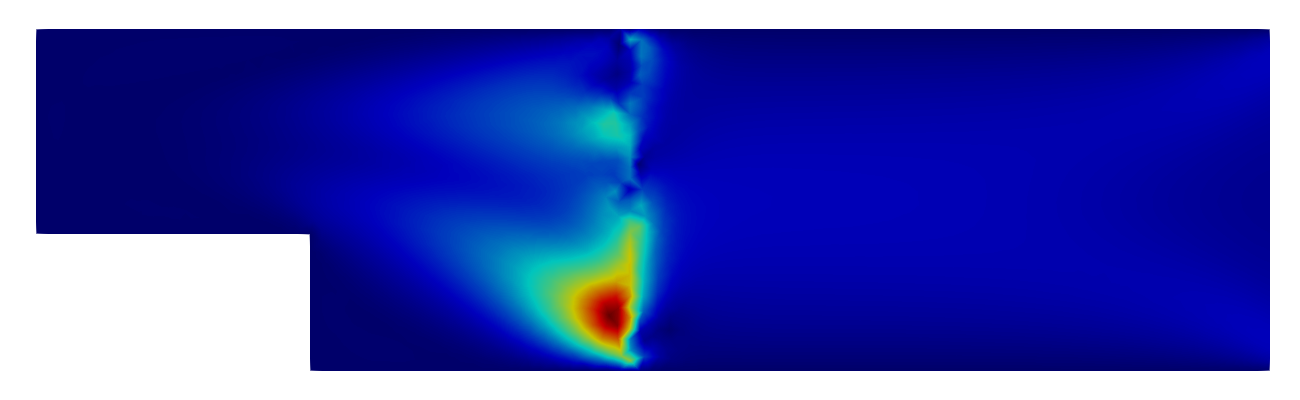}
         \caption{\reviewerA{The fourth mode}}
         \label{fig:mode_x_4_bfs}
    \end{subfigure}
    \caption{The first POD modes for the adjoint velocities $\xi_1$ and $\xi_2$ (subdomain functions are glued together for visualisation purposes).}
    \label{fig:pod_modes_x_bfs}
\end{figure}

Figures \ref{fig:pod_modes_u_bfs}--\ref{fig:pod_modes_x_bfs} represent the first four POD modes for each of the variables $u_1, u_2, s_1, s_2,$ $p_1, p_2, \xi_1$ and $\xi_2$. We stress that the POD modes were obtained separately for each component and the resulting figures are obtained by gluing the subdomain function just for the sake of visualisation.

Figure~\ref{fig:pod_modes_u_bfs} shows the first modes for the fluid velocities $u_1$ and $u_2$: in particular, notice that the modes corresponding to $u_1$ (on the left section of the domain) are zero at the inlet boundary due to the use of lifting function. 

In Figure~\ref{fig:pod_modes_s_bfs}, we can see the first four modes for $s_1$ and $s_2$: here, the corresponding functions are mostly localised inside the domains $\Omega_1$ and $\Omega_2$ thanks to the homogeneous conditions at the boundaries and the non-zero forcing term coming from the pressure. 

Figure~\ref{fig:pod_modes_p_bfs} represents the first modes for the pressures $p_1$ and $p_2$: we point out the signs of the oscillation behaviour, which suggests that the supremiser enrichment might be needed to assure stability of the reduced--order solution. 
Finally, Figure~\ref{fig:pod_modes_x_bfs} shows the first four modes for the adjoint variables  $\xi_1$ and $\xi_2$: note that they are concentrated only around the interface $\Gamma_0$ because the only nonzero contribution in the adjoint equations is coming from the source terms, which are defined solely on the interface $\Gamma_0$.

\begin{figure}[H]
    \centering
    \begin{subfigure}[b]{0.49\textwidth}
        \includegraphics[width=\textwidth]{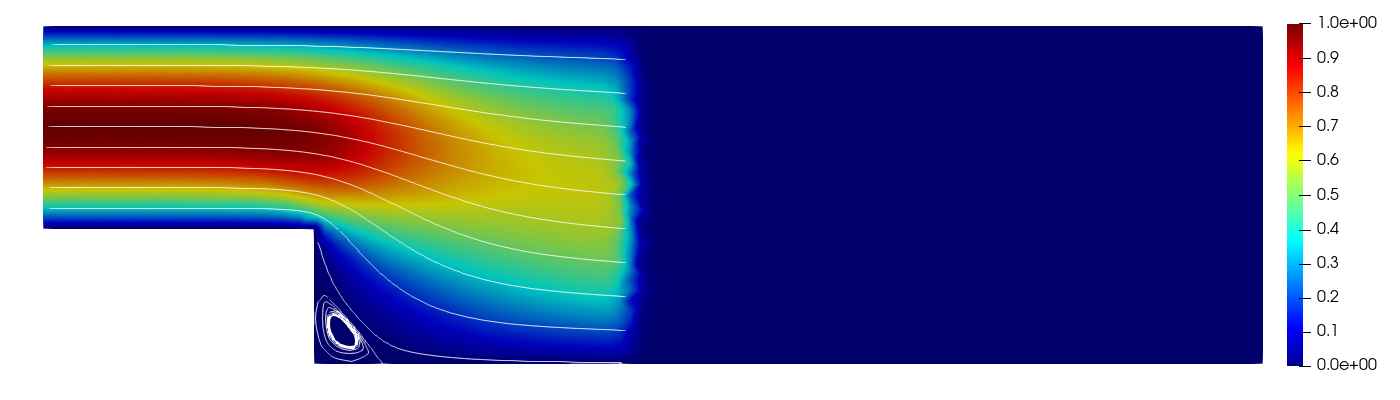}
        \caption{Iteration 0}
         \label{fig:truth_u_11_0_bfs}
    \end{subfigure}
    \hfill
    \begin{subfigure}[b]{0.49\textwidth}
        \includegraphics[width=\textwidth]{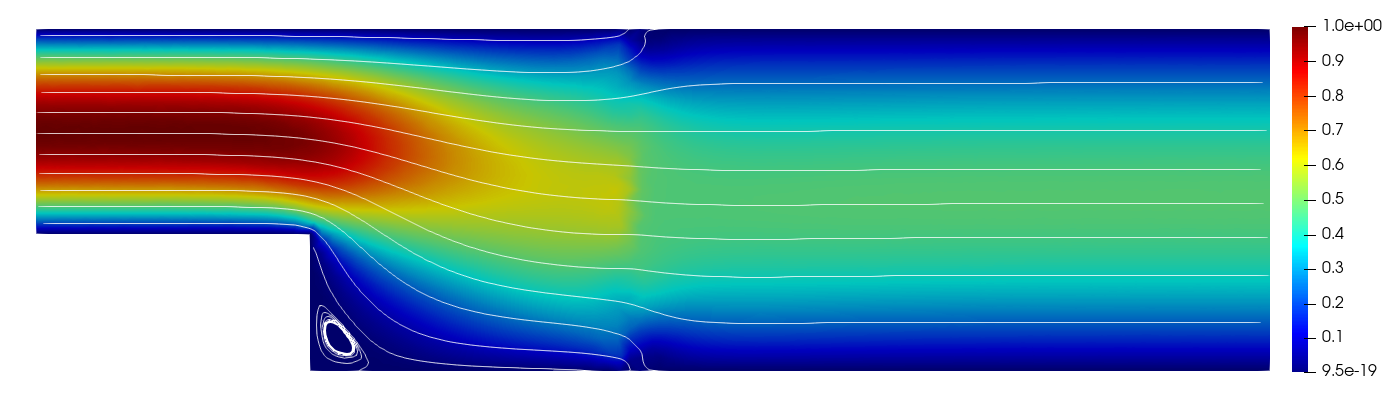}
         \caption{Iteration 5}
         \label{fig:truth_u_11_5_bfs}
    \end{subfigure}
    \begin{subfigure}[b]{0.49\textwidth}
        \includegraphics[width=\textwidth]{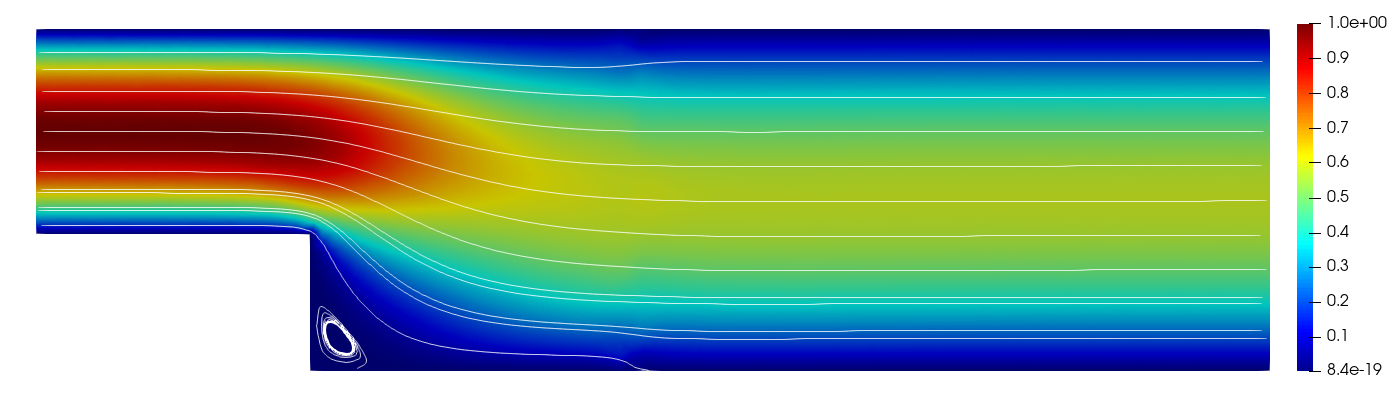}
        \caption{Iteration 10}
         \label{fig:truth_u_11_10_bfs}
    \end{subfigure}
    \hfill
    \begin{subfigure}[b]{0.49\textwidth}
        \includegraphics[width=\textwidth]{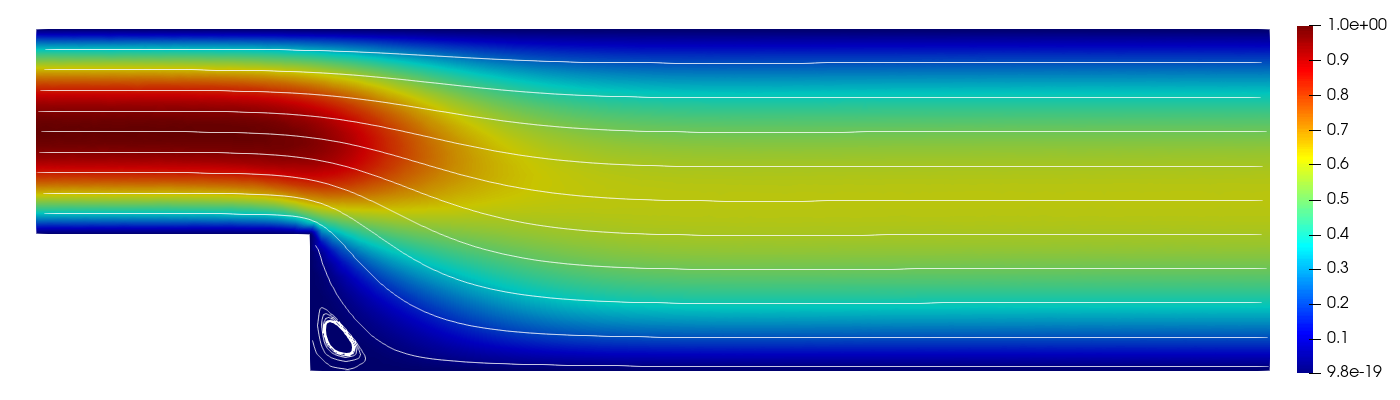}
         \caption{Iteration 40}
         \label{fig:truth_u_11_40_bfs}
    \end{subfigure}
    \caption{High--fidelity solution for the velocities $u_1$ and $u_2$. Values of the parameters $\bar U=1$, $\nu=1$ \reviewerA{and $Re=3$}}
    \label{fig:truth_u_11_bfs}
\end{figure}

 \begin{figure}[H]
    \centering
    \begin{subfigure}[b]{0.49\textwidth}
        \includegraphics[width=\textwidth]{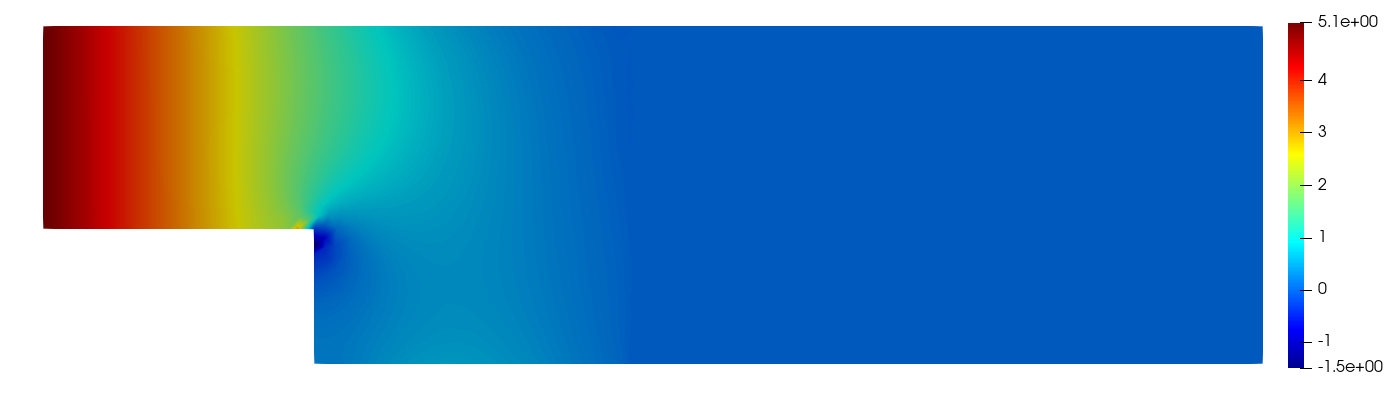}
        \caption{Iteration 0}
         \label{fig:truth_p_11_0_bfs}
    \end{subfigure}
    \hfill
    \begin{subfigure}[b]{0.49\textwidth}
        \includegraphics[width=\textwidth]{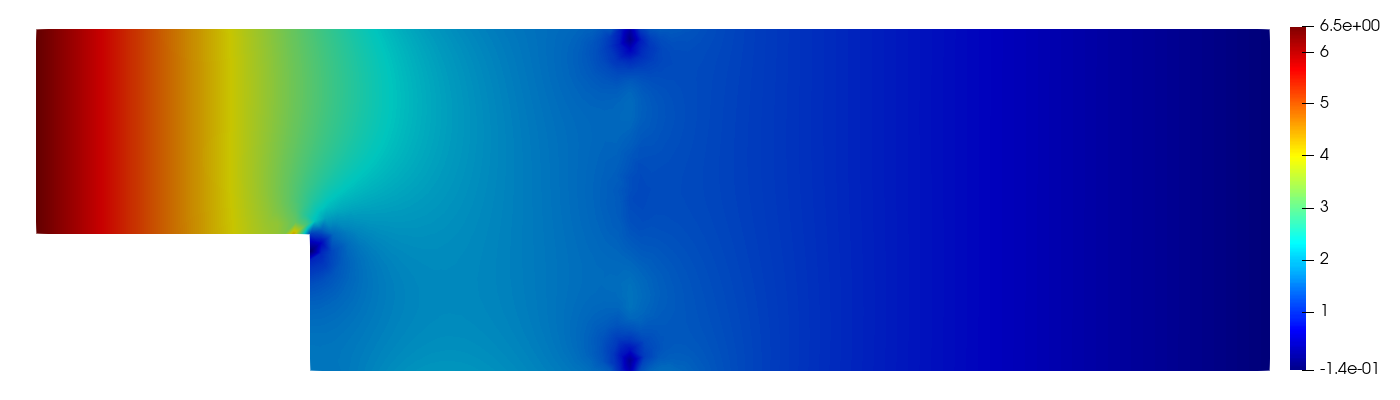}
         \caption{Iteration 5}
         \label{fig:truth_p_11_5_bfs}        
    \end{subfigure}
    \begin{subfigure}[b]{0.49\textwidth}
        \includegraphics[width=\textwidth]{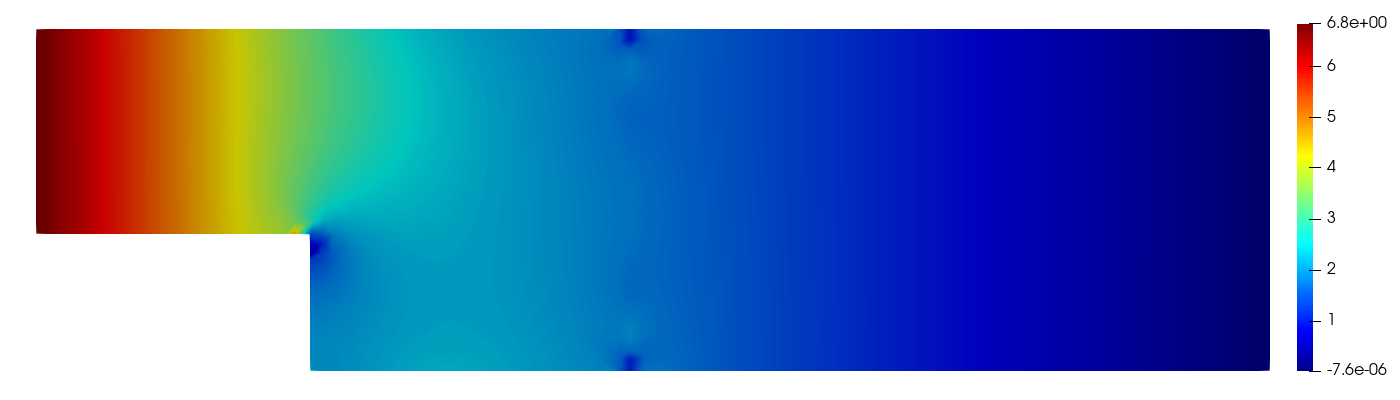}
        \caption{Iteration 10}
         \label{fig:truth_p_11_10_bfs}
    \end{subfigure}
    \hfill
    \begin{subfigure}[b]{0.49\textwidth}
        \includegraphics[width=\textwidth]{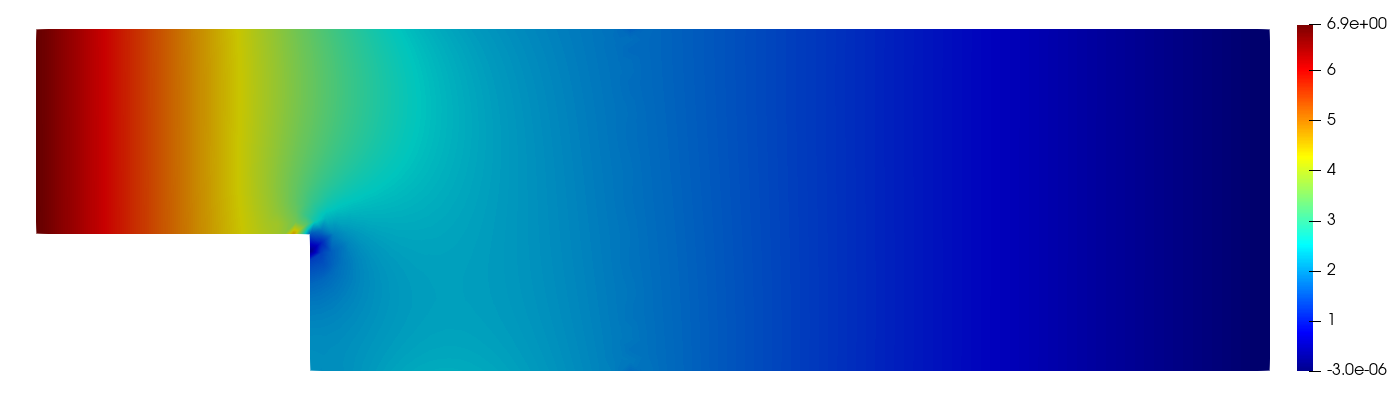}
         \caption{Iteration 40}
         \label{fig:truth_p_11_40_bfs}
    \end{subfigure}
    \caption{High--fidelity solution for the pressures $p_1$ and $p_2$. Values of the parameters $\bar U=1$, $\nu=1$ \reviewerA{and $Re=3$}}
    \label{fig:truth_p_11_bfs}
\end{figure}

\begin{figure}[H]
    \centering
    \begin{subfigure}[b]{0.49\textwidth}
        \includegraphics[width=\textwidth]{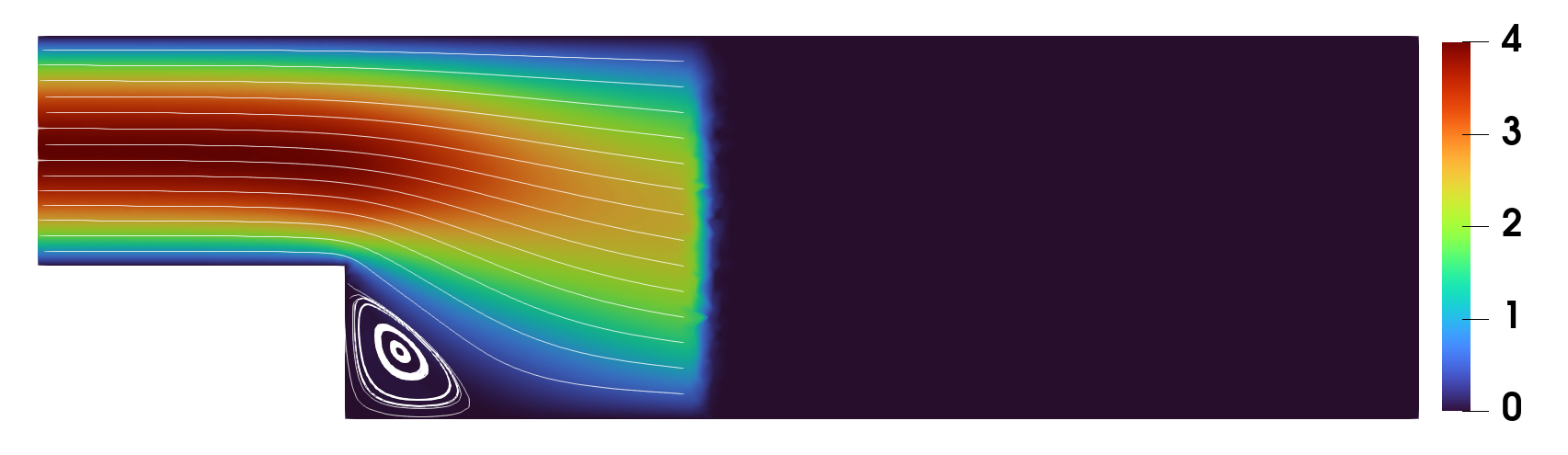}
        \caption{Iteration 0}
         \label{fig:truth_u_22_0_bfs}
    \end{subfigure}
    \hfill
    \begin{subfigure}[b]{0.49\textwidth}
        \includegraphics[width=\textwidth]{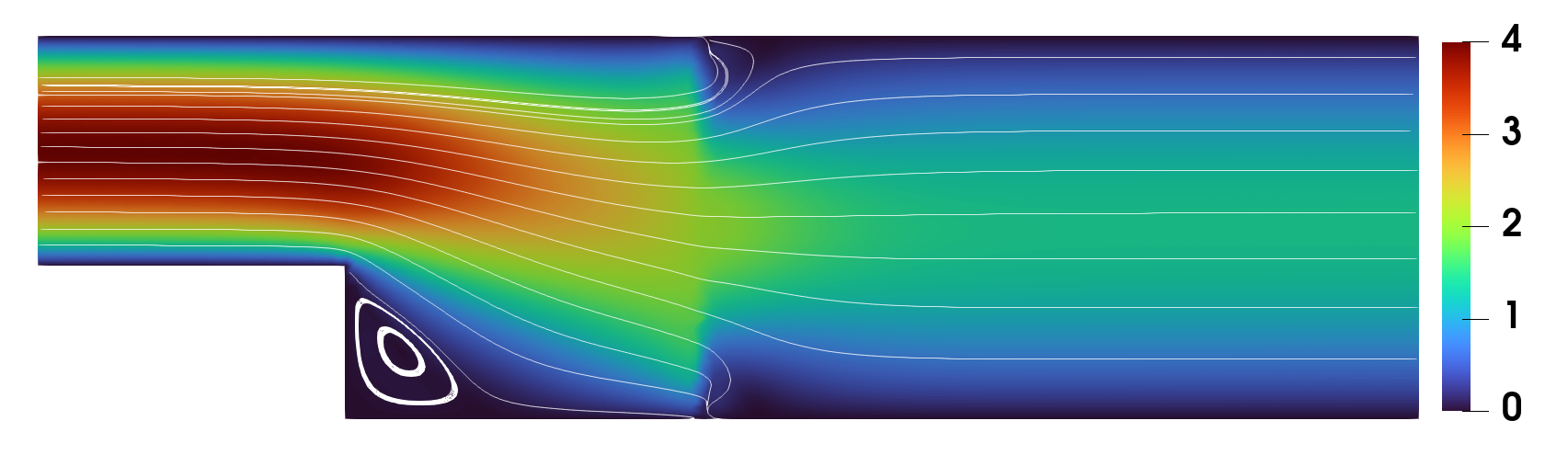}
         \caption{Iteration 5}
         \label{fig:truth_u_22_5_bfs}

    \end{subfigure}
    
    \begin{subfigure}[b]{0.49\textwidth}
        \includegraphics[width=\textwidth]{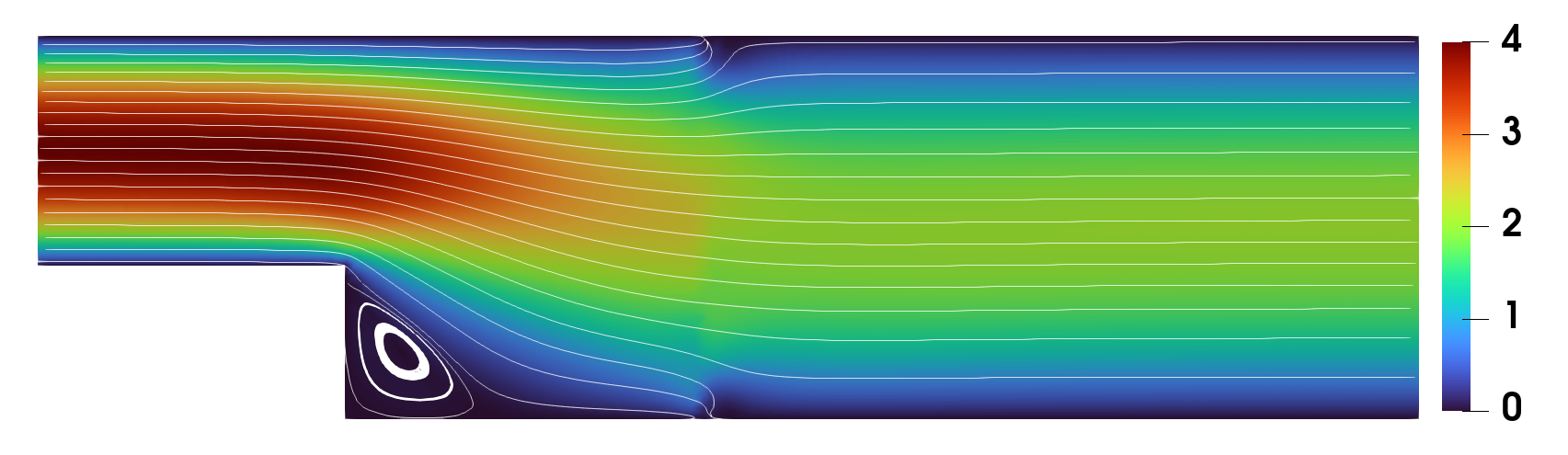}
        \caption{Iteration 10}
         \label{fig:truth_u_22_10_bfs}
    \end{subfigure}
    \hfill
    \begin{subfigure}[b]{0.49\textwidth}
        \includegraphics[width=\textwidth]{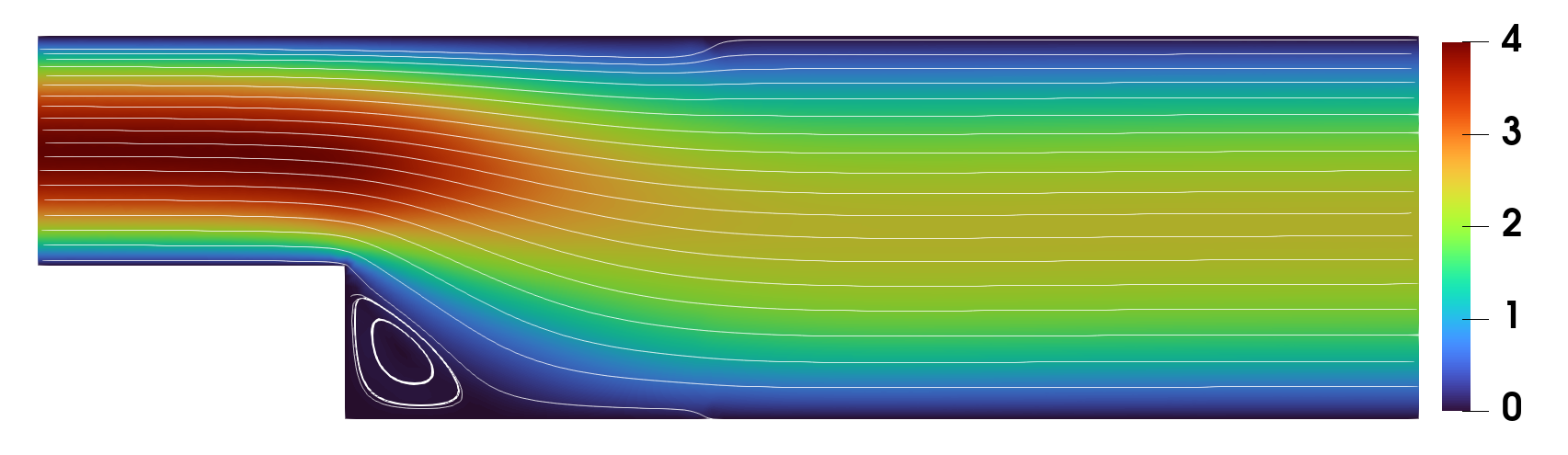}
         \caption{Iteration 40}
         \label{fig:truth_u_22_40_bfs}  
    \end{subfigure}
    \caption{High--fidelity solution for the velocities $u_1$ and $u_2$. Values of the parameters \reviewerA{$\bar U=4$, $\nu=0.75$ and $Re\approx 19$}}
    \label{fig:truth_u_22_bfs}
\end{figure}

\begin{figure}[H]
    \centering
    \begin{subfigure}[b]{0.49\textwidth}
        \includegraphics[width=\textwidth]{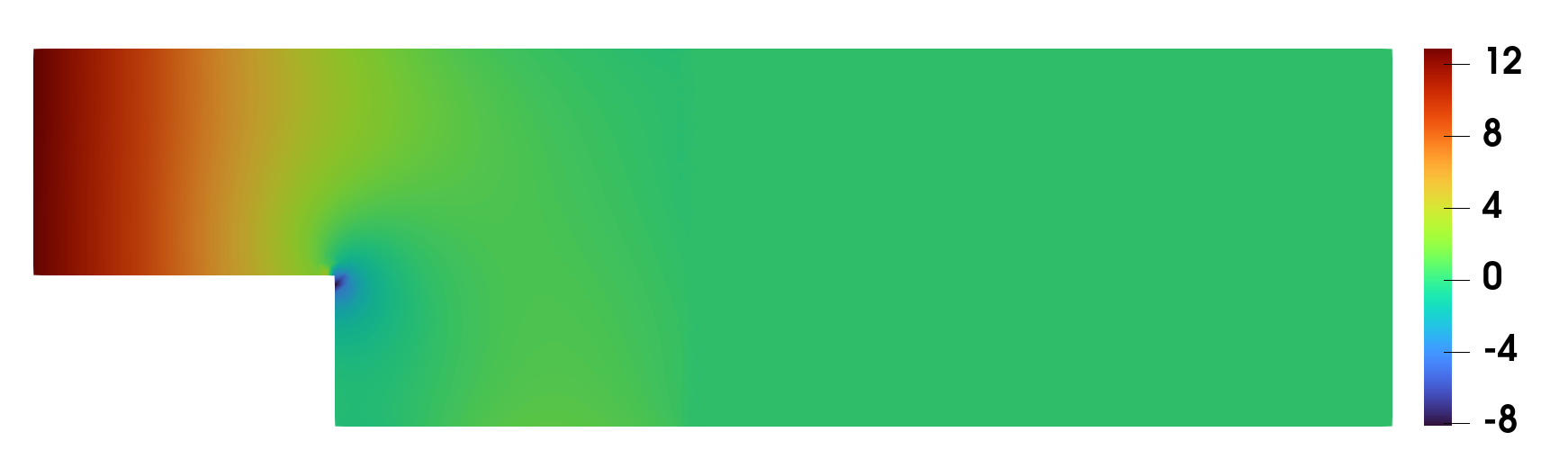}
        \caption{Iteration 0}
         \label{fig:truth_p_22_0_bfs}
    \end{subfigure}
    \hfill
    \begin{subfigure}[b]{0.49\textwidth}
        \includegraphics[width=\textwidth]{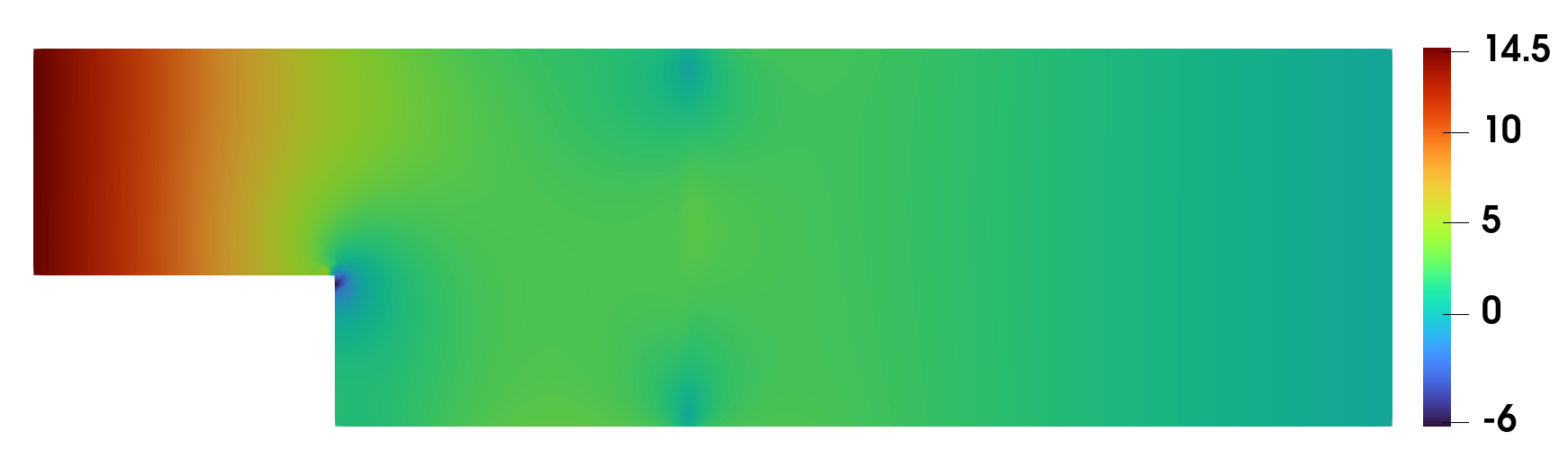}
         \caption{Iteration 5}
         \label{fig:truth_p_22_5_bfs} 
    \end{subfigure}
    \begin{subfigure}[b]{0.49\textwidth}
        \includegraphics[width=\textwidth]{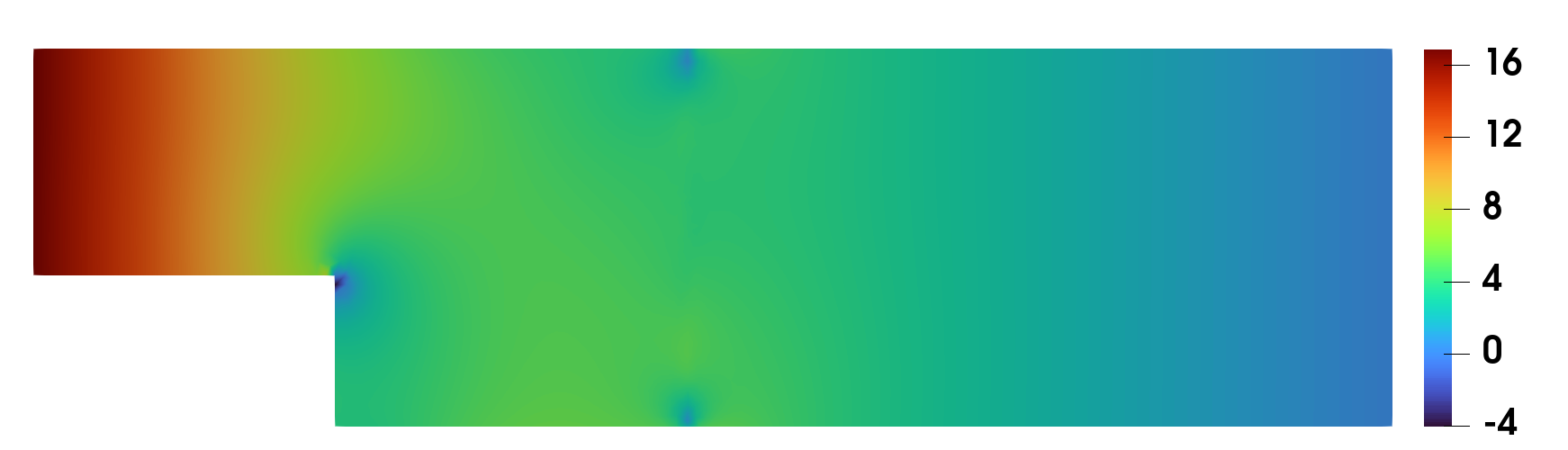}
        \caption{Iteration 10}
         \label{fig:truth_p_22_10_bfs}
    \end{subfigure}
    \hfill
    \begin{subfigure}[b]{0.49\textwidth}
        \includegraphics[width=\textwidth]{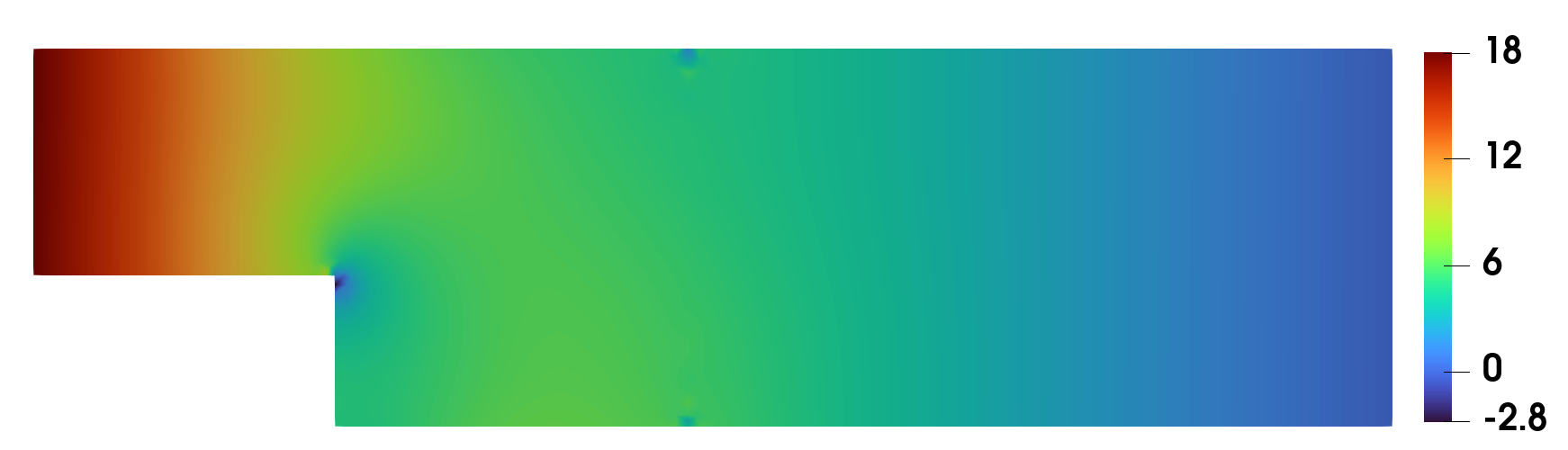}
         \caption{Iteration 40}
         \label{fig:truth_p_22_40_bfs}         
    \end{subfigure}
    \caption{High--fidelity solution for the pressures $p_1$ and $p_2$. Values of the parameters \reviewerA{$\bar U=4$, $\nu=0.75$ and $Re\approx 19$}}
    \label{fig:truth_p_22_bfs}
\end{figure}

Figures~\ref{fig:truth_u_11_bfs}-\ref{fig:truth_p_22_bfs} represent the high--fidelity solutions for two different values of the parameters $(\bar U, \nu)=(1,1)$, \reviewerA{resulting in $Re=3$}, and $(\bar U, \nu)=(4.5,0.7)$ \reviewerA{with $Re\approx 19$}. The solutions were obtained by carrying out 40 optimisation iterations via L--BFGS--B algorithm. Figures \ref{fig:truth_u_11_bfs} and \ref{fig:truth_u_22_bfs} show the intermediate solutions at iteration 0, 5, 10 and 40 for the fluid velocities $u_1$ and $u_2$, whereas Figures \ref{fig:truth_p_11_bfs} and \ref{fig:truth_p_22_bfs} show the corresponding pressures $p_1$ and $p_2$. The final solution is taken to be the 40th iteration optimisation solution in which we can observe a continuity between subdomain solutions at the interface $\Gamma_0$. Moreover, it can be noticed that the solution for parameters $(\bar U, \nu)=(1,1)$ looks continuous already at iteration 10, which suggests that the convergence of the optimisation algorithm might depend on the Reynolds number.

\begin{table}[H]
    \centering
    
    \begin{tabular}{|c|c|c|}
      \hline
        \textbf{Iteration} & \textbf{Functional Value}  &\textbf{Gradient norm }    \\
        \hline  
         0 & $4.8 \cdot 10^{-1}$ & $4.1 \cdot 10^{-1}$ \\
         5 & $6.0 \cdot 10^{-2}$ & $2.2 \cdot 10^{-1}$  \\ 
         10 & $5.0 \cdot 10^{-3}$ & $3.3 \cdot 10^{-2}$  \\ 
         40 & $1.7 \cdot 10^{-4}$ & $2.4 \cdot 10^{-3}$  \\
         \hline  
    \end{tabular}
    \caption{Functional values and the gradient norm for the \reviewerA{FOM} optimisation solution at the parameter values $\bar U=1$, $\nu=1$ \reviewerA{and $Re=3$}}
    \label{tab:truth_1_func_bfs}
\end{table}

\begin{table}[H]
    \centering 
     \begin{adjustbox}{max width=\textwidth}
    \begin{tabular}{|c|c|c|c|c|c|c|c|c|}
      \hline
         \textbf{Iteration}  &\multicolumn{2}{|c|}{\textbf{Abs. error $u_{h}$}} &\multicolumn{2}{|c|}{\textbf{Rel. error $u_{h}$}} &\multicolumn{2}{|c|}{\textbf{Abs. error $p_{h}$} }&\multicolumn{2}{|c|}{\textbf{Rel. error $p_{h}$}}   \\ \hline  
         & $\Omega_1$&$\Omega_2$& $\Omega_1$&$\Omega_2$& $\Omega_1$&$\Omega_2$& $\Omega_1$&$\Omega_2$\\ \hline
         0 & 0.0302 & 2.9935 & 0.0088 & 1.0000 & 10.6515 & 7.0679 & 0.5056 & 1.0000 \\
        5 & 0.1020 & 0.6279 & 0.0297 & 0.2098 & 2.4520 & 1.5317 & 0.1164 & 0.2167 \\
        10 & 0.0384 & 0.1355 & 0.0112 & 0.0453 & 0.5807 & 0.3793 & 0.0276 & 0.0537 \\
        40 & 0.0184 & 0.0583 & 0.0053 & 0.0195 & 0.2670 & 0.1827 & 0.0127 & 0.0259 \\
         \hline  
    \end{tabular}
    \end{adjustbox}
     \caption{Absolute and relative errors of the \reviewerA{FOM} optimisation solution with respect to the monolithic solution at the parameter values $\bar U=1$, $\nu=1$ \reviewerA{and $Re=3$}}
    \label{tab:truth_1_errors_bfs}
\end{table}

We present additional details in Tables \ref{tab:truth_1_func_bfs} - \ref{tab:truth_2_errors_bfs}. In particular, in Tables \ref{tab:truth_1_func_bfs} and \ref{tab:truth_2_func_bfs}, we list the values for the functional $\mathcal J_\gamma$ and the $L^2(\Gamma_0)$-norm  of the gradient $\frac{d \mathcal J_\gamma}{dg}$ at the different iteration of the optimisation procedure, while Table \ref{tab:truth_1_errors_bfs} contains the absolute and relative errors with respect to the monolithic (entire--domain) solutions $u_h, p_h$, i.e.,

\begin{eqnarray*}
    \text{Abs. error} \ u_{h} := || u_{i,h} - u_h ||_{L^2(\Omega_i)} & \text{on domain $\Omega_i$}, 
\\
    \text{Rel. error} \ u_{h} : = \frac{  || u_{i,h} - u_h ||_{L^2(\Omega_i)}}{|| u_h ||_{L^2(\Omega_i)}} & \text{on domain $\Omega_i$}, 
 \end{eqnarray*}
    \begin{eqnarray*}
     \text{Abs. error} \ p_{h} := || p_{i,h} - p_h ||_{L^2(\Omega_i)} & \text{on domain $\Omega_i$}, 
\\
    \text{Rel. error} \ p_{h} : = \frac{||p_{i,h} - p_h ||_{L^2(\Omega_i)}}{|| p_h ||_{L^2(\Omega_i)}} & \text{on domain $\Omega_i$}, 
 \end{eqnarray*}
for $i=1,2$.

\begin{table}[H]
    \centering  
    \begin{tabular}{|c|c|c|}
      \hline
        \textbf{Iteration} & \textbf{Functional Value}  &\textbf{Gradient norm }    \\
        \hline  
         0 & $7.902 $ & $2.213$ \\
         5 & $1.956$ & $1.210 $  \\ 
         10 & $0.403$ & $2.132$  \\ 
         40 & $0.007$ & $0.069$  \\
         \hline  
    \end{tabular}
    \caption{Functional values and the gradient norm for the \reviewerA{FOM} optimisation solution at parameter values \reviewerA{$\bar U=4$, $\nu=0.75$ and $Re\approx 19$}}
    \label{tab:truth_2_func_bfs}
\end{table}

\begin{table}[H]
    \centering
    \begin{adjustbox}{max width=\textwidth}
    \begin{tabular}{|c|c|c|c|c|c|c|c|c|}
      \hline
         \textbf{Iteration}  &\multicolumn{2}{|c|}{\textbf{Abs. error $u_{h}$}} &\multicolumn{2}{|c|}{\textbf{Rel. error $u_{h}$}} &\multicolumn{2}{|c|}{\textbf{Abs. error $p_{h}$} }&\multicolumn{2}{|c|}{\textbf{Rel. error $p_{h}$}}   \\ \hline  
         & $\Omega_1$&$\Omega_2$& $\Omega_1$&$\Omega_2$& $\Omega_1$&$\Omega_2$& $\Omega_1$&$\Omega_2$\\ \hline
        0 & 0.2520 & 11.9830 & 0.0181 & 1.0000 & 31.6121 & 21.1630 & 0.5859 & 1.0000 \\
        5 & 0.6639 & 5.0075 & 0.0478 & 0.4179 & 20.7060 & 10.2359 & 0.3838 & 0.4837 \\
        10 & 0.2704 & 1.3722 & 0.0195 & 0.1145 & 6.7317 & 2.8262 & 0.1248 & 0.1335 \\
        40 & 0.0865 & 0.2566 & 0.0062 & 0.0214 & 1.4498 & 0.6443 & 0.0269 & 0.0304 \\
         \hline  
    \end{tabular}
    \end{adjustbox}
     \caption{Absolute and relative errors of the \reviewerA{FOM} optimisation solution with respect to the monolithic solution at the parameter values \reviewerA{$\bar U=4$, $\nu=0.75$ and $Re\approx 19$}}
    \label{tab:truth_2_errors_bfs}
\end{table}



\begin{table}[H]
    \centering
    
    \begin{tabular}{|c|c|c|}
      \hline
        \textbf{Iteration} & \textbf{Functional Value}  &\textbf{Gradient norm }    \\
        \hline  
         0 & $4.8 \cdot 10^{-1}$  & 0.391 \\
         5 & $5.4 \cdot 10^{-3}$ & 0.047   \\ 
         10 & $3.6 \cdot 10^{-4}$ & 0.015  \\ 
         \hline  
    \end{tabular}
    \caption{Functional values and the gradient norm for the \reviewerA{ROM} optimisation solution at parameter values $\bar U=1$, $\nu=1$ \reviewerA{and $Re=3$}}
    \label{tab:reduced_1_func_bfs}
\end{table}

\begin{table}[H]
    \centering
   \begin{adjustbox}{max width=\textwidth}
    \begin{tabular}{|c|c|c|c|c|c|c|c|c|}
      \hline
         \textbf{Iteration}  &\multicolumn{2}{|c|}{\textbf{Abs. error $u_{N}$}} &\multicolumn{2}{|c|}{\textbf{Rel. error $u_{N}$}} &\multicolumn{2}{|c|}{\textbf{Abs. error $p_{N}$} }&\multicolumn{2}{|c|}{\textbf{Rel. error $p_{N}$}}   \\ \hline  
         & $\Omega_1$&$\Omega_2$& $\Omega_1$&$\Omega_2$& $\Omega_1$&$\Omega_2$& $\Omega_1$&$\Omega_2$\\ \hline
           0 & 0.0284 & 2.9935 & 0.0083 & 1.0000& 10.9522 & 7.0679 & 0.5198 & 1.0000 \\
        5 & 0.0746 & 0.1956 & 0.0217 & 0.0653 & 0.8548 & 0.5672 & 0.0406 & 0.0803 \\
        10 & 0.0135 & 0.0357 & 0.0039 & 0.0119 & 0.1714 & 0.1186 & 0.0081 & 0.0168 \\
         \hline  
    \end{tabular}
    \end{adjustbox}
     \caption{Absolute and relative errors of the \reviewerA{ROM} optimisation solution with respect to the monolithic solution at the parameter values $\bar U=1$, $\nu=1$ \reviewerA{and $Re=3$}}
    \label{tab:reduced_1_errors_bfs}
\end{table}

\begin{figure}[H]
    \centering
    \begin{subfigure}[b]{0.49\textwidth}
        \includegraphics[width=\textwidth]{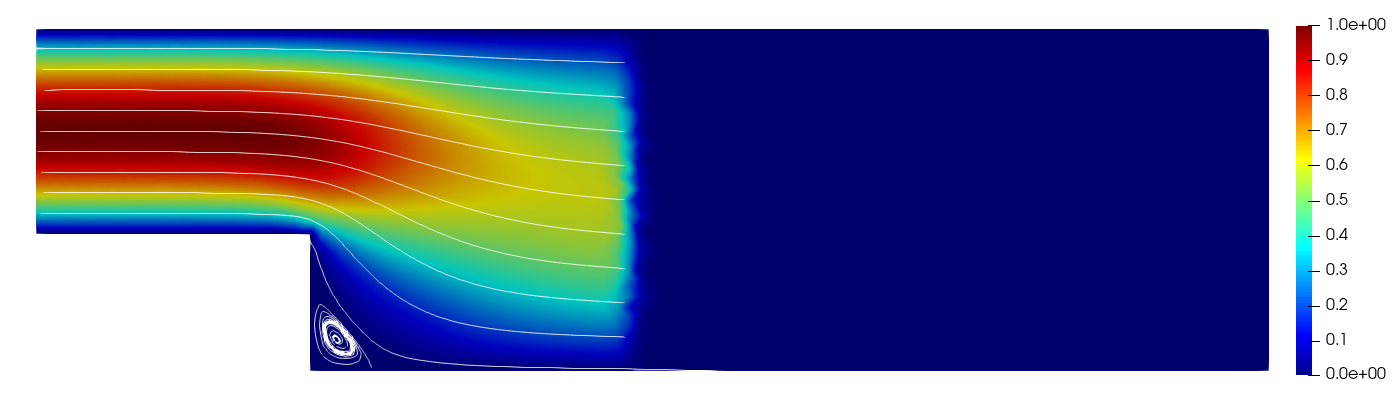}
        \caption{Iteration 0}
         \label{fig:reduced_u_11_0_bfs}
    \end{subfigure}
    \hfill
    \begin{subfigure}[b]{0.49\textwidth}
        \includegraphics[width=\textwidth]{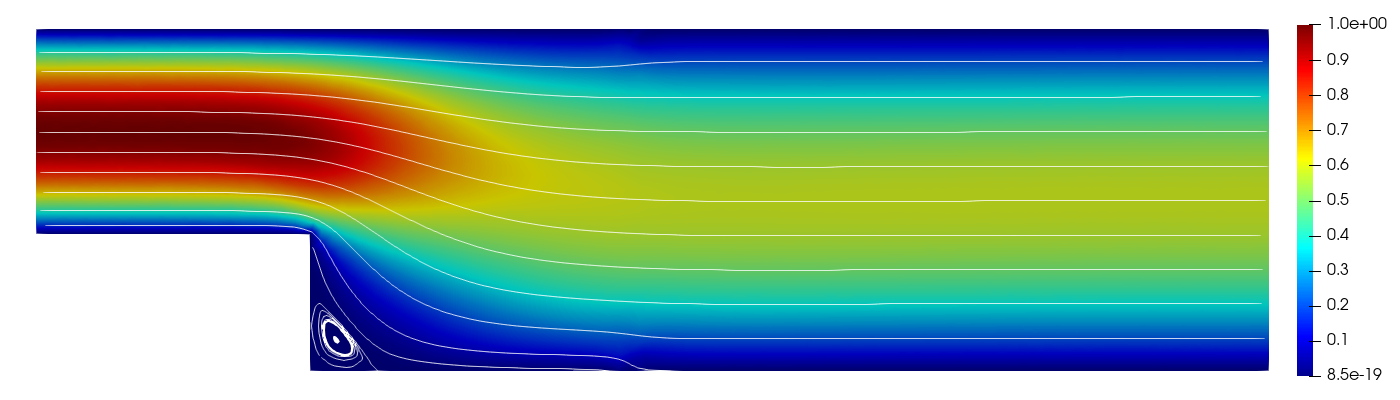}
         \caption{Iteration 5}
         \label{fig:reduced_u_11_5_bfs}

    \end{subfigure}
    
    \begin{subfigure}[b]{0.49\textwidth}
        \includegraphics[width=\textwidth]{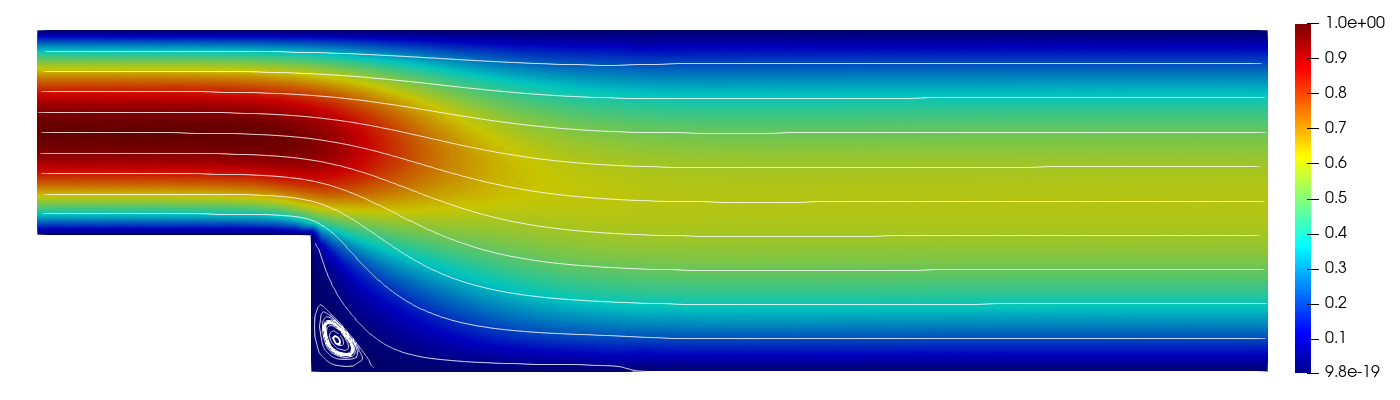}
        \caption{Iteration 10}
         \label{fig:reduced_u_11_10_bfs}
    \end{subfigure}

    \caption{Reduced order solution for the velocities $u_1$ and $u_2$. Values of the parameters $\bar U=1$, $\nu=1$ and \reviewerA{$Re=3$}. \reviewerA{Number of POD modes: 10 - for each state variable, each supremiser and the control, 30 -- for both adjoint velocities}}
    \label{fig:reduced_u_11_bfs}
\end{figure}

\begin{figure}[H]
    \centering
    \begin{subfigure}[b]{0.49\textwidth}
        \includegraphics[width=\textwidth]{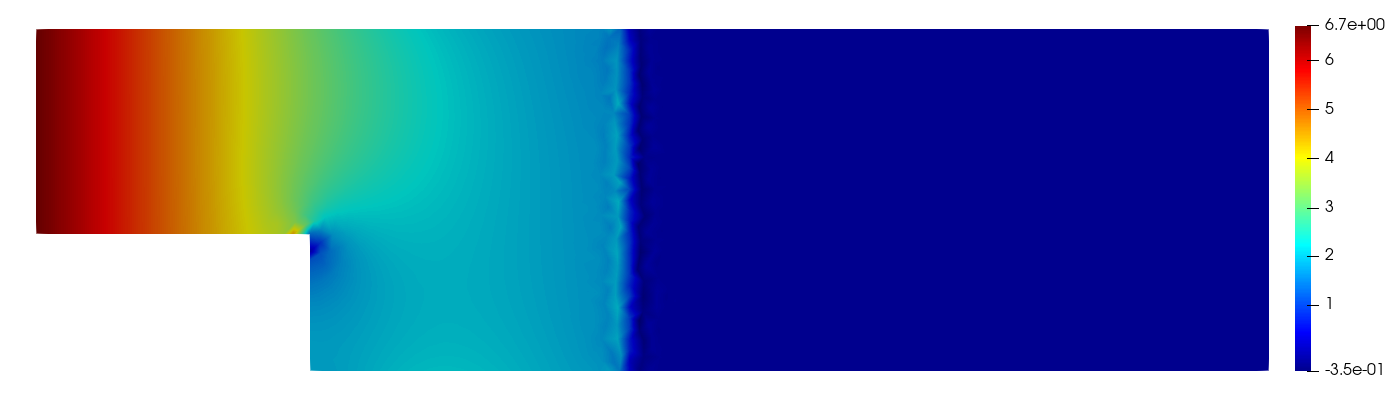}
        \caption{Iteration 0}
         \label{fig:reduced_p_11_0_bfs}
    \end{subfigure}
    \hfill
    \begin{subfigure}[b]{0.49\textwidth}
        \includegraphics[width=\textwidth]{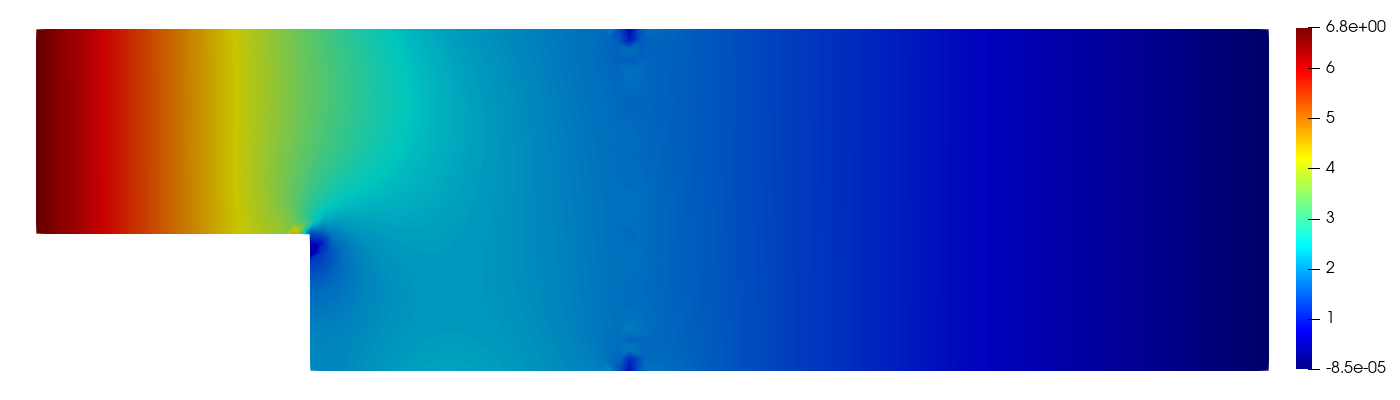}
         \caption{Iteration 5}
         \label{fig:reduced_p_11_5_bfs}

    \end{subfigure}
    
    \begin{subfigure}[b]{0.49\textwidth}
        \includegraphics[width=\textwidth]{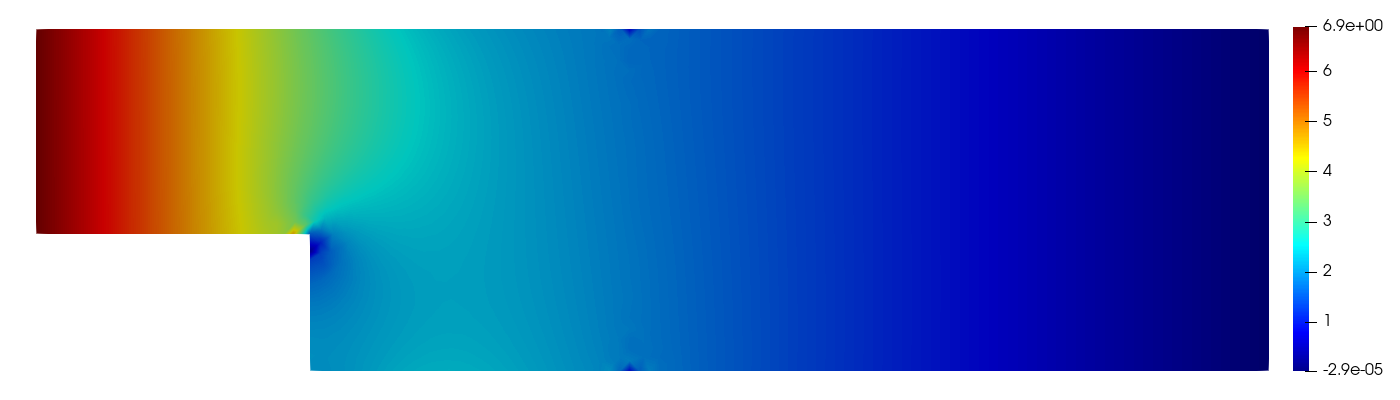}
        \caption{Iteration 10}
         \label{fig:reduced_p_11_10_bfs}
    \end{subfigure}

    \caption{Reduced order solution for the pressures $p_1$ and $p_2$. Values of the parameters $\bar U=1$, $\nu=1$ and \reviewerA{$Re=3$}.  \reviewerA{Number of POD modes: 10 - for each state variable, each supremiser and the control, 39 -- for both adjoint velocities}}
    \label{fig:reduced_p_11_bfs}
\end{figure}

\begin{figure}[H]
    \centering
    \begin{subfigure}[b]{0.49\textwidth}
        \includegraphics[width=\textwidth]{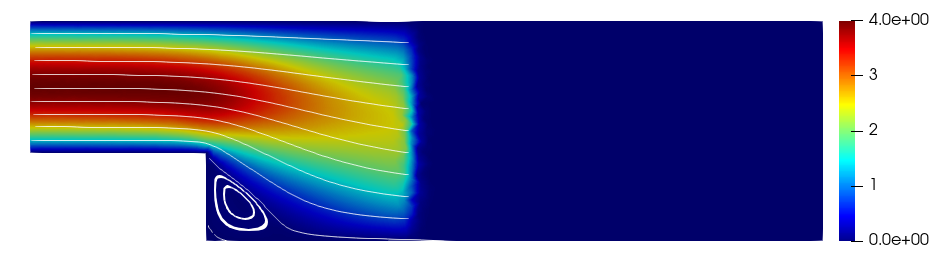}
        \caption{Iteration 0}
         \label{fig:reduced_u_22_0_bfs}
    \end{subfigure}
    \hfill
    \begin{subfigure}[b]{0.49\textwidth}
        \includegraphics[width=\textwidth]{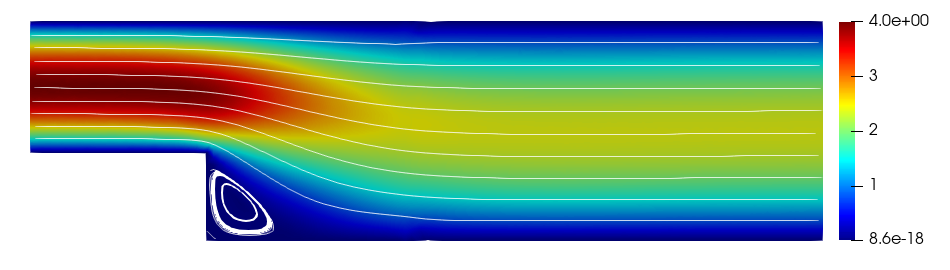}
         \caption{Iteration 5}
         \label{fig:reduced_u_22_5_bfs}
    \end{subfigure}
    \begin{subfigure}[b]{0.49\textwidth}
        \includegraphics[width=\textwidth]{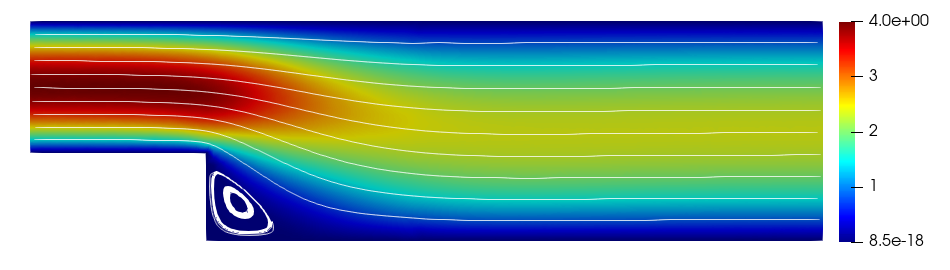}
        \caption{Iteration 10}
         \label{fig:reduced_u_22_10_bfs}
    \end{subfigure}
    \caption{Reduced order solution for the velocities $u_1$ and $u_2$. Values of the parameters $\bar U=4,$ $\nu=0.75$ and \reviewerA{$Re\approx 19$}. \reviewerA{Number of POD modes: 10 - for each state variable, each supremiser and the control, 39 -- for both adjoint velocities}}
    \label{fig:reduced_u_22_bfs}
\end{figure}

\begin{figure}[H]
    \centering
    \begin{subfigure}[b]{0.49\textwidth}
        \includegraphics[width=\textwidth]{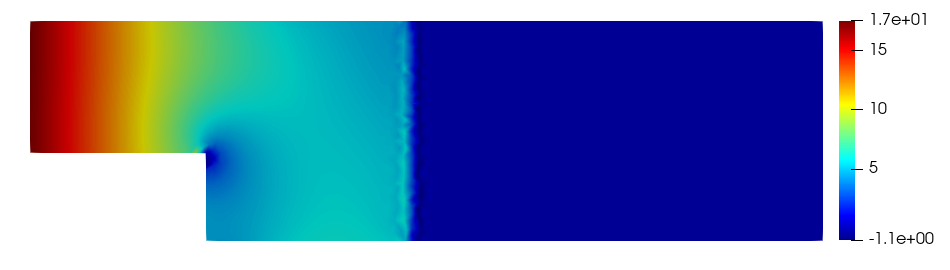}
        \caption{Iteration 0}
         \label{fig:reduced_p_22_0_bfs}
    \end{subfigure}
    \hfill
    \begin{subfigure}[b]{0.49\textwidth}
        \includegraphics[width=\textwidth]{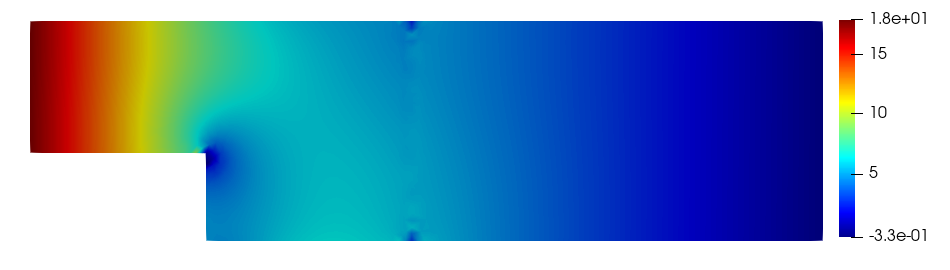}
         \caption{Iteration 5}
         \label{fig:reduced_p_22_5_bfs}

    \end{subfigure}
    
    \begin{subfigure}[b]{0.49\textwidth}
        \includegraphics[width=\textwidth]{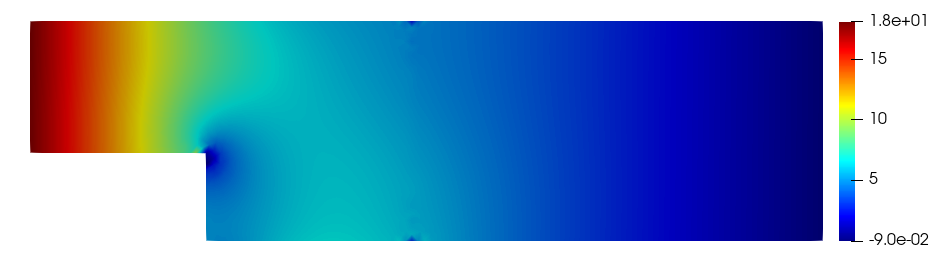}
        \caption{Iteration 10}
         \label{fig:reduced_p_22_10_bfs}
    \end{subfigure}

    \caption{Reduced order solution for the pressures $p_1$ and $p_2$. Values of the parameters $\bar U=4$, $\nu=0.75$ \reviewerA{and $Re\approx 19$}. \reviewerA{Number of POD modes: 10 - for each state variable, each supremiser and the control, 39 -- for both adjoint velocities}}
    \label{fig:reduced_p_22_bfs}
\end{figure}

Figures \ref{fig:reduced_u_11_bfs} -- \ref{fig:reduced_p_22_bfs} \reviewerA{represent} the reduced--order solutions for two different values of the parameters $(\bar U, \nu)=(1,1)$ and \reviewerA{$Re=3$} and $(\bar U, \nu)=(4,0.75)$ \reviewerA{and $Re\approx 19$}. In each of the cases, we choose the following number of the reduced basis functions: $N_{u_1} = N_{s_1} =  N_{p_1} = N_{u_2} = N_{s_2} =  N_{p_2} = N_g = 10$ and $N_{\xi_1} = N_{\xi_2} = 30$. As was previously anticipated, we use a higher number for the adjoint variables $\xi_1$ and $\xi_2$ since they show much slower decay of the singular values (see Figure \ref{fig:singlular_values_bfs}). The solutions were obtained by carrying out 10 optimisation iterations of L--BFGS--B algorithm. Figures \ref{fig:reduced_u_11_bfs} and \ref{fig:reduced_u_22_bfs} show the intermediate solutions at iteration 0, 5 and 10 for the fluid velocities $u_1$ and $u_2$, whereas Figures \ref{fig:reduced_p_11_bfs} and \ref{fig:reduced_p_22_bfs} show the corresponding pressures $p_1$ and $p_2$. The final solution, at the 10th iteration, shows continuity between subdomain solutions at the interface $\Gamma_0$. 

\begin{table}[H]
    \centering
    
    \begin{tabular}{|c|c|c|}
      \hline
        \textbf{Iteration} & \textbf{Functional Value}  &\textbf{Gradient norm }    \\
        \hline  
         0 & 7.869  & 2.120 \\
         5 & 0.107  & 0.401   \\ 
         10 & 0.060 & 0.555  \\ 
         \hline  
    \end{tabular}
    \caption{Functional values and the gradient norm for the \reviewerA{ROM} optimisation solution at parameter values $\bar U=4$, $\nu=0.75$ \reviewerA{and $Re\approx 19$}}
    \label{tab:reduced_2_func_bfs}
\end{table}

We present additional details in Tables \ref{tab:reduced_1_func_bfs} - \ref{tab:reduced_2_errors_bfs}. In particular, in Tables \ref{tab:reduced_1_func_bfs} and \ref{tab:reduced_2_func_bfs}, we list the values for the functional $\mathcal J_\gamma$ and the $L^2(\Gamma_0)$-norm  of the gradient $\frac{d \mathcal J_\gamma}{dg}$ at the different iteration of the optimisation procedure, while Table \ref{tab:reduced_1_errors_bfs} and Table \ref{tab:reduced_2_errors_bfs} contain the absolute and relative errors with respect to the monolithic (entire--domain) solutions $u_h, p_h$, i.e.

\begin{eqnarray*}
    \text{Abs. error} \ u_{N} := || u_{i,N} - u_h ||_{L^2(\Omega_i)}& \text{on domain $\Omega_i$},  \\
    \text{Rel. error} \ u_{N} : = \frac{  || u_{i,N} - u_h ||_{L^2(\Omega_i)}}{|| u_h ||_{L^2(\Omega_i)}}& \text{on domain $\Omega_i$}, \\
     \text{Abs. error} \ p_{N} := || p_{i,N} - p_h ||_{L^2(\Omega_i)}& \text{on domain $\Omega_i$}, \\
    \text{Rel. error} \ p_{N} : = \frac{|| p_{i,N} - p_h ||_{L^2(\Omega_i)}}{|| p_h ||_{L^2(\Omega_i)}} & \text{on domain $\Omega_i$},
\end{eqnarray*}
for $i=1,2.$

\begin{table}[H]
    \centering
    \begin{adjustbox}{max width=\textwidth}
    \begin{tabular}{|c|c|c|c|c|c|c|c|c|}
      \hline
         \textbf{Iteration}  &\multicolumn{2}{|c|}{\textbf{Abs. error $u_{N}$}} &\multicolumn{2}{|c|}{\textbf{Rel. error $u_{N}$}} &\multicolumn{2}{|c|}{\textbf{Abs. error $p_{N}$} }&\multicolumn{2}{|c|}{\textbf{Rel. error $p_{N}$}}   \\ \hline  
         & $\Omega_1$&$\Omega_2$& $\Omega_1$&$\Omega_2$& $\Omega_1$&$\Omega_2$& $\Omega_1$&$\Omega_2$\\ \hline
          0 & 0.1782 & 11.9830 & 0.0128 & 1.0000 & 32.5149 & 21.1630 & 0.6026 & 1.0000 \\
        5 & 0.2826 & 0.8724 & 0.0204 & 0.0728 & 4.1633 & 1.9392 & 0.0772 & 0.0916 \\
        10 & 0.1910 & 0.3826 & 0.0138 & 0.0319 & 0.6725 & 0.7453 & 0.0125 & 0.0352 \\
         \hline  
    \end{tabular}        
    \end{adjustbox}
     \caption{Absolute and relative errors of the \reviewerA{ROM} optimisation solution with respect to the monolithic solution at the parameter values $\bar U=4$, $\nu=0.75$ \reviewerA{and $Re\approx 19$}}
    \label{tab:reduced_2_errors_bfs}
\end{table}

    

    

\begin{table}[H]
    \centering
    \begin{tabular}{|c|c|c|c|c|c|}
      \hline
      \multicolumn{2}{|c|}{\textbf{Parameter value}}  &\multicolumn{2}{|c|}{\textbf{Velocity relative error}} &\multicolumn{2}{|c|}{\textbf{Pressure relative error}}   \\ \hline  
      $\bar U$  & $\nu$ &  $\Omega_1$ & $\Omega_2$ & $\Omega_1$  & $\Omega_2$  \\  \hline  
      $1$ & $1$ & $0.024$ & $0.032$ & $0.005$ & $0.012$ \\ 
       $4$ & $0.75$ & $0.019$ & $0.059$ & $0.021$ & $0.046$ \\ \hline
    \end{tabular}
     \caption{\reviewerB{Relative errors between FOM and ROM solutions (in terms of $H^1$--norm for the velocity fields and $L^2$--norm for the pressure fields)}}
    \label{tab:fom_vs_rom_bfs}
\end{table}

Analysing the results, we are able to see that the reduced basis method gives us a  solution as accurate as the high--fidelity one. The reduced--order approximation of the optimisation problem at hand allowed us to reduce the dimension of the high-fidelity optimisation functional by more than 10 times and enabled us to use 4 times fewer iterations in the optimisation algorithm (each optimisation step requires at least one solve of the state and the adjoint equations). We also note that the fact that we chose a bigger number of the reduced basis functions for the adjoint variables $\xi_1$ and $\xi_2$ is not supposed to affect the computational costs much since the adjoint problem is linear and does not require multiple Newton iteration to be solved so that the biggest computational effort still lies in the nonlinear Navier--Stokes equations and the optimisation process.

\reviewerB{Additionally, in Table \ref{tab:fom_vs_rom_bfs} we provide a comparison between full--order and reduced--order models in terms of the relative errors between ROM solutions with respect to the corresponding FOM solutions. Comparing the convergence results for different models -- monolithic vs. DD--FOM,  monolithic vs. DD--ROM, and  DD--FOM vs. DD--ROM -- it can be seen that the DD--ROM method gives a more accurate solution with respect to DD-FOM. We believe that this is due to the optimisation process: the  DD-ROM is much less sensitive to the initial guess in the optimisation procedure and much fewer iterations are needed for the optimisation algorithm to converge. Nevertheless, errors between DD--FOM and DD--ROM are comparable to the ones with respect to the monolithic solution.}

\reviewerB{
\begin{remark}[High Reynolds and uniqueness of the solution]
    As it is evident from Table~\ref{table:offline_bfs}, the Reynolds number reported for this test case is quite small. This is due to the fact that the optimisation solver diverges for higher Reynolds numbers. The authors suspect that this issue is mostly due to the bifurcation effect (known as the ``Coanda effect'' or ``wall hugging effect'' of these types of simulations). One of the reasons to support this argument is that the range of Reynolds numbers for which the optimisation solver converges changes (though not very significantly) when the interface is moved closer to the beginning or the end of the channel. This problem is very complicated in itself and is addressed, for instance, in \cite{seydel2009practical,kuznetsov1998elements,caloz1997numerical,hess2019localized,pichi2019reduced,pichi2022driving,pintore2021efficient,quaini2016symmetry, panitz1972flow}. In particular, in \cite{pichi2022driving}, it is shown that for a similar test already for $Re \approx 78$ there is non--uniqueness of the solution.
\end{remark}
}

\subsection{Lid-driven cavity flow test case}
\label{cavity_flow}

\begin{figure}

\begin{subfigure}[t]{0.49\textwidth}
    \centering

\begin{tikzpicture}[scale=4.5]

    \draw[draw=red,very thick] (0,1) -- (0.5, 1) node [anchor=south] {$\Gamma_{lid}$}  -- (1,1);
    \draw[draw=blue,very thick] (0,0) --  (0, 0.5) node [anchor=east] {$\Gamma_{wall}$}-- (0,1);
    \draw[draw=blue,very thick](0,0) -- (0.5, 0.) node [anchor=north] {$\Gamma_{wall}$}-- (1,0);
    \draw[draw=blue,very thick] (1,0) -- (1, 0.5) node [anchor=west]{$\Gamma_{wall}$} -- (1,1);

    \fill[fill=blue!10] (0,0) -- (1,0) -- (1,1) -- (0,1) -- (0,0);
    
    \draw (0.5, 0.5) node {$\Omega$};
\end{tikzpicture}
\caption{Physical domain\label{fig:Mono_domain_cavity}}
\captionsetup{justification=centering}
\end{subfigure}
\hfill        
\begin{subfigure}[t]{0.49\textwidth}
\centering    
\begin{tikzpicture}[scale=4.5]

    \fill[fill=blue] (0,0) -- (1,0) -- (1,0.5) -- (0,0.5) -- (0,0);
    
    \draw (0.5, 0.25) node {$\Omega_1$};
    
    \fill[fill=red] (0,0.5) -- (1,0.5) -- (1,1) -- (0,1) -- (0,0.5);
    
    \draw (0.5, 0.75) node {$\Omega_2$};
    
    \draw (0.5, 0) [anchor=north]node {\textcolor{white}{$\Omega_2$}};

\end{tikzpicture}
\caption{Domain splitting\label{fig:dd_domain_cavity}}
\captionsetup{justification=centering}
\end{subfigure}
\caption{Lid-driven cavity flow geometry}
\label{fig:geometry_cavity}
\end{figure}
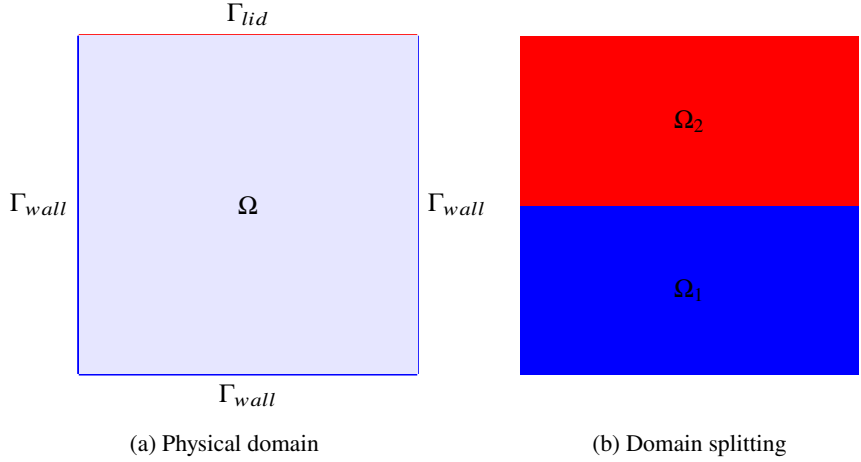

In this section, we provide the numerical simulation for the lid-driven cavity flow test case. Figure \ref{fig:Mono_domain_cavity} represents the physical domain of interest - the unit square. The split into two domains is performed by dissecting the domain by a median horizontal line as shown in Figure \ref{fig:dd_domain_cavity}.

We impose homogeneous Dirichlet boundary conditions on the part of the boundary $\Gamma_{wall}$ for the fluid velocity and the nonzero horizontal constant velocity on the lid boundary $\Gamma_{lid}$: $u_{lid} = \left( \bar U, 0 \right)$; the values of $\bar U$ are reported in Table \ref{table:offline_cavity}.

Two physical parameters are considered: viscosity $\nu$ and the magnitude $\bar U$ of the lid velocity profile $u_{in}$. Details of the offline stage and the finite-element discretisation are summarised in Table \ref{table:offline_cavity}. High-fidelity solutions are obtained by carrying out the minimisation in the space of dimension equal to the number of degrees of freedom at the interface, which is 138 in our test case. The best performance has been achieved by using the limited-memory Broyden–Fletcher–Goldfarb–Shanno (L-BFGS-B) optimisation algorithm, and two stopping criteria are applied: either the maximal number of iteration $ \ It_{max}$ is reached or the gradient norm of the target functional is less than the given tolerance $Tol_{opt}$.

\begin{table}[H]
\begin{center}
\begin{tabular}{ c c }

\hline
    &   \\
    Physical parameters & $2: \nu, \bar U$ \\
    Range $\nu$ & [0.05, 2] \\
    Range $\bar U$ & [0.5, 10] \\
    \reviewerAB{Resulting $Re$ number} & \reviewerAB{[0.25, 200]} \\
    &   \\
    FE velocity order & 2 \\
    FE pressure order & 1 \\
    Total number of FE dofs & 14,867 \\
    Number of FE dofs at the interface & 138 \\
    & \\
    Optimisation algorithm & L-BFGS-B \\
    $It_{max}$ & 100 \\
    $Tol_{opt}$ & $10^{-6}$ \\
    & \\
    $M$ & 300 \\
    $N_{max}$ & 100 \\
    \hline
\end{tabular}
\caption{Computational details of the offline stage. \label{table:offline_cavity}}
\end{center}
\end{table}

 \begin{figure}
    \centering
    \begin{subfigure}[b]{0.49\textwidth}
        \includegraphics[width=\textwidth]{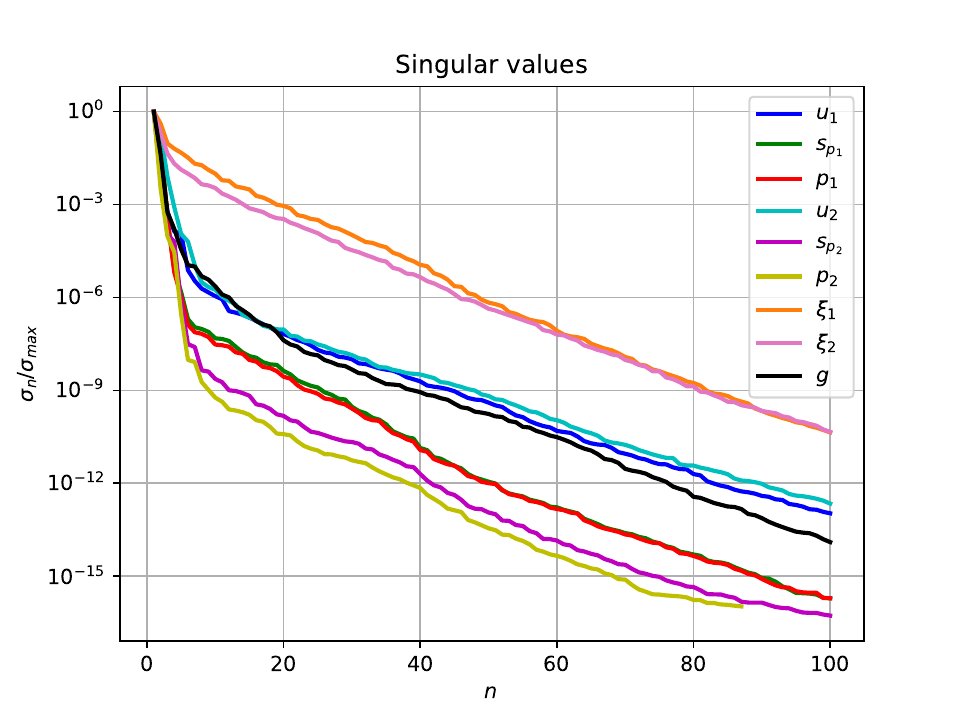}
        \caption{POD singular values as a function of number $n$ of POD modes (log scaling in $y$-direction)}
         \label{fig:singlular_values_cavity}
    \end{subfigure}
    \hfill
    \begin{subfigure}[b]{0.49\textwidth}
        \includegraphics[width=\textwidth]{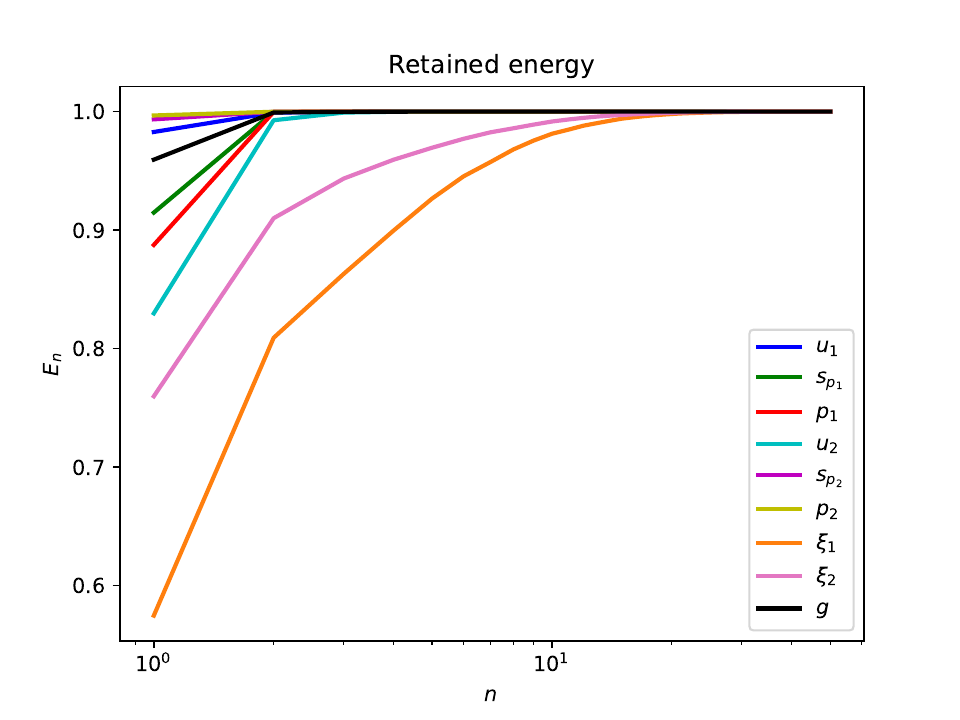}
        \caption{Energy retained by the first $N_{max}$ POD modes (log scaling in $x$-direction)}
         \label{fig:retained_energy_cavity}
    \end{subfigure}
        
    \caption{Results of the offline stage: POD singular eigenvalue decay (a) and retained energy (b) of the first $N_{max}$ POD modes}
    \label{fig:pod_modes_cavity}
\end{figure}

Snapshots are sampled from a training set of $M$ parameters uniformly distributed in the 2-dimensional parameter space, and the first $N_{max}$ POD modes have been retained. Figure~\ref{fig:singlular_values_cavity} shows POD singular values for all the state, the adjoint and the control variables. As it can be seen, the POD singular values corresponding to the adjoint velocities $\xi_1$ and $\xi_2$ feature a slower decay compared to the one for the other variables. In Figure~\ref{fig:retained_energy_cavity}, we can see the behaviour of the energy $E_n$ retained by the first $N$  modes for different components of the solution. Note that, as it was in the previous numerical case, a higher number of modes is needed to correctly represent the adjoint variables $\xi_1$ and $\xi_2$.    

\begin{figure}[H]
    \centering
    \begin{subfigure}[b]{0.32\textwidth}
        \includegraphics[width=\textwidth]{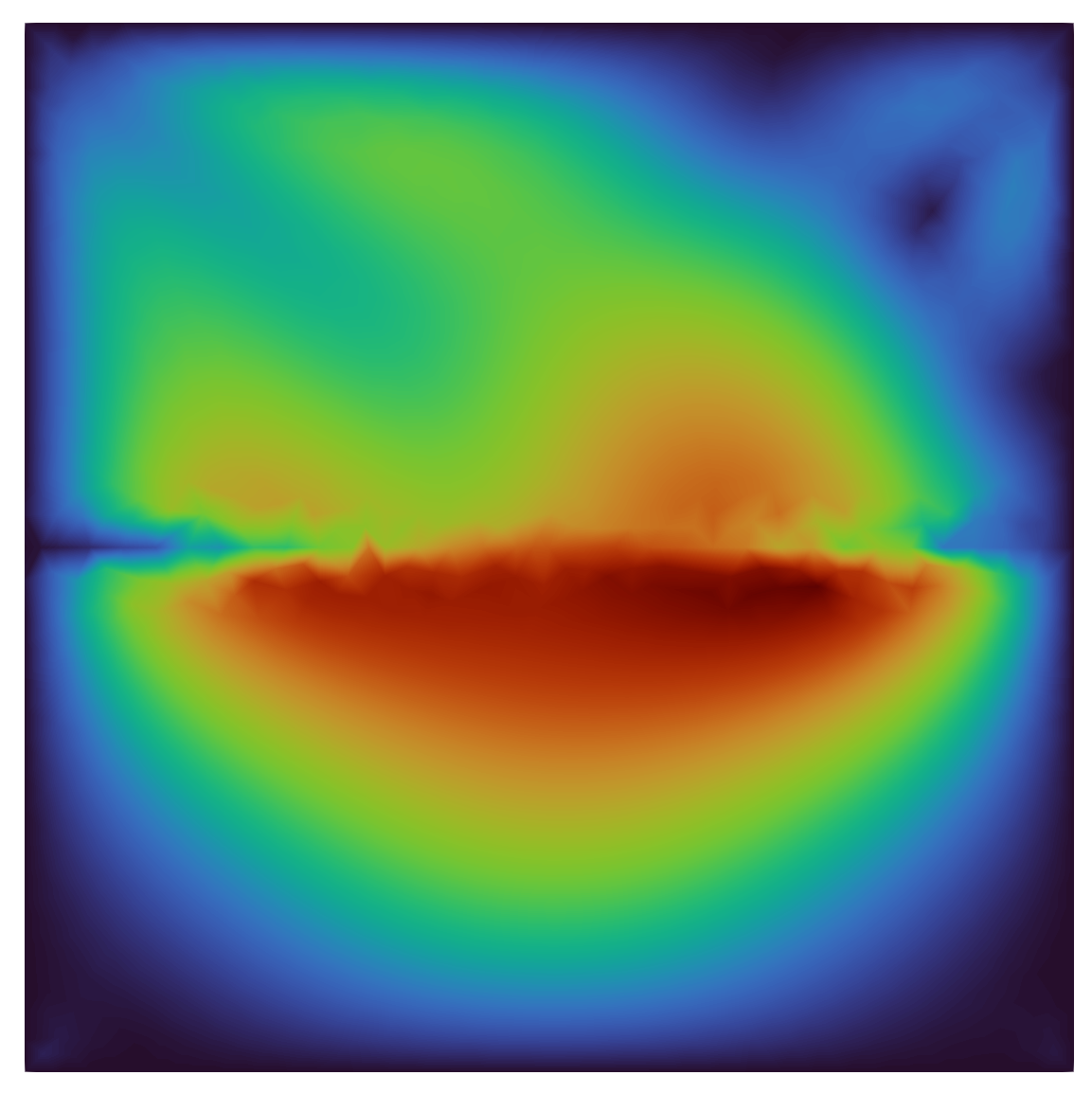}
         \label{fig:mode_u_1_cavity}
    \end{subfigure}
    \hfill
    \begin{subfigure}[b]{0.32\textwidth}
        \includegraphics[width=\textwidth]{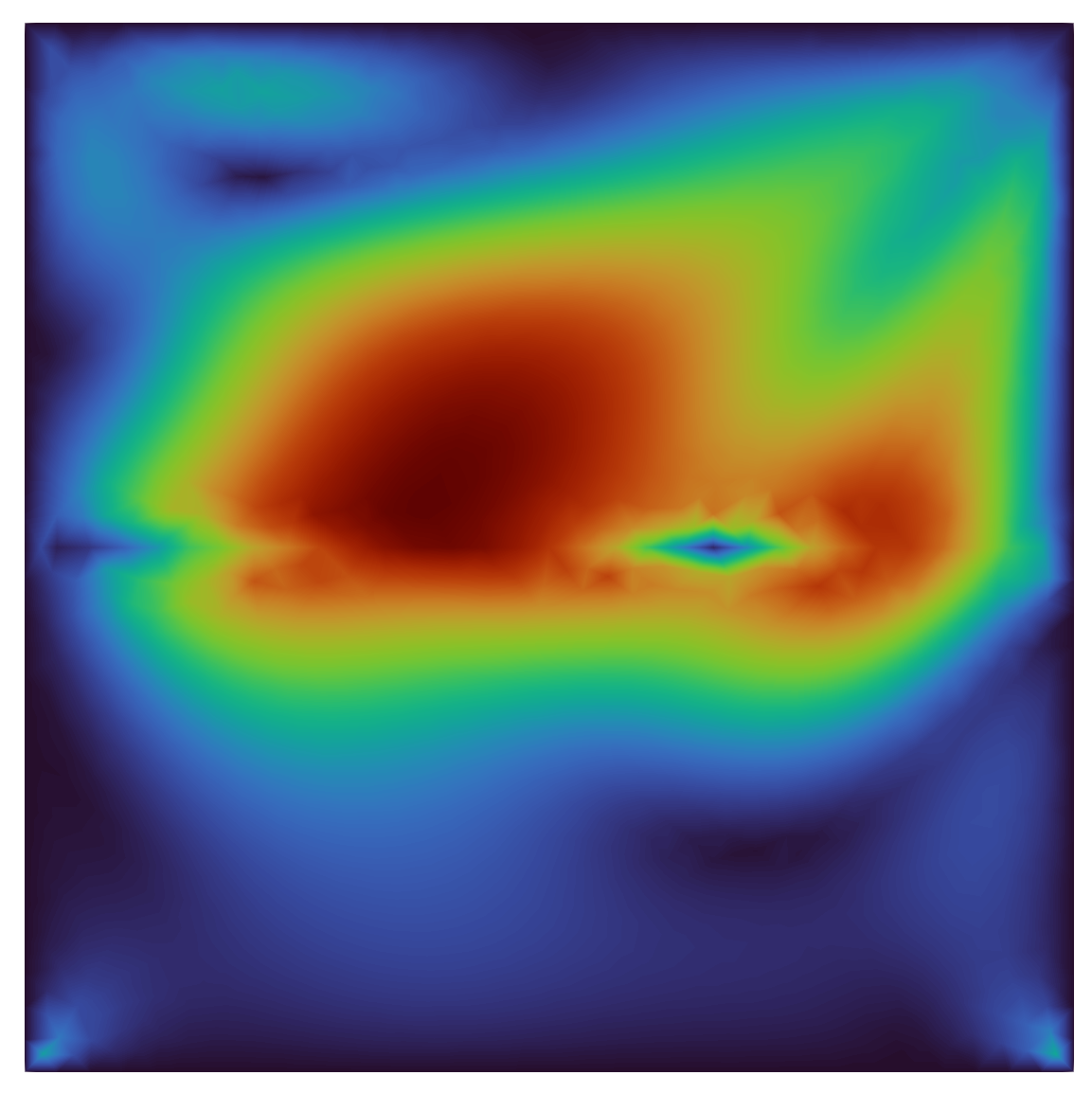}
         \label{fig:mode_u_2_cavity}     
    \end{subfigure}
    \hfill
    \begin{subfigure}[b]{0.32\textwidth}
        \includegraphics[width=\textwidth]{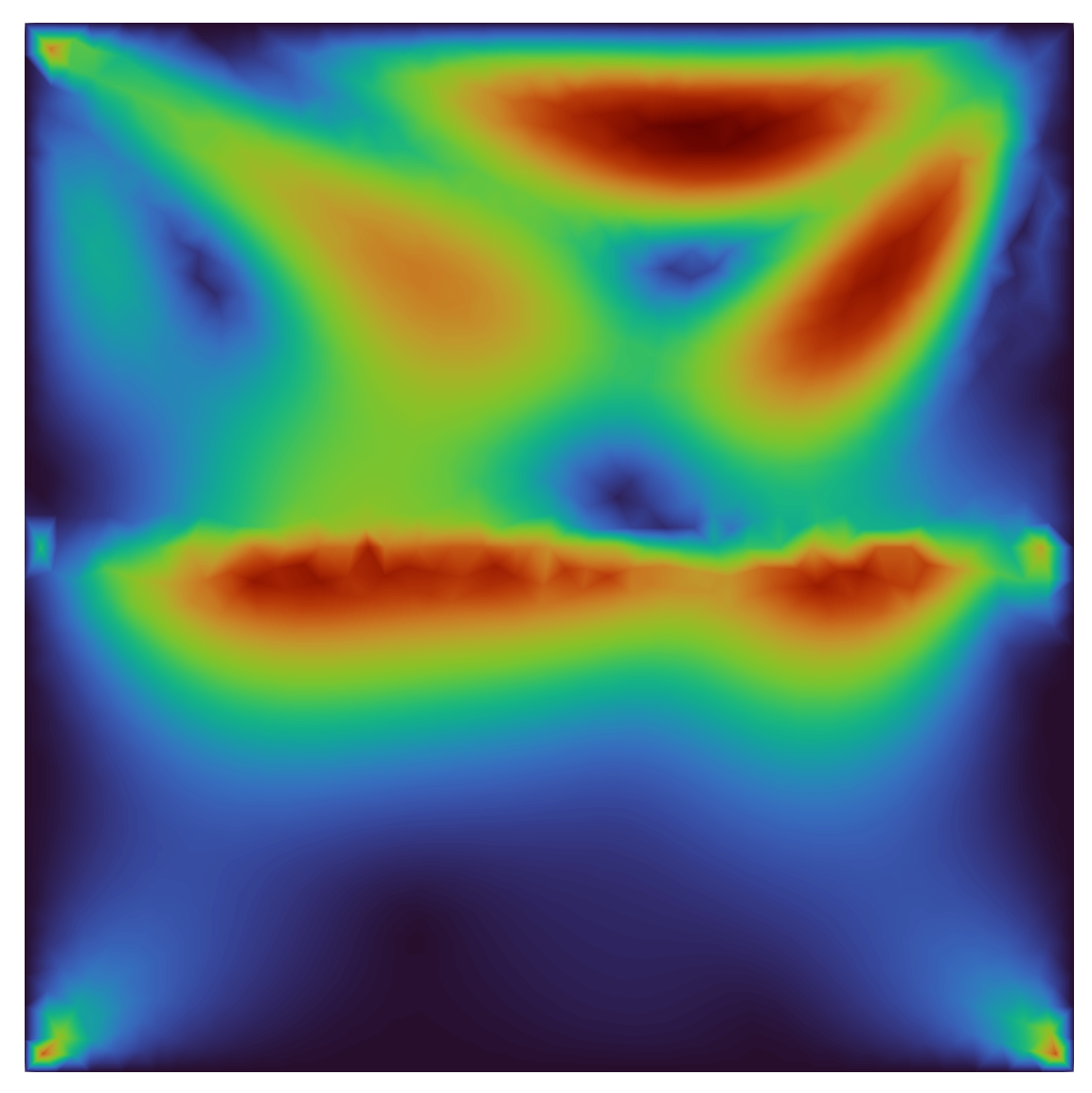}       
         \label{fig:mode_u_3_cavity}
    \end{subfigure}
        
    \caption{The first POD modes for the velocities $u_1$ and $u_2$ (subdomain functions are glued together for visualisation purposes).}
    \label{fig:pod_modes_u_cavity}
\end{figure}

\begin{figure}[H]
    \centering
    \begin{subfigure}[b]{0.32\textwidth}
        \includegraphics[width=\textwidth]{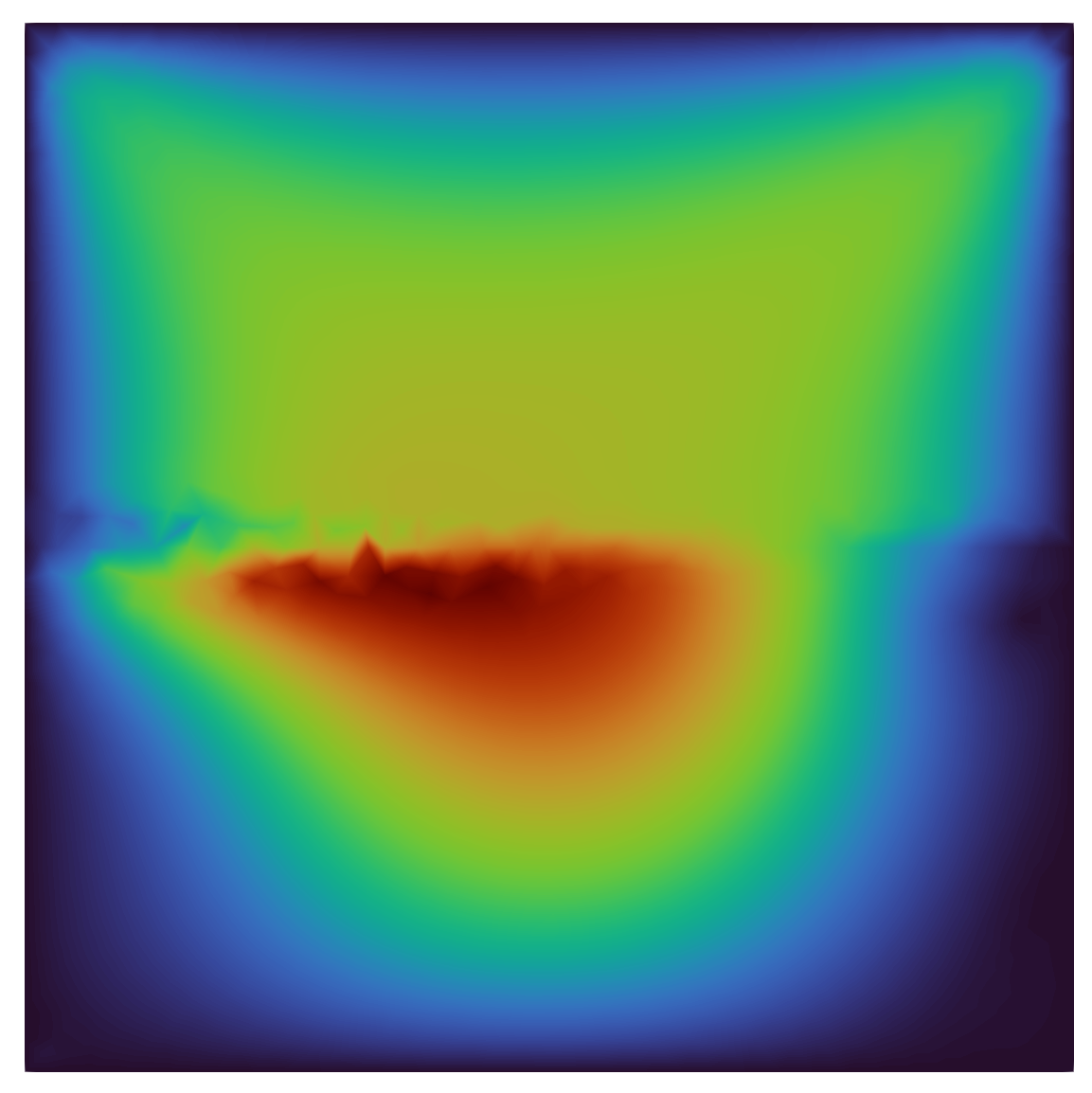}
         \label{fig:mode_sp_1_cavity}
    \end{subfigure}
    \hfill
    \begin{subfigure}[b]{0.32\textwidth}
        \includegraphics[width=\textwidth]{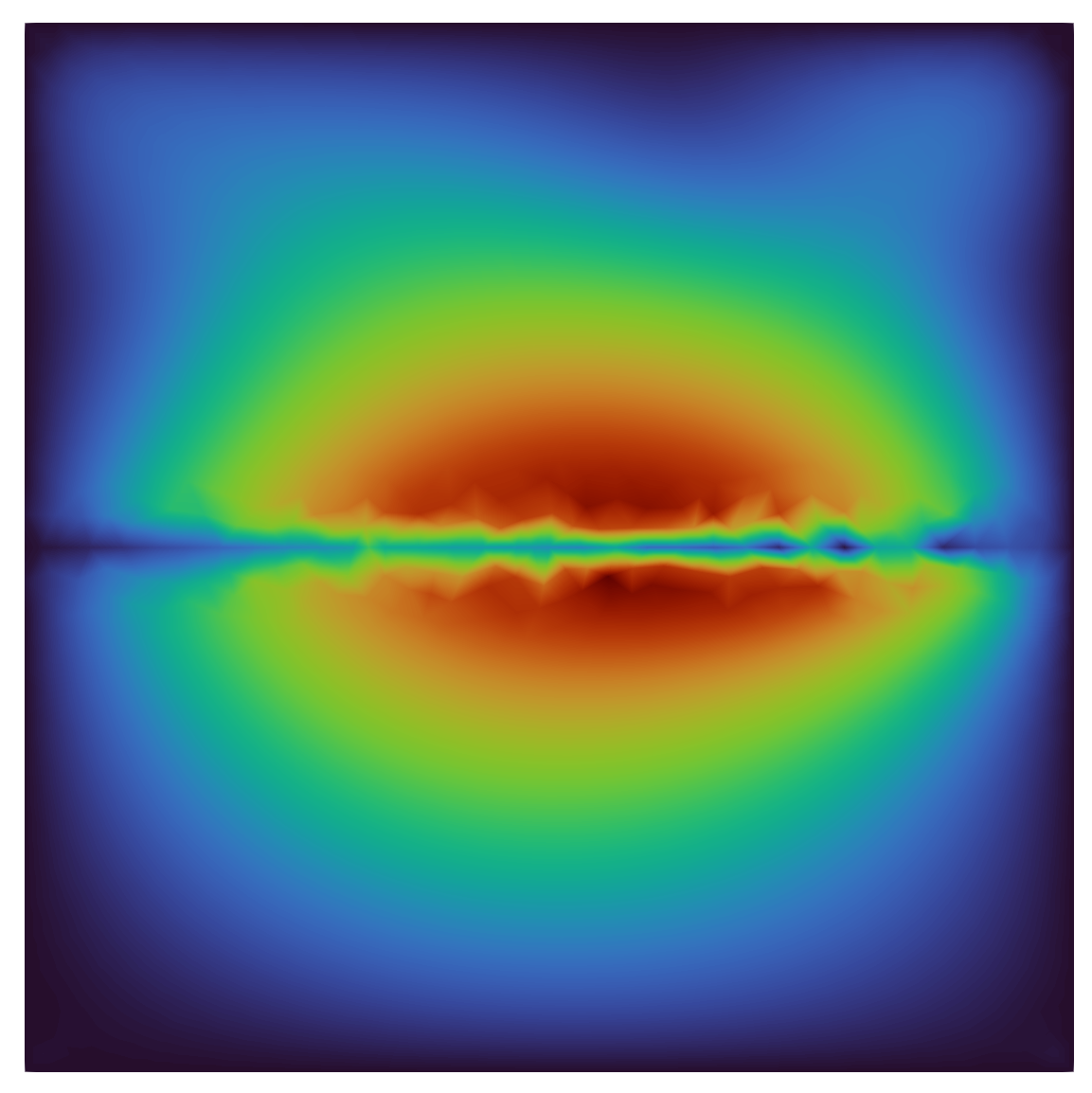}
         \label{fig:mode_sp_2_cavity}
    \end{subfigure}
    \hfill
    \begin{subfigure}[b]{0.32\textwidth}
        \includegraphics[width=\textwidth]{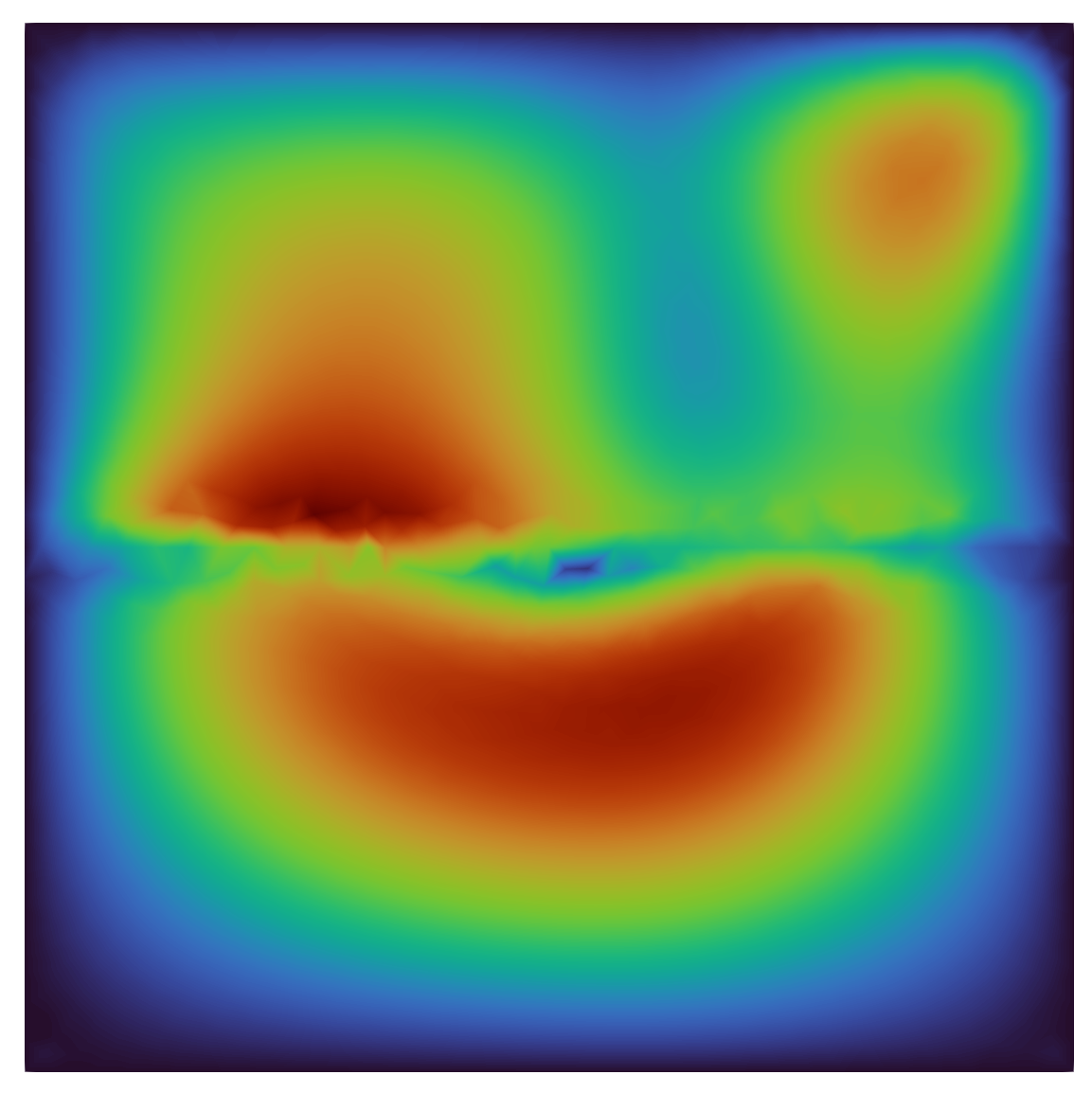}
         \label{fig:mode_sp_3_cavity}
    \end{subfigure}
    \caption{The first POD modes for the supremiser variables $s_1$ and $s_2$ (subdomain functions are glued together for visualisation purposes).}
    \label{fig:pod_modes_sp_cavity}
\end{figure}

\begin{figure}[H]
    \centering
    \begin{subfigure}[b]{0.32\textwidth}
        \includegraphics[width=\textwidth]{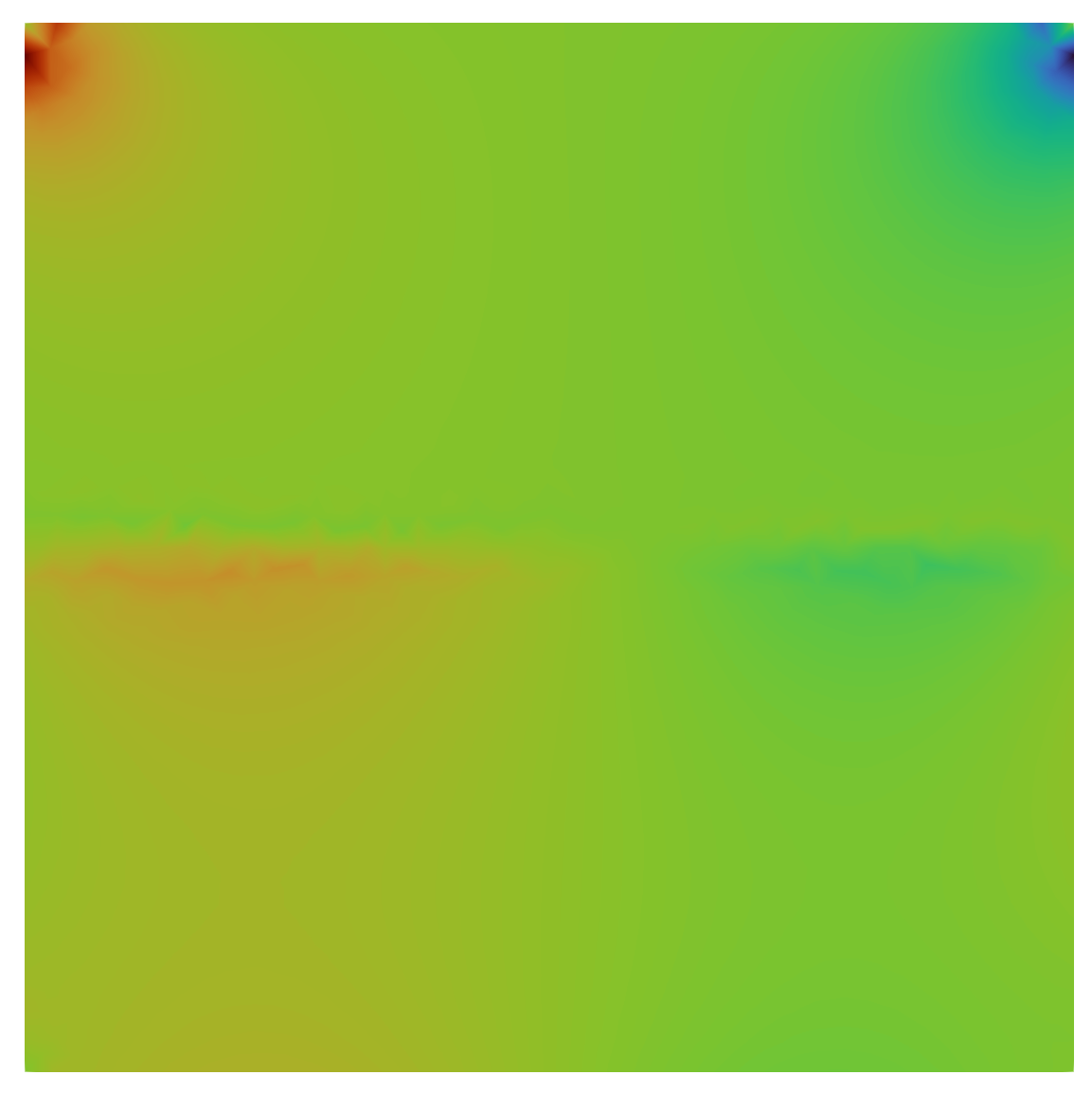}
         \label{fig:mode_p_1_cavity}
    \end{subfigure}
    \hfill
    \begin{subfigure}[b]{0.32\textwidth}
        \includegraphics[width=\textwidth]{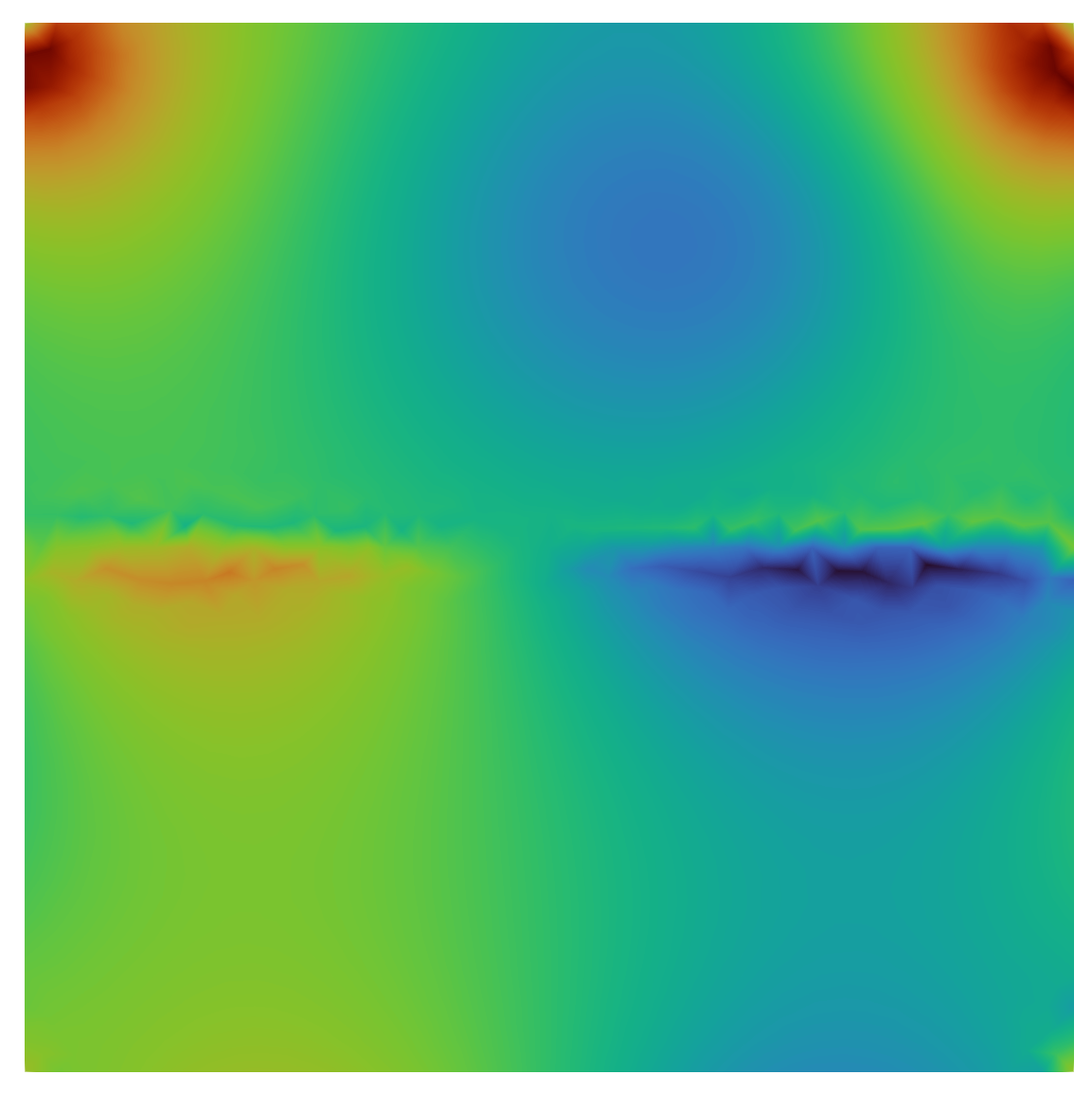}
         \label{fig:mode_p_2_cavity}     
    \end{subfigure}
    \hfill
    \begin{subfigure}[b]{0.32\textwidth}
        \includegraphics[width=\textwidth]{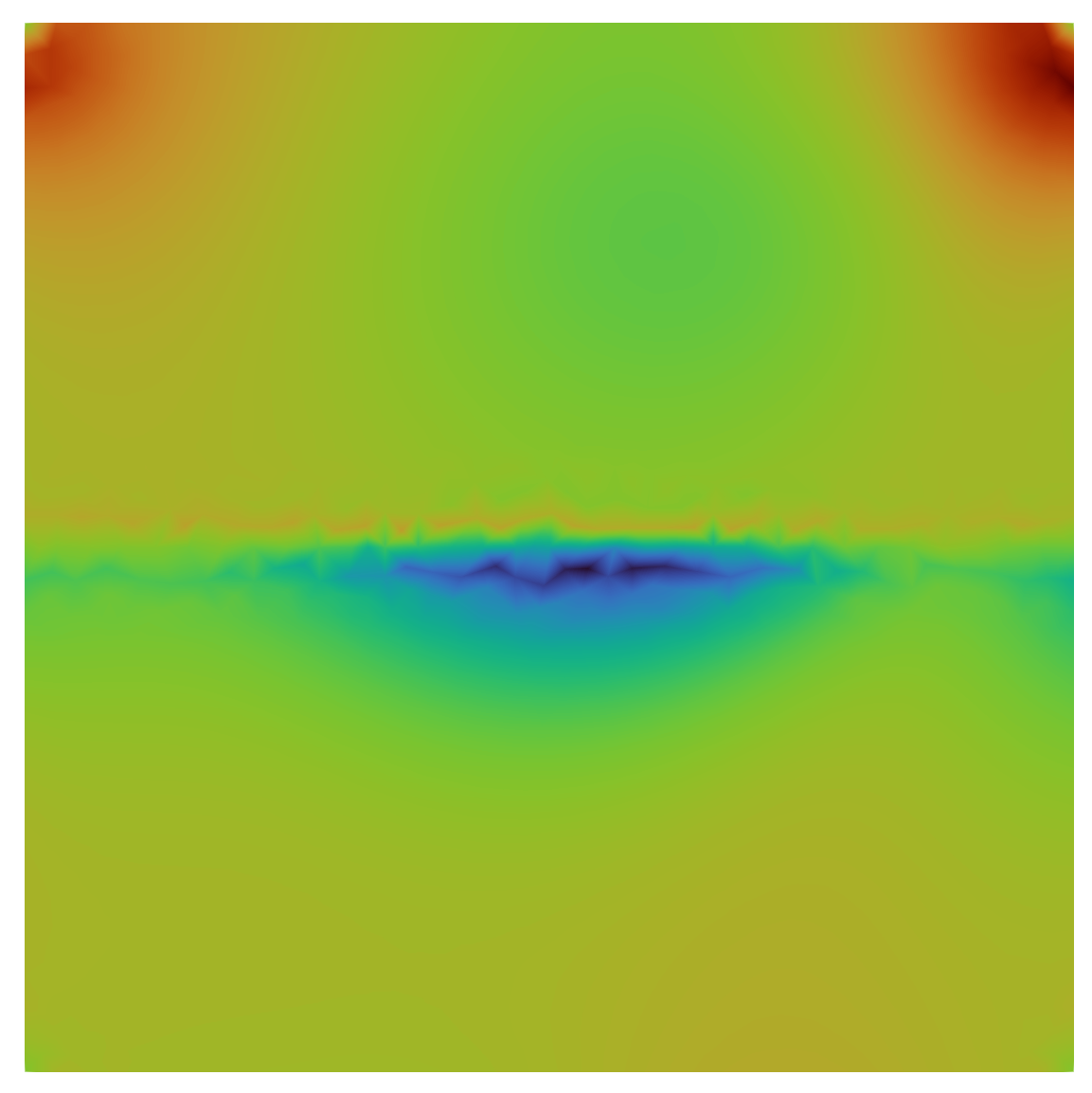}       
         \label{fig:mode_p_3_cavity}
    \end{subfigure}
        
    \caption{The first POD modes for the pressures $p_1$ and $p_2$ (subdomain functions are glued together for visualisation purposes).}
    \label{fig:pod_modes_p_cavity}
\end{figure}

\begin{figure}[H]
    \centering
    \begin{subfigure}[b]{0.32\textwidth}
        \includegraphics[width=\textwidth]{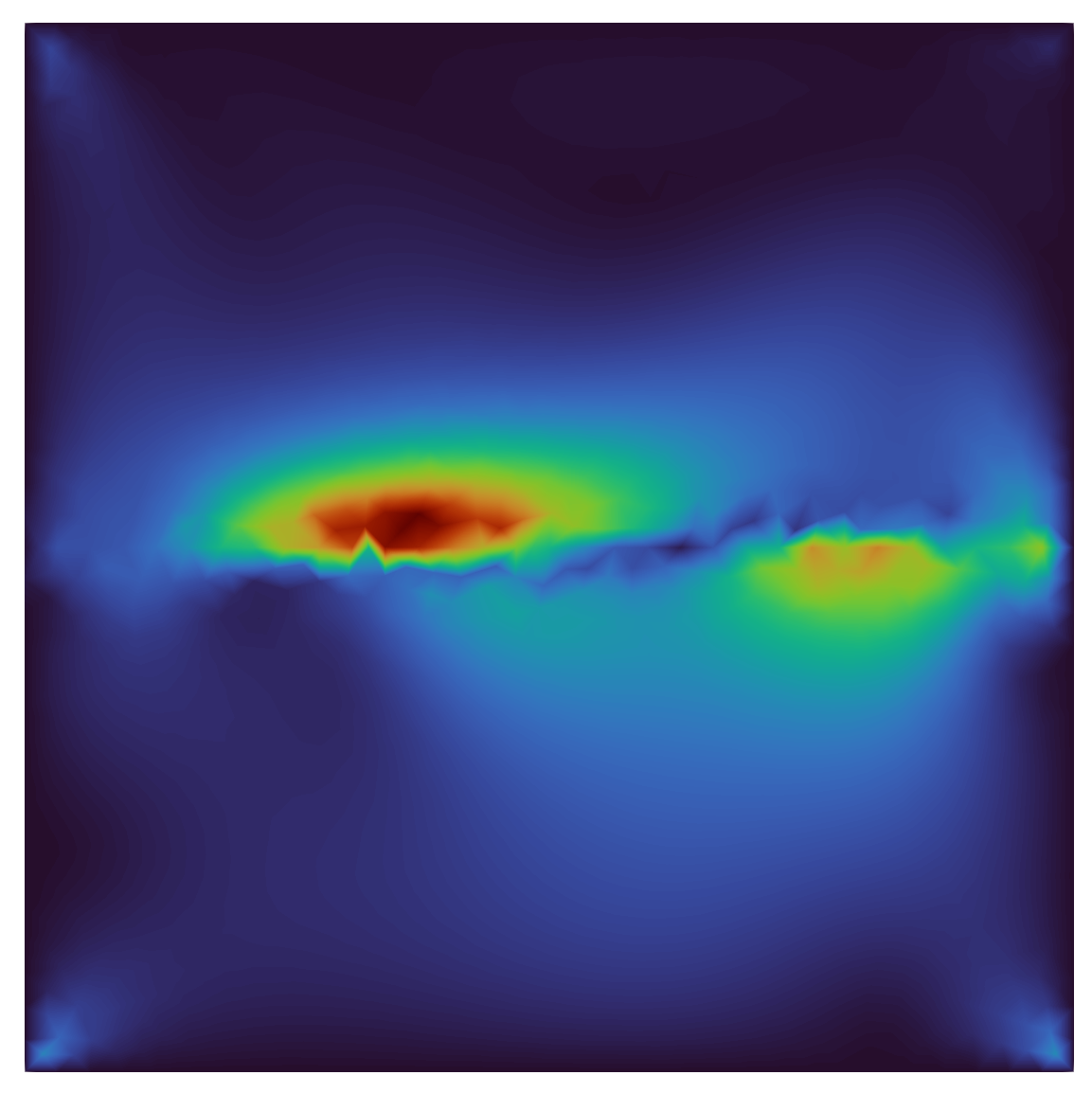}
         \label{fig:mode_x_1_cavity}
    \end{subfigure}
    \hfill
    \begin{subfigure}[b]{0.32\textwidth}
        \includegraphics[width=\textwidth]{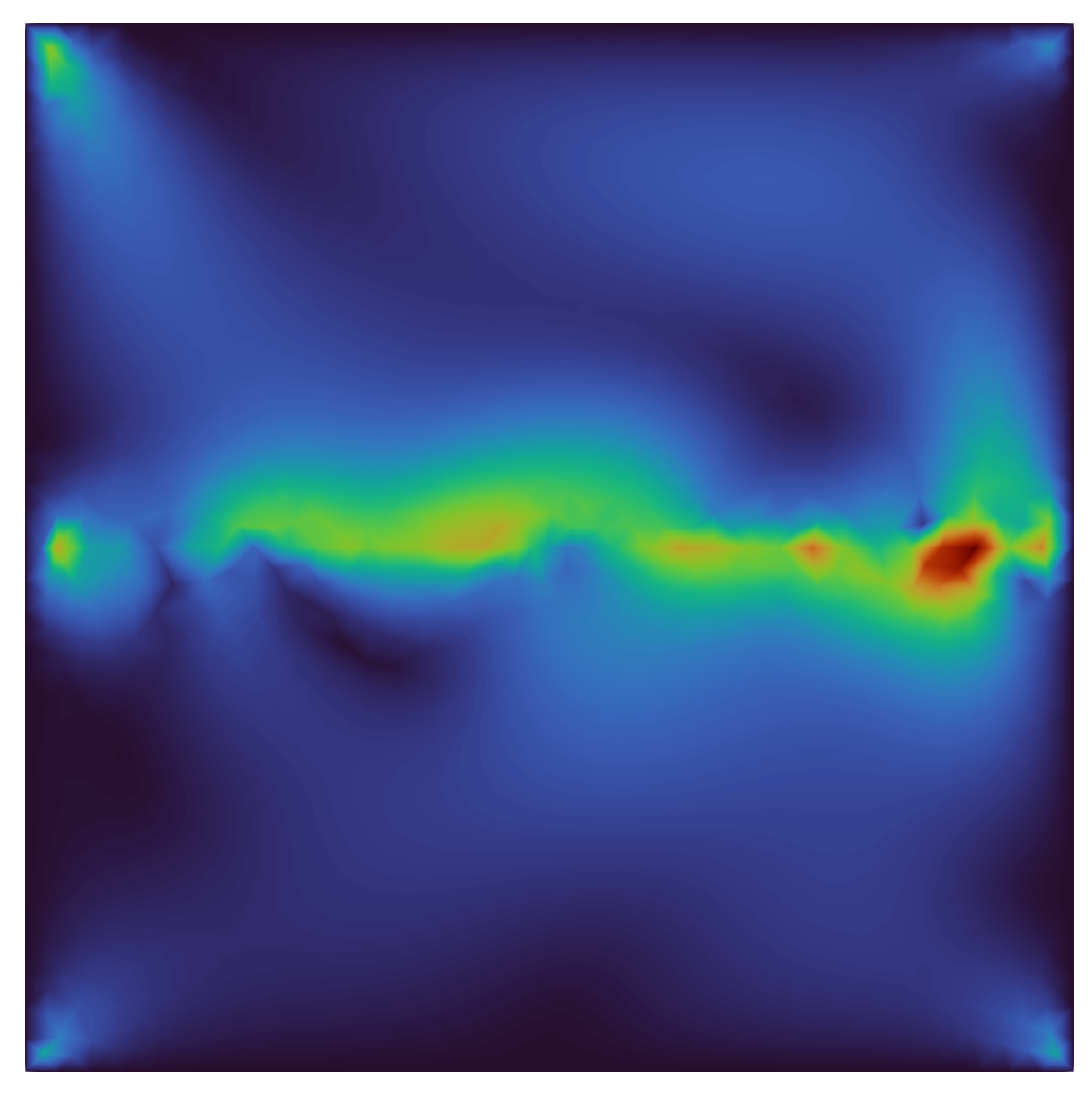}
         \label{fig:mode_x_2_cavity}
    \end{subfigure}
    \hfill
    \begin{subfigure}[b]{0.32\textwidth}
        \includegraphics[width=\textwidth]{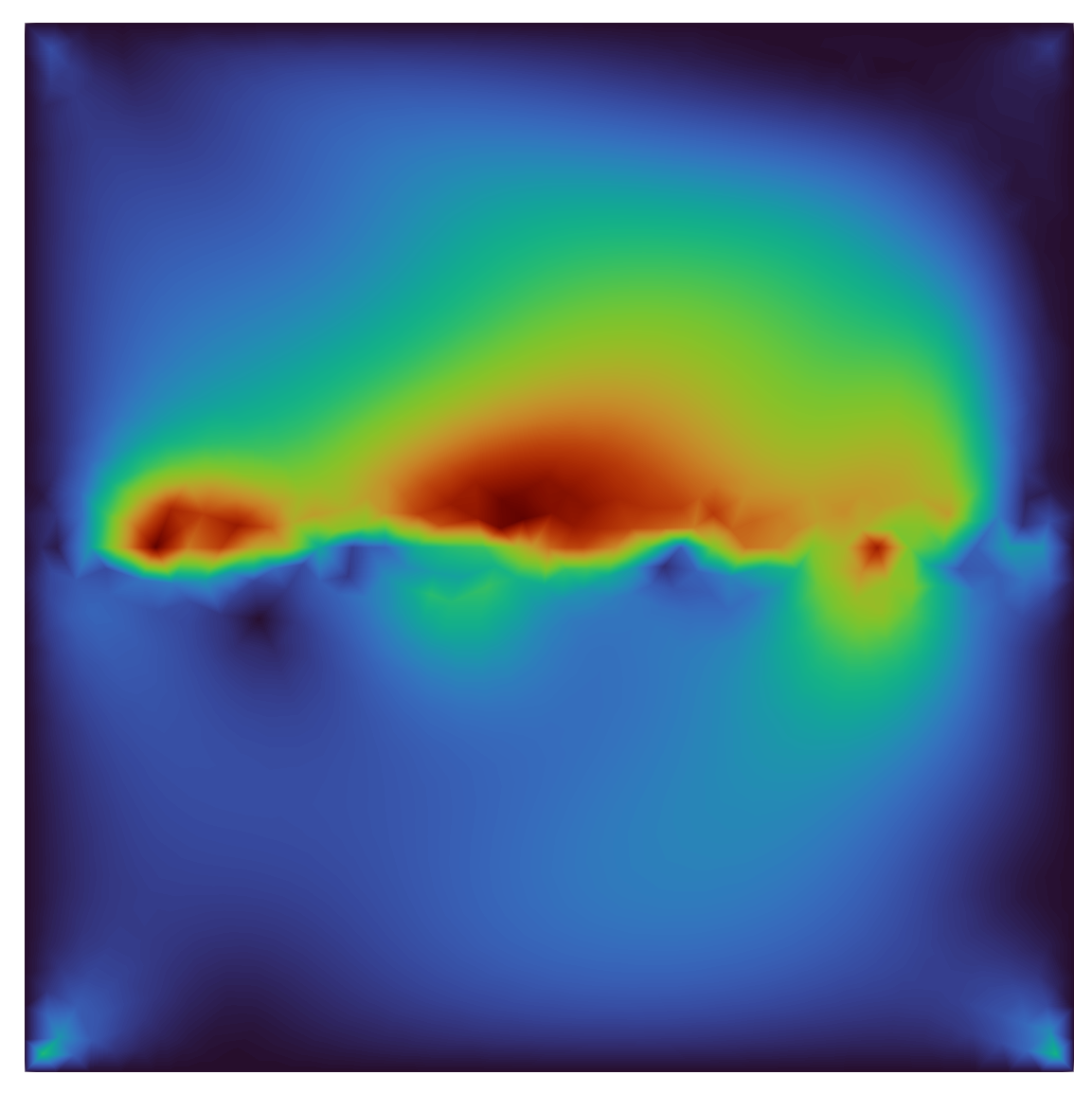}
         \label{fig:mode_x_3_cavity}
    \end{subfigure}
    \caption{The first POD modes for the adjoint velocities $\xi_1$ and $\xi_2$ (subdomain functions are glued together for visualisation purposes).}
    \label{fig:pod_modes_x_cavity}
\end{figure}

Figures~\ref{fig:pod_modes_u_cavity}--\ref{fig:pod_modes_x_cavity} represent first three POD modes for the variables $u_1, u_2, s_1, s_2$, $p_1$, $p_2$ and $\xi_1,\xi_2$. We stress that the POD modes were obtained separately for each component and the resulting figures are obtained by gluing the subdomain functions just for the sake of visualisation. 

Figure~\ref{fig:pod_modes_u_cavity} shows the first modes for the fluid velocities $u_1$ and $u_2$. In particular, we notice that the modes corresponding to $u_2$ (on the upper section of the domain) are zero at the lid boundary due to the use of lifting function. Figure~\ref{fig:pod_modes_x_cavity} shows the first three modes for the adjoint variables  $\xi_1$ and $\xi_2$: note that they are concentrated only around the interface $\Gamma_0$ because the only nonzero contribution in the adjoint equations is coming from the source terms, which are defined solely on the interface $\Gamma_0$.

\begin{figure}[H]
    \centering
    \begin{subfigure}[b]{0.32\textwidth}
        \includegraphics[width=\textwidth]{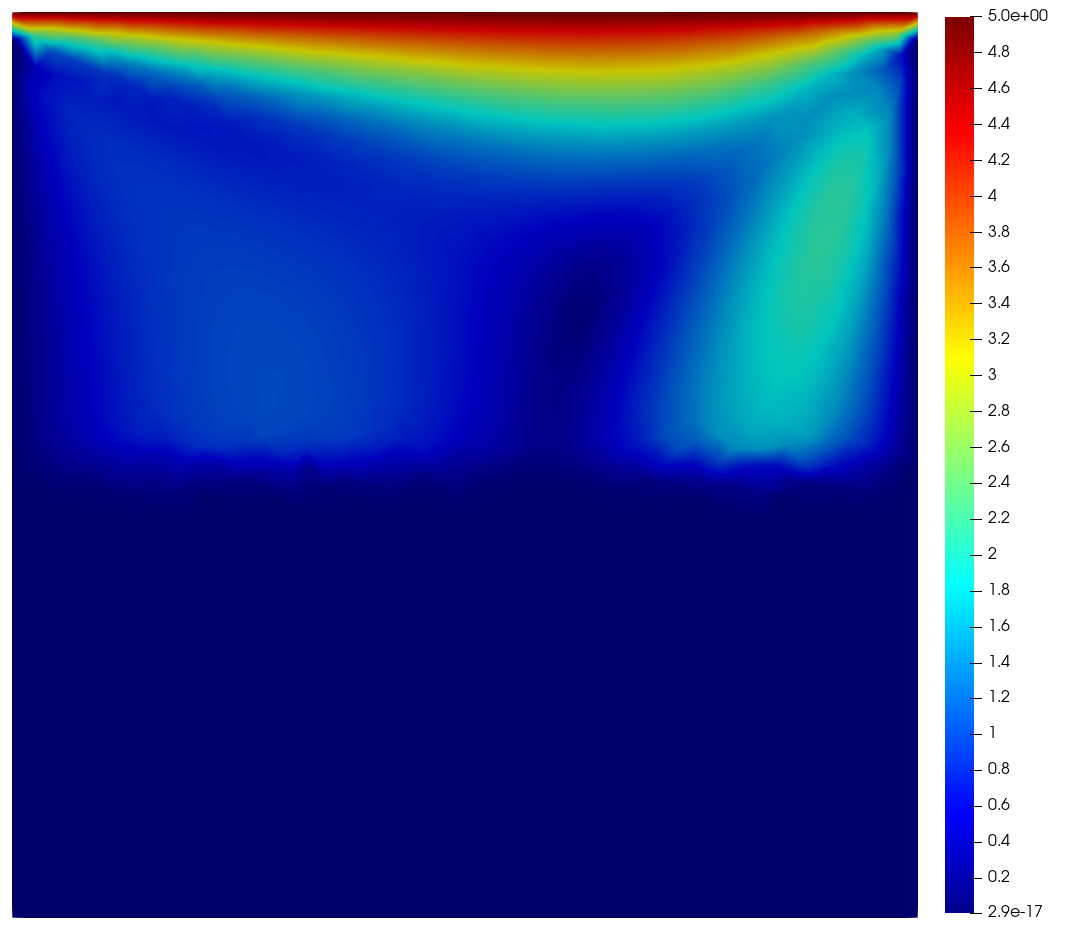}
        \caption{Iteration 0}
         \label{fig:truth_u_5005_0_cavity}
    \end{subfigure}
    \hfill
    \begin{subfigure}[b]{0.32\textwidth}
        \includegraphics[width=\textwidth]{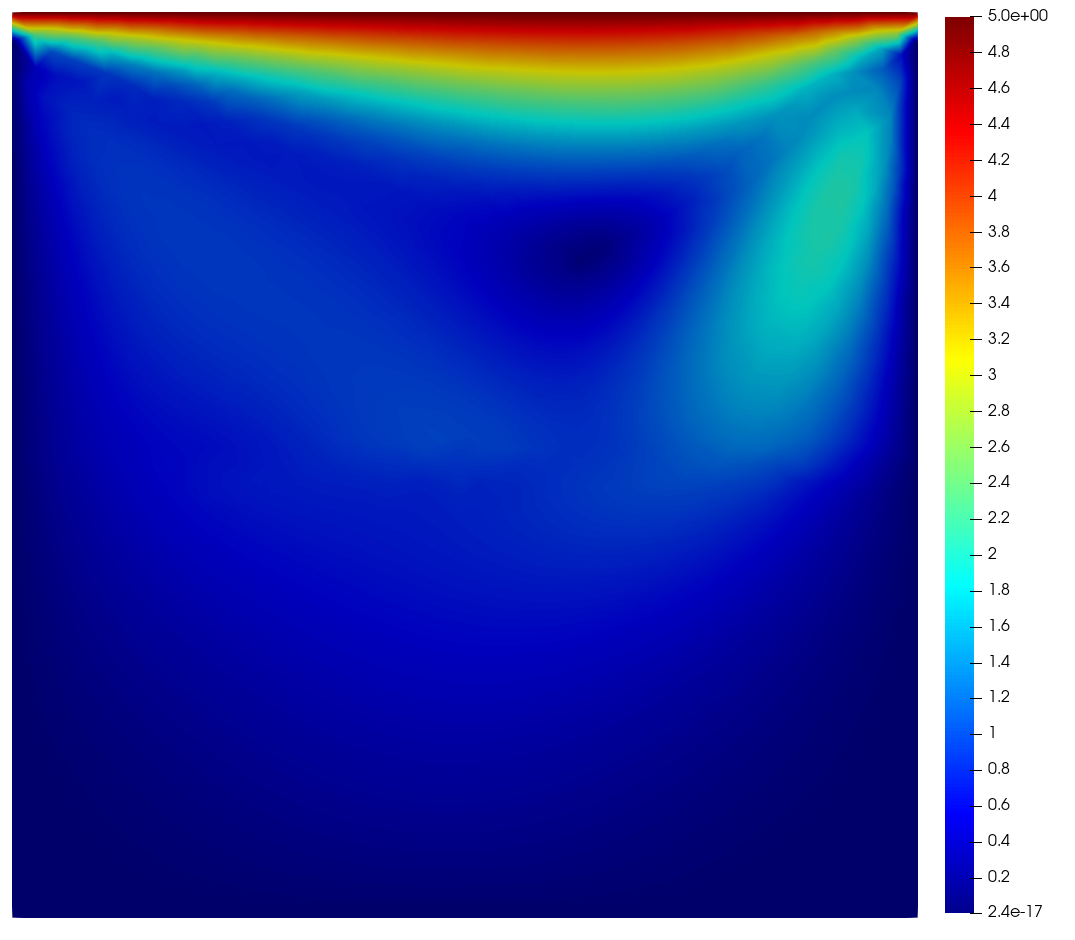}
         \caption{Iteration 5}
         \label{fig:truth_u_5005_5_cavity}

    \end{subfigure}
    \hfill
    \begin{subfigure}[b]{0.32\textwidth}
        \includegraphics[width=\textwidth]{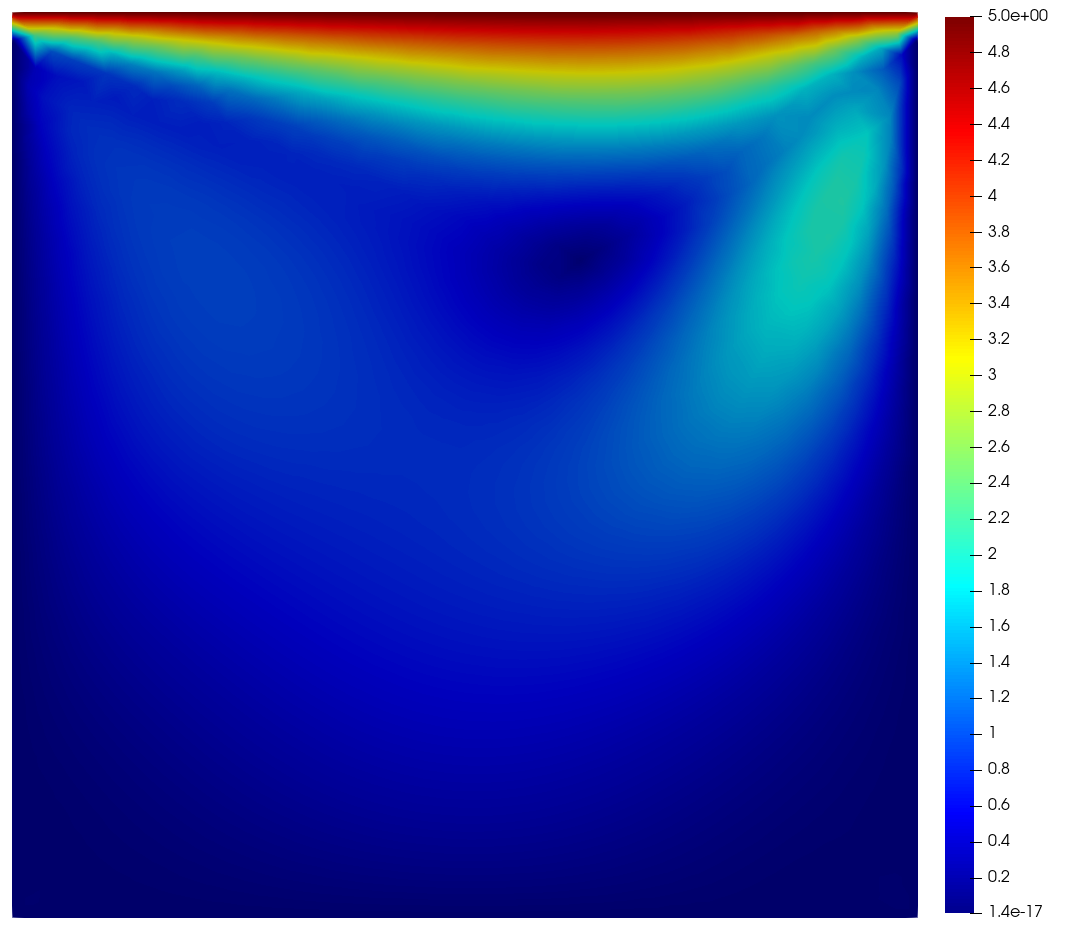}
         \caption{Iteration 25}
         \label{fig:truth_u_5005_25_cavity}

    \end{subfigure}
        
    \caption{High--fidelity solution for the velocities $u_1$ and $u_2$. Values of the parameters $\bar U=5$ and $\nu=0.05$}
    \label{fig:truth_u_5005_cavity}
\end{figure}

\begin{figure}[H]
    \centering
    \begin{subfigure}[b]{0.32\textwidth}
        \includegraphics[width=\textwidth]{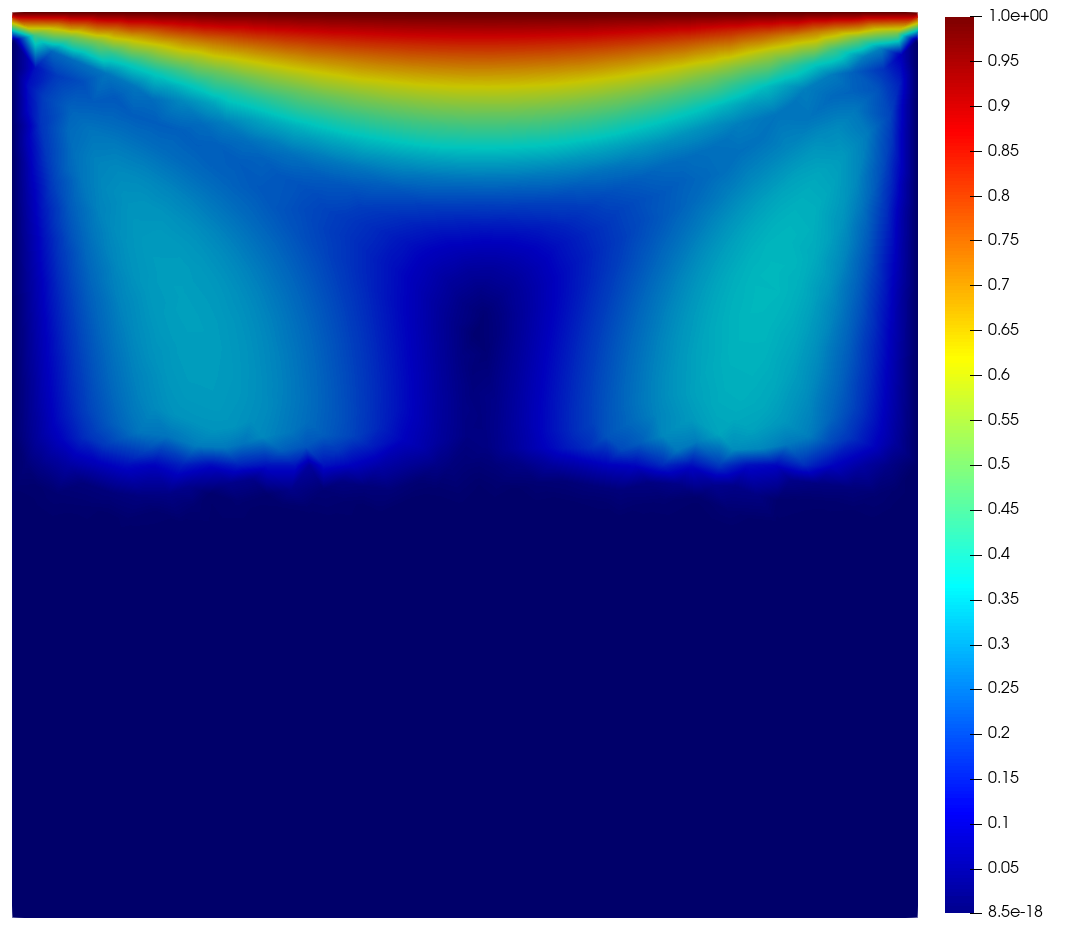}
        \caption{Iteration 0}
         \label{fig:truth_u_101_0_cavity}
    \end{subfigure}
    \hfill
    \begin{subfigure}[b]{0.32\textwidth}
        \includegraphics[width=\textwidth]{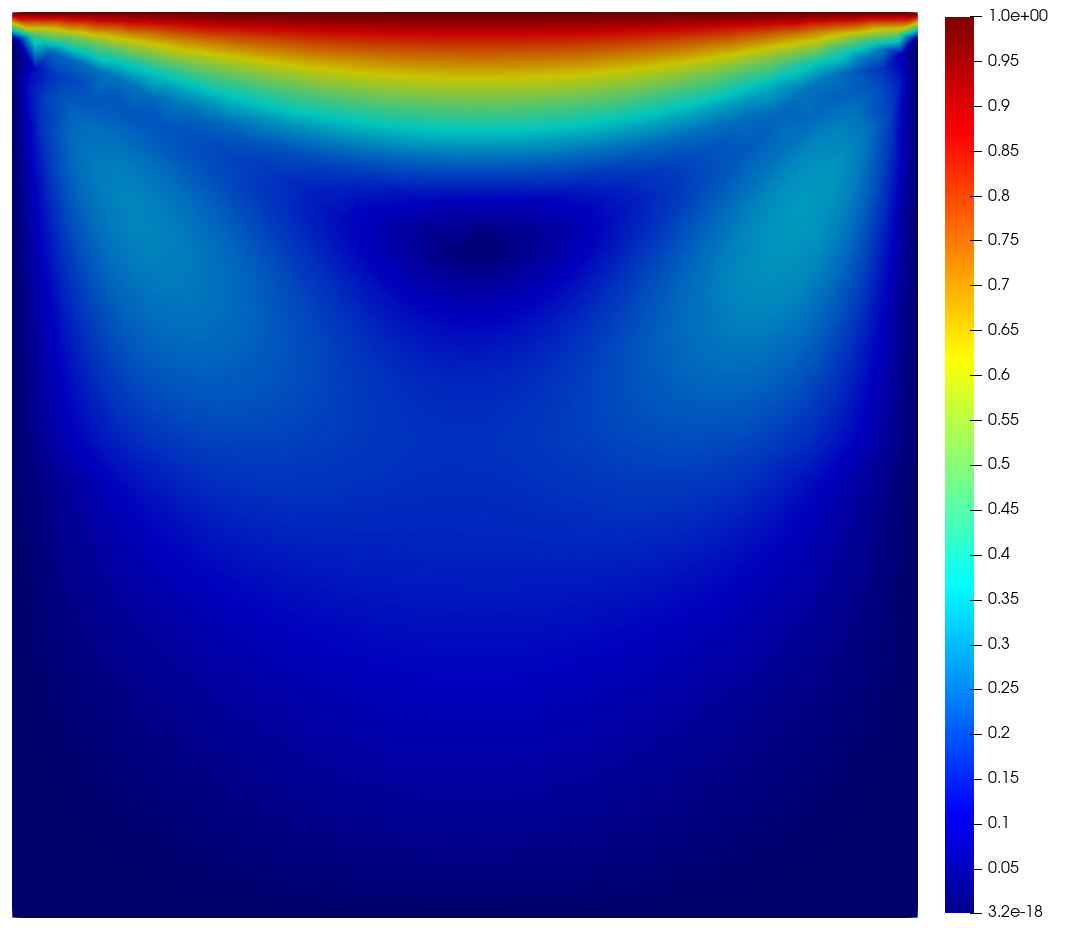}
         \caption{Iteration 5}
         \label{fig:truth_u_101_5_cavity}
    \end{subfigure}
    \hfill
    \begin{subfigure}[b]{0.32\textwidth}
        \includegraphics[width=\textwidth]{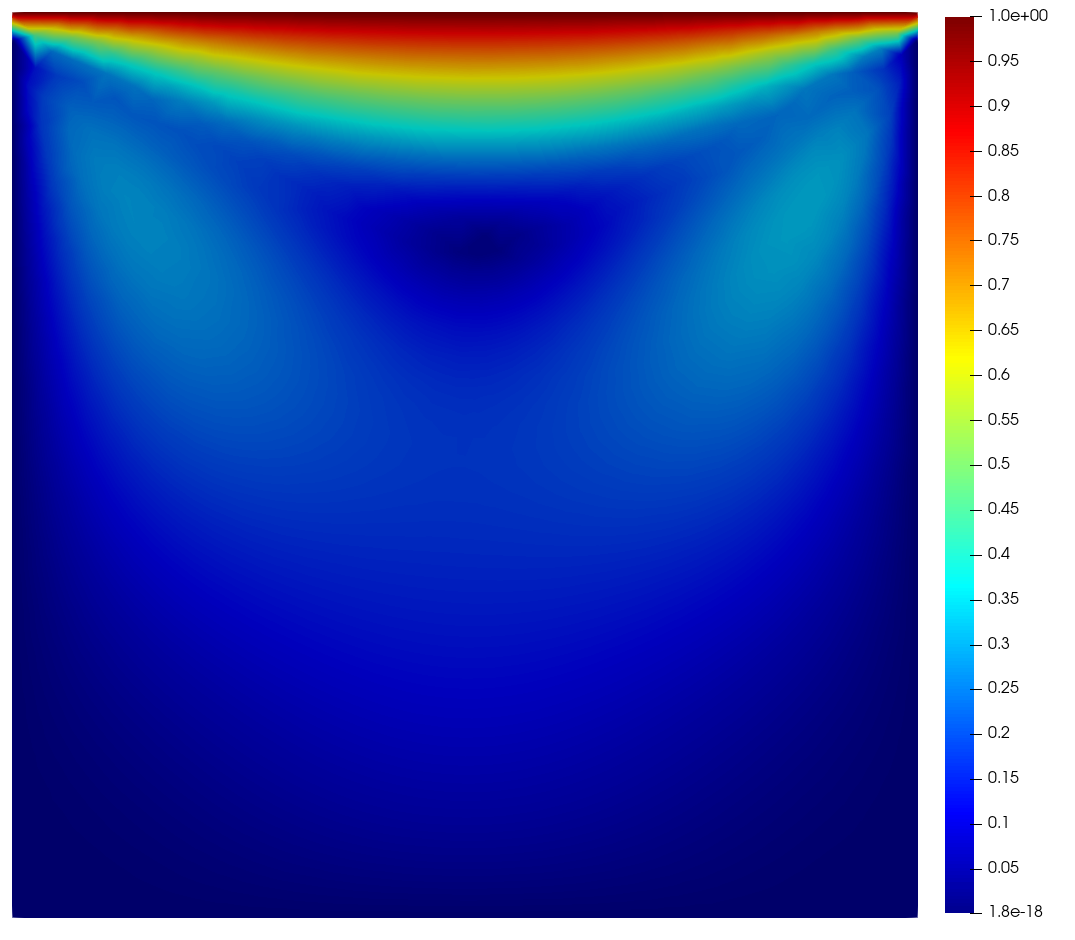}
         \caption{Iteration 25}
         \label{fig:truth_u_101_25_cavity}
    \end{subfigure}
        
    \caption{High--fidelity \reviewerA{FOM} solution for the velocities $u_1$ and $u_2$. Values of the parameters $\bar U=1$, $\nu=0.1$ \reviewerA{and with $Re=10$}}
    \label{fig:truth_u_101_cavity}
\end{figure}

Figures \ref{fig:truth_u_5005_cavity} and \ref{fig:truth_u_101_cavity} represent the high--fidelity solutions for two different values of the parameters $(\bar U, \nu)=(5,0.05)$, \reviewerA{with $Re=100$,} and $(\bar U, \nu)=(1,0.1)$, \reviewerA{with $Re=10$}. The solutions were obtained by carrying out 25 optimisation iterations via L--BFGS--B algorithm. Figures \ref{fig:truth_u_5005_cavity} and \ref{fig:truth_u_101_cavity} show the intermediate solutions at iteration 0, 5 and 25 for the fluid velocities $u_1$ and $u_2$. The final solution is taken to be the 25-iteration optimisation solution as we can observe a continuity between subdomain solutions at the interface $\Gamma_0$. We present additional details in Tables \ref{tab:truth_5005_func_cavity} - \ref{tab:truth_101_errors_cavity}. In particular, in Tables \ref{tab:truth_5005_func_cavity} and \ref{tab:truth_101_func_cavity}, we list the values for the functional $\mathcal J_\gamma$ and the $L^2(\Gamma_0)$-norm of the gradient $\frac{d \mathcal J_\gamma}{dg}$ at the different iteration of the optimisation procedure, while Table \ref{tab:truth_5005_errors_cavity} and Table \ref{tab:truth_101_errors_cavity} contain the absolute and relative errors with respect to the monolithic(entire--domain) solutions $u_h, p_h$.

\begin{table}[H]
    \centering
    
    \begin{tabular}{|c|c|c|}
      \hline
        \textbf{Iteration} & \textbf{Functional Value}  &\textbf{Gradient norm }    \\
        \hline  
        0 & $4.4 \cdot 10^{-1}$ & $3.398 $ \\
         5 & $3.0 \cdot 10^{-2}$ & $1.001 $  \\ 
         10 & $3.5 \cdot 10^{-3}$ & $0.171 $  \\ 
         25 & $8.7 \cdot 10^{-5}$ & $0.016 $  \\
         \hline  
    \end{tabular}
    \caption{Functional values and the gradient norm for the \reviewerA{FOM} optimisation solution at parameter values $\bar U=5$, $\nu=0.05$ \reviewerA{and with $Re=100$}}
    \label{tab:truth_5005_func_cavity}
\end{table}

\begin{table}[H]
    \centering
   
        \begin{adjustbox}{max width=\textwidth}
    \begin{tabular}{|c|c|c|c|c|c|c|c|c|}
      \hline
         \textbf{Iteration}  &\multicolumn{2}{|c|}{\textbf{Abs. error $u_{N}$}} &\multicolumn{2}{|c|}{\textbf{Rel. error $u_{N}$}} &\multicolumn{2}{|c|}{\textbf{Abs. error $p_{N}$} }&\multicolumn{2}{|c|}{\textbf{Rel. error $p_{N}$}}   \\ \hline  
         & $\Omega_1$&$\Omega_2$& $\Omega_1$&$\Omega_2$& $\Omega_1$&$\Omega_2$& $\Omega_1$&$\Omega_2$\\ \hline
        0 & 0.3411 & 0.1949 & 1.0000 & 0.1653 & 0.2689 & 0.3149 & 1.0000 & 0.2330 \\
        5 & 0.0623 & 0.0613 & 0.1826 & 0.0520 & 0.0531 & 0.0575 & 0.3633 & 0.0426 \\
        10 & 0.0114 & 0.0136 & 0.0334 & 0.0116 & 0.0184 & 0.0206 & 0.1256 & 0.0153 \\
        25 & 0.0051 & 0.0062 & 0.0151 & 0.0053 & 0.0143 & 0.0147 & 0.0980 & 0.0109 \\   
         \hline  
    \end{tabular}
    \end{adjustbox}
    \caption{Absolute and relative errors of the \reviewerA{FOM} optimisation solution with respect to the monolithic solution at the parameter value $\bar U=5$, $\nu=0.05$ \reviewerA{and with $Re=100$}}
    \label{tab:truth_5005_errors_cavity}
\end{table}

\begin{table}[ht]
    \centering
    
    \begin{tabular}{|c|c|c|}
      \hline
        \textbf{Iteration} & \textbf{Functional Value}  &\textbf{Gradient norm }    \\
        \hline  
         0 & $2.5 \cdot 10^{-2}$ & $4.1 \cdot 10^{-1}$ \\
         5 & $7.4 \cdot 10^{-5}$ & $1.4 \cdot 10^{-2}$  \\ 
         10 & $3.3 \cdot 10^{-6}$ & $9.1 \cdot 10^{-4}$  \\ 
         25 & $7.0 \cdot 10^{-7}$ & $3.9 \cdot 10^{-4}$  \\
         \hline  
    \end{tabular}
    \caption{Functional values and the gradient norm for the \reviewerA{FOM} optimisation solution at the parameter values $\bar U=1$, $\nu=0.1$ \reviewerA{and with $Re=10$}}
    \label{tab:truth_101_func_cavity}
\end{table}

\begin{table}[ht]
    \centering
    
    \begin{adjustbox}{max width=\textwidth}
    \begin{tabular}{|c|c|c|c|c|c|c|c|c|}
      \hline
         \textbf{Iteration}  &\multicolumn{2}{|c|}{\textbf{Abs. error $u_{N}$}} &\multicolumn{2}{|c|}{\textbf{Rel. error $u_{N}$}} &\multicolumn{2}{|c|}{\textbf{Abs. error $p_{N}$} }&\multicolumn{2}{|c|}{\textbf{Rel. error $p_{N}$}}   \\ \hline  
         & $\Omega_1$&$\Omega_2$& $\Omega_1$&$\Omega_2$& $\Omega_1$&$\Omega_2$& $\Omega_1$&$\Omega_2$\\ \hline
         0 & 0.0668 & 0.0589 & 1.0000 & 0.2416 & 0.0349 & 0.0411 & 1.0000 & 0.0956 \\
        5 & 0.0032 & 0.0028 & 0.0483 & 0.0114 & 0.0036 & 0.0036 & 0.1100 & 0.0084 \\
        10 & 0.0006 & 0.0006 & 0.0095 & 0.0027 & 0.0024 & 0.0023 & 0.0733 & 0.0054 \\
        25 & 0.0005 & 0.0005 & 0.0069 & 0.0019 & 0.0021 & 0.0021 & 0.0663 & 0.0048 \\
         \hline  
    \end{tabular}
    \end{adjustbox}
     \caption{Absolute and relative errors of the optimisation \reviewerA{FOM} solution with respect to the monolithic solution at the parameter value $\bar U=1$, $\nu=0.1$ \reviewerA{and with $Re=10$}}
    \label{tab:truth_101_errors_cavity}
\end{table}

\begin{figure}[H]
    \centering
    \begin{subfigure}[b]{0.32\textwidth}
        \includegraphics[width=\textwidth]{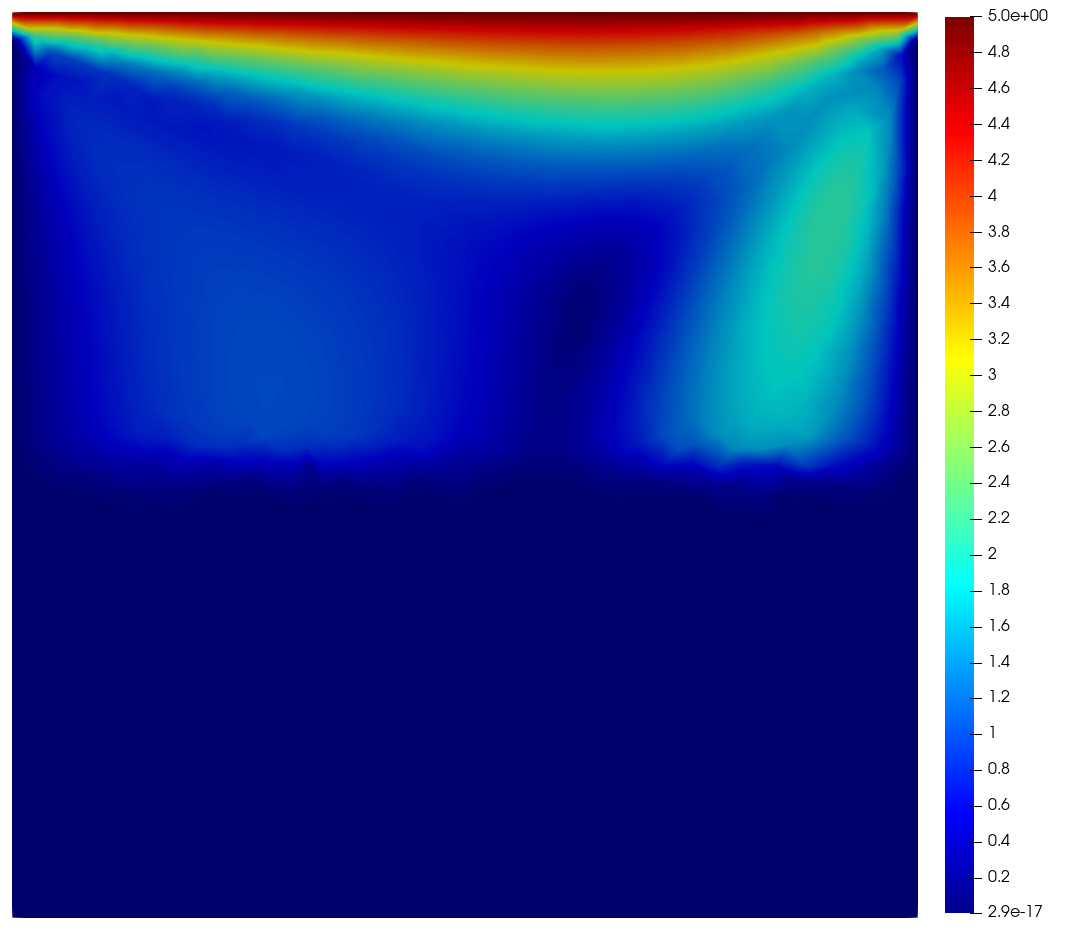}
        \caption{Iteration 0}
         \label{fig:reduced_u_5005_0_cavity}
    \end{subfigure}
    \hfill
    \begin{subfigure}[b]{0.32\textwidth}
        \includegraphics[width=\textwidth]{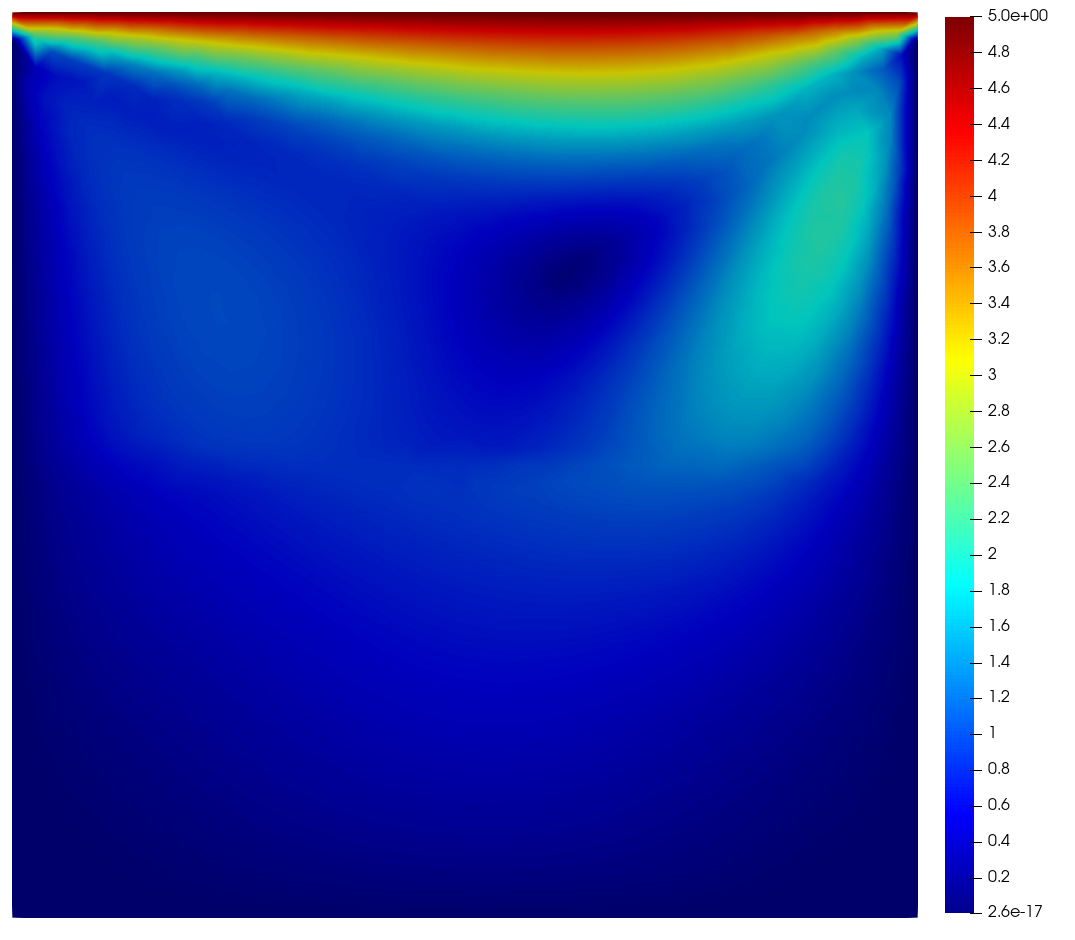}
         \caption{Iteration 3}
         \label{fig:reduced_u_5005_3_cavity}
         
    \end{subfigure}
    \hfill
    \begin{subfigure}[b]{0.32\textwidth}
        \includegraphics[width=\textwidth]{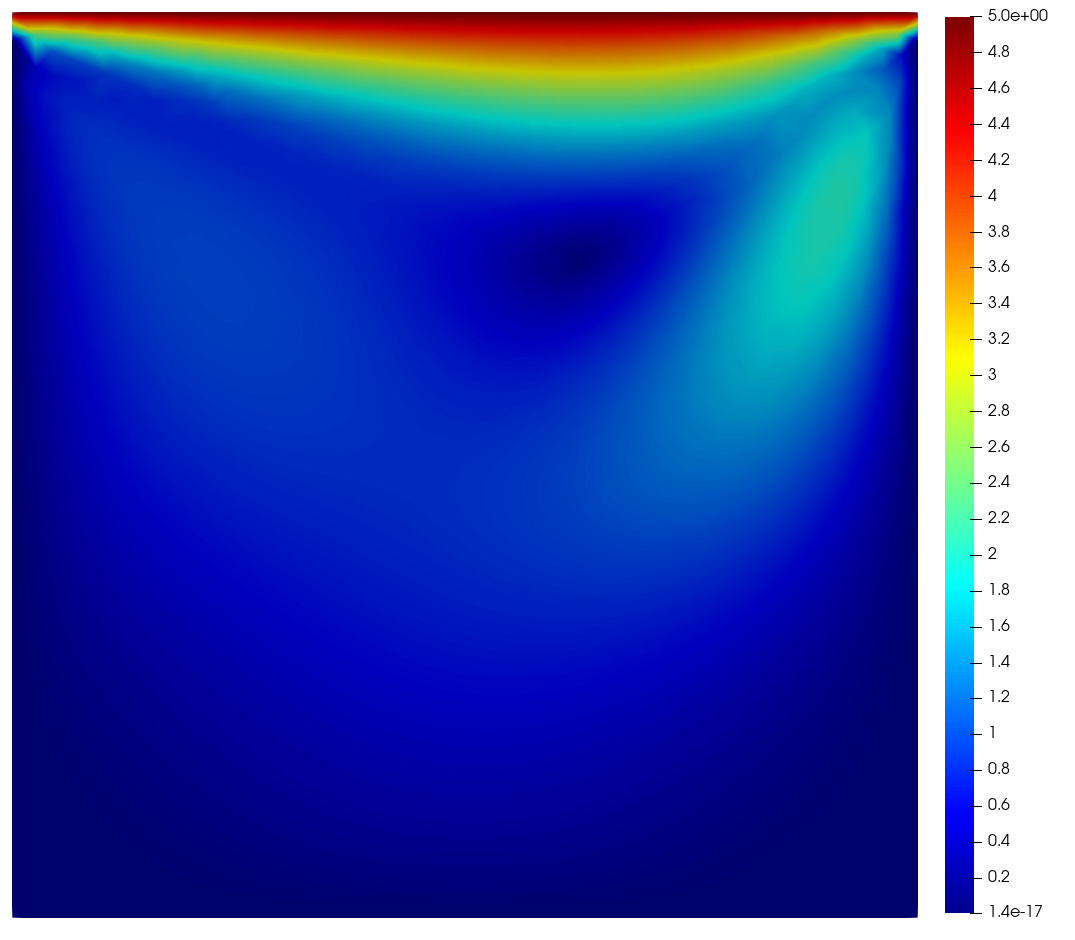}
         \caption{Iteration 10}
         \label{fig:reduced_u_5005_15_cavity}

    \end{subfigure}
        
    \caption{Reduced-order solution for the velocities $u_1$ and $u_2$. Values of the parameters $\bar U=5$, $\nu=0.05$ \reviewerA{and with $Re=100$}. \reviewerA{Number of POD modes: 10 - for each state variable, each supremiser and the control, 15 -- for both adjoint velocities}}
    \label{fig:reduced_u_5005_cavity}
\end{figure}

\begin{figure}[H]
    \centering
    \begin{subfigure}[b]{0.32\textwidth}
        \includegraphics[width=\textwidth]{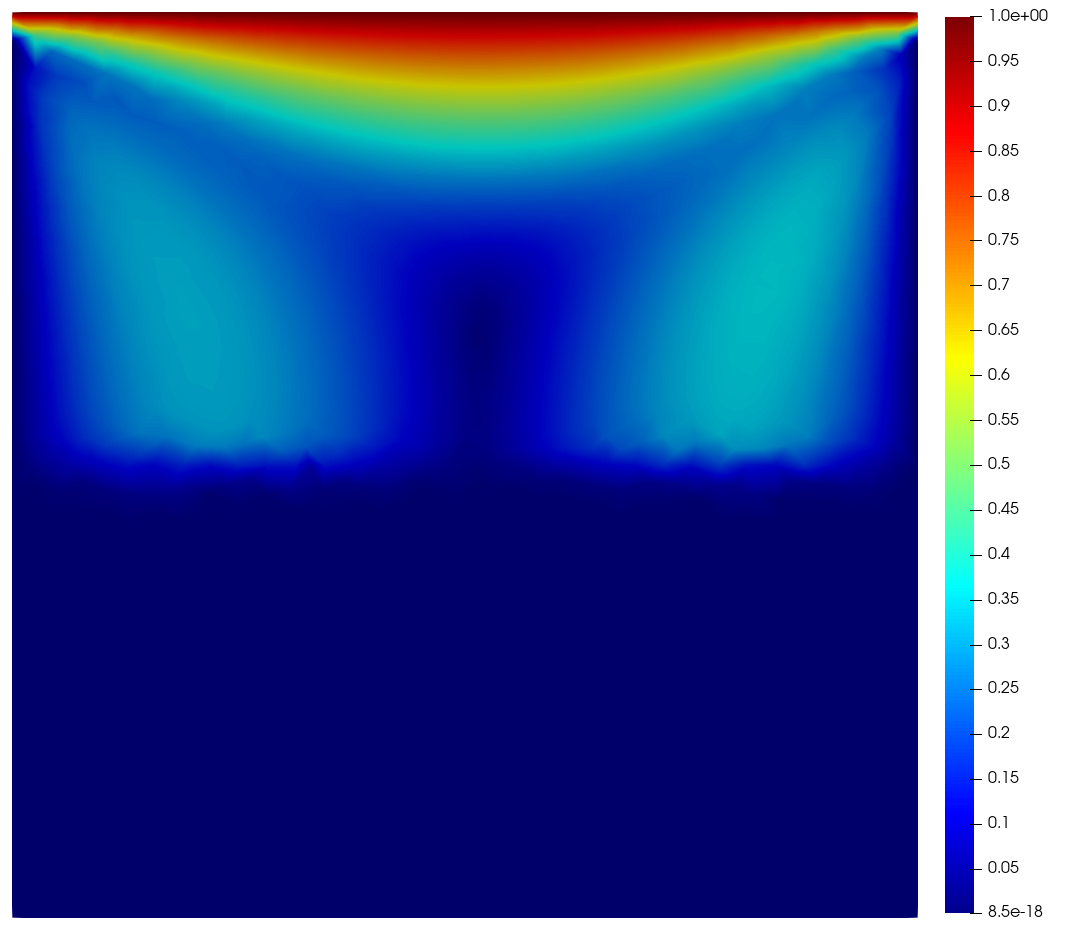}
        \caption{Iteration 0}
         \label{fig:reduced_u_101_0_cavity}
    \end{subfigure}
    \hfill
    \begin{subfigure}[b]{0.32\textwidth}
        \includegraphics[width=\textwidth]{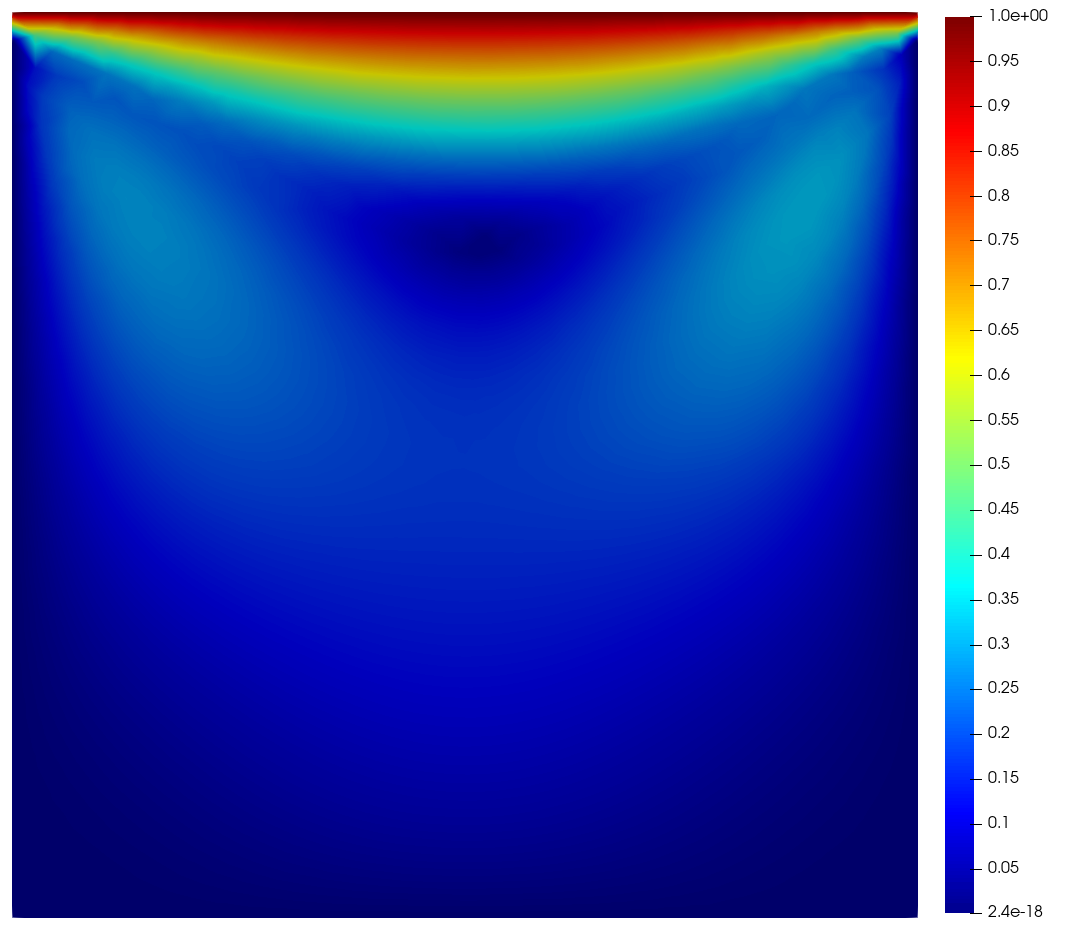}
         \caption{Iteration 3}
         \label{fig:reduced_u_101_3_cavity}
         
    \end{subfigure}
    \hfill
    \begin{subfigure}[b]{0.32\textwidth}
        \includegraphics[width=\textwidth]{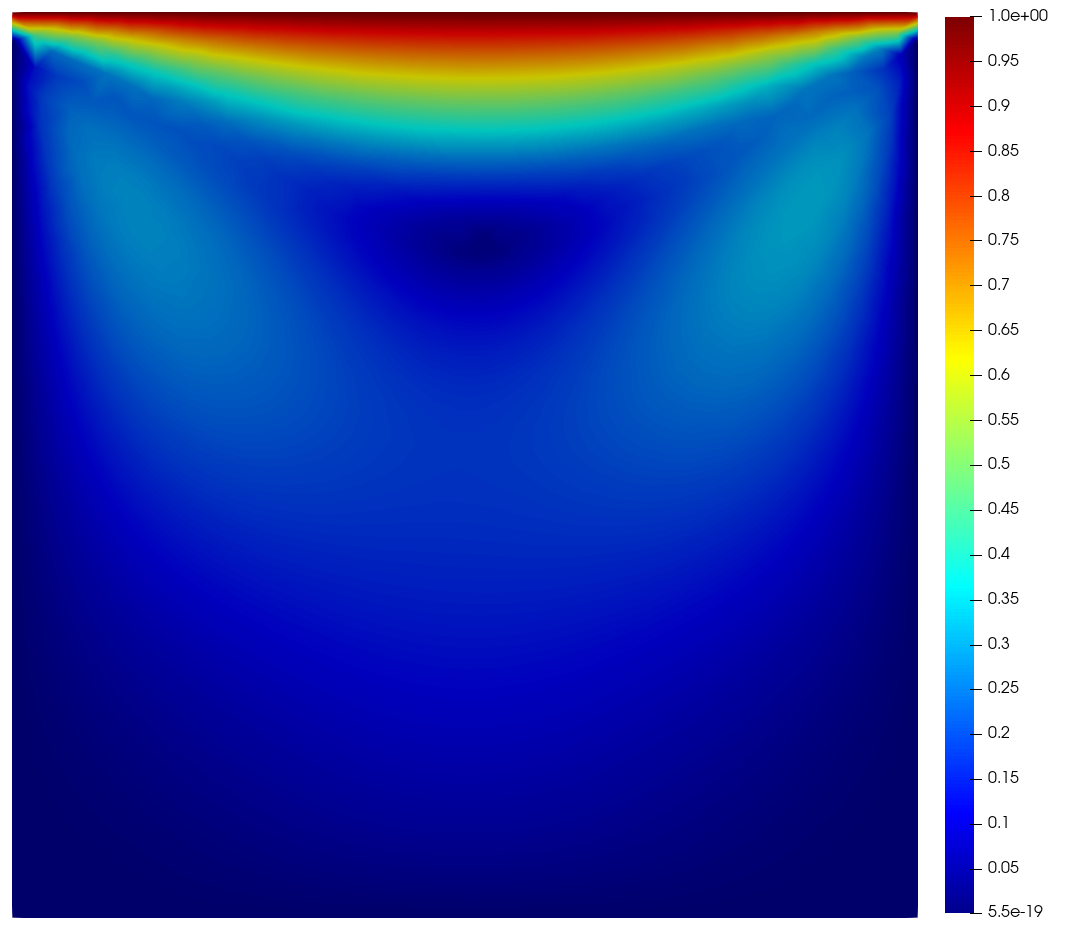}
         \caption{Iteration 10}
         \label{fig:reduced_u_101_10_cavity}

    \end{subfigure}
        
    \caption{Reduced-order solution for the velocities $u_1$ and $u_2$. Values of the parameters $\bar U=1$, $\nu=0.1$ \reviewerA{and with $Re=10$}. \reviewerA{Number of POD modes: 10 - for each state variable, each supremiser and the control, 15 -- for both adjoint velocities}}
    \label{fig:reduced_u_101_cavity}
\end{figure}

Figures \ref{fig:reduced_u_5005_cavity} -- \ref{fig:reduced_u_101_cavity} represent the reduced--order solutions for two different values of the parameters $(\bar U, \nu)=(5,0.05)$ and $(\bar U, \nu)=(1,0.1)$. For both cases, we choose the following number of the reduced basis functions: $N_{u_1} = N_{s_1} =  N_{p_1} = N_{u_2} = N_{s_2} =  N_{p_2} = N_g = 10$, whereas for the adjoint velocities we choose $N_{\xi_1} = N_{\xi_2} = 15$. As it was mentioned before we use a higher number for the adjoint variables $\xi_1$ and $\xi_2$ since they show much slower decay of the singular values (see Figure \ref{fig:singlular_values_cavity}). Figure \ref{fig:reduced_u_5005_cavity}  shows the intermediate solutions at iteration 0, 3 and 15 for the fluid velocities $u_1$ and $u_2$ corresponding to the parameter value $(\bar U, \nu)=(5,0.05)$, and Figure \ref{fig:reduced_u_101_cavity} shows the velocities $u_1$ and $u_2$ for the parameter value $(\bar U, \nu)=(1,0.1)$. The final solutions are taken to be the 10-iteration optimisation solution.

\begin{table}[ht]
    \centering
    
    \begin{tabular}{|c|c|c|}
      \hline
        \textbf{Iteration} & \textbf{Functional Value}  &\textbf{Gradient norm }    \\
        \hline  
        0 & $4.8 \cdot 10^{-1}$ & $3.153 $ \\
         3 & $2.4 \cdot 10^{-2}$ & $1.634 $  \\ 
         10 & $7.2 \cdot 10^{-5}$ & $0.023 $  \\
         \hline  
    \end{tabular}
    \caption{Functional values and the gradient norm for the \reviewerA{ROM} optimisation solution at parameter values $\bar U=5$, $\nu=0.05$ \reviewerA{and with $Re=100$}}
    \label{tab:reduced_5005_func_cavity}
\end{table}

\begin{table}[ht]
    \centering
    
        \begin{adjustbox}{max width=\textwidth}
    \begin{tabular}{|c|c|c|c|c|c|c|c|c|}
      \hline
         \textbf{Iteration}  &\multicolumn{2}{|c|}{\textbf{Abs. error $u_{N}$}} &\multicolumn{2}{|c|}{\textbf{Rel. error $u_{N}$}} &\multicolumn{2}{|c|}{\textbf{Abs. error $p_{N}$} }&\multicolumn{2}{|c|}{\textbf{Rel. error $p_{N}$}}   \\ \hline  
         & $\Omega_1$&$\Omega_2$& $\Omega_1$&$\Omega_2$& $\Omega_1$&$\Omega_2$& $\Omega_1$&$\Omega_2$\\ \hline
        0 & 0.3411 & 0.1796 & 1.0000 & 0.1523 & 0.2431 & 0.2519 & 1.0000 & 0.1864 \\
        3 & 0.0512 & 0.0552 & 0.1501 & 0.0468 & 0.0531 & 0.0646 & 0.3634 & 0.0478 \\
        10 & 0.0050 & 0.0056 & 0.0147 & 0.0047 & 0.0139 & 0.0139 & 0.0956 & 0.0103 \\
         \hline  
    \end{tabular}
    \end{adjustbox}
    \caption{Absolute and relative errors of the \reviewerA{ROM} optimisation solution with respect to the monolithic solution at the parameter values $\bar U=5$, $\nu=0.05$ \reviewerA{and with $Re=100$}}
    \label{tab:reduced_5005_errors_cavity}
\end{table}

We present additional details in Tables \ref{tab:reduced_5005_func_cavity} - \ref{tab:reduced_101_errors_cavity}. In particular, in Tables \ref{tab:reduced_5005_func_cavity} and \ref{tab:reduced_101_func_cavity} we list the values for the functional $\mathcal J_\gamma$ and the $L^2(\Gamma_0)$-norm  of the gradient $\frac{d \mathcal J_\gamma}{dg}$ at the different iteration of the optimisation procedure, while Table \ref{tab:reduced_5005_errors_cavity} and Table \ref{tab:reduced_101_errors_cavity} contain the $L^2$-relative errors with respect to the monolithic (the entire--domain) solutions $u_h, p_h$.

Analyzing the results, we are able to see that the reduced basis method gives us a solution as accurate as the high--fidelity model. The reduced--order approximation of the optimisation problem at hand allowed us to reduce the dimension of the high-fidelity optimisation functional by more than 10-20 times and enabled us to use half optimisation algorithm iterations (each optimisation step requires at least one solve of the state and the adjoint equations).

\begin{table}[ht]
    \centering
    
    \begin{tabular}{|c|c|c|}
      \hline
        \textbf{Iteration} & \textbf{Functional Value}  &\textbf{Gradient norm }    \\
        \hline  
         0 & $2.6 \cdot 10^{-2}$ & $2.6 \cdot 10^{-1}$ \\
         3 & $1.5 \cdot 10^{-5}$ & $1.0 \cdot 10^{-2}$  \\ 
        10 & $7.1 \cdot 10^{-7}$ & $1.2 \cdot 10^{-3}$  \\
         \hline  
    \end{tabular}
    \caption{Functional values and the gradient norm for the \reviewerA{ROM} optimisation solution at the parameter values $\bar U=1$, $\nu=0.1$ \reviewerA{and with $Re=10$}}
    \label{tab:reduced_101_func_cavity}
\end{table}

\begin{table}[ht]
    \centering
    
        \begin{adjustbox}{max width=\textwidth}
    \begin{tabular}{|c|c|c|c|c|c|c|c|c|}
      \hline
         \textbf{Iteration}  &\multicolumn{2}{|c|}{\textbf{Abs. error $u_{N}$}} &\multicolumn{2}{|c|}{\textbf{Rel. error $u_{N}$}} &\multicolumn{2}{|c|}{\textbf{Abs. error $p_{N}$} }&\multicolumn{2}{|c|}{\textbf{Rel. error $p_{N}$}}   \\ \hline  
         & $\Omega_1$&$\Omega_2$& $\Omega_1$&$\Omega_2$& $\Omega_1$&$\Omega_2$& $\Omega_1$&$\Omega_2$\\ \hline
        0 & 0.0668 & 0.0591 & 1.0000 & 0.2424 & 0.0349 & 0.0403 & 1.0000 & 0.0936 \\
        3 & 0.0010 & 0.0019 & 0.0155 & 0.0076 & 0.0024 & 0.0020 & 0.0752 & 0.0047 \\
        10 & 0.0004 & 0.0004 & 0.0066 & 0.0017 & 0.0020 & 0.0019 & 0.0621 & 0.0045 \\
         \hline  
    \end{tabular}
    \end{adjustbox}
     \caption{Absolute and relative errors of the \reviewerA{ROM}  optimisation solution with respect to the monolithic solution at the parameter values $\bar U=1$, $\nu=0.1$ \reviewerA{and with $Re=10$}}
    \label{tab:reduced_101_errors_cavity}
\end{table}

\begin{figure}
    \centering
    \begin{subfigure}[b]{0.49\textwidth}
        \includegraphics[width=\textwidth]{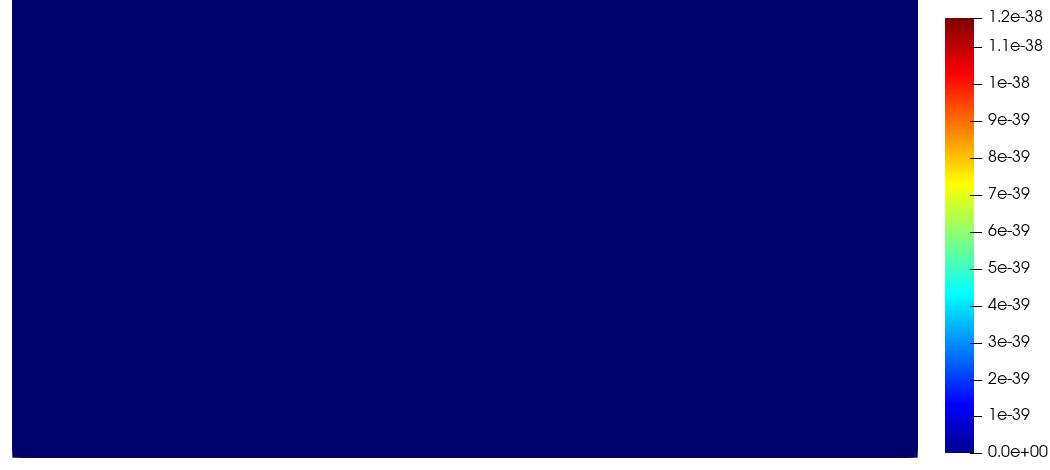}
        \caption{Iteration 0}
         \label{fig:reduced_u1_5005_0_cavity}
    \end{subfigure}
    \hfill
    \begin{subfigure}[b]{0.49\textwidth}
        \includegraphics[width=\textwidth]{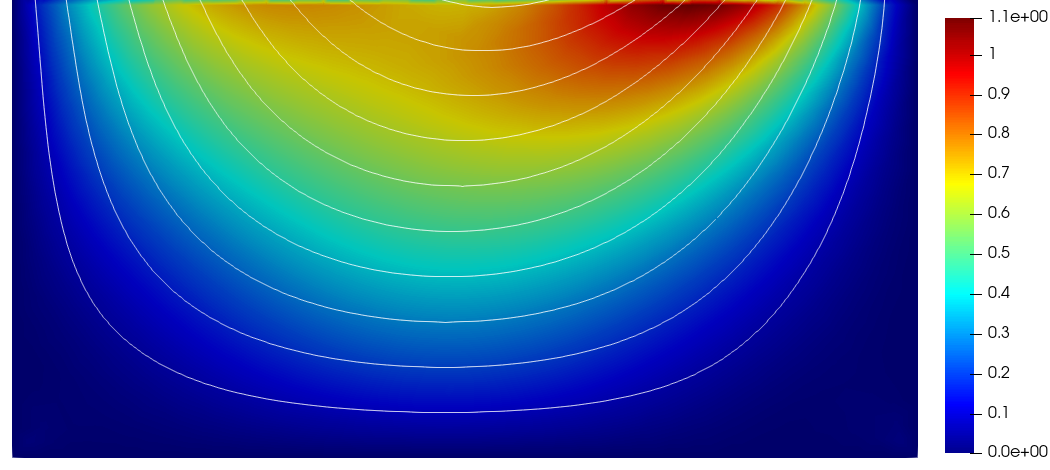}
         \caption{Iteration 5}
         \label{fig:reduced_u1_5005_5_cavity}

    \end{subfigure}
    
    \begin{subfigure}[b]{0.49\textwidth}
        \includegraphics[width=\textwidth]{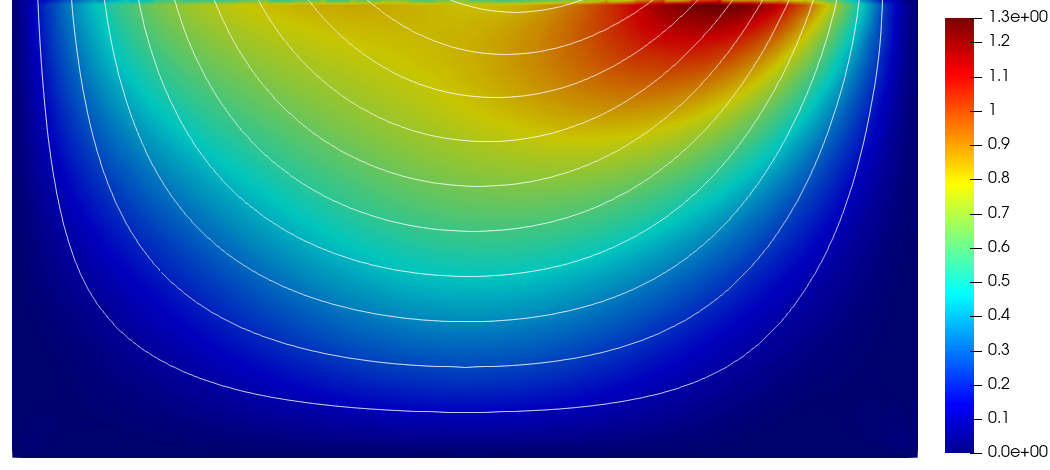}
        \caption{Iteration 10}
         \label{fig:reduced_u1_5005_10_cavity}
    \end{subfigure}
    \hfill
    \begin{subfigure}[b]{0.49\textwidth}
        \includegraphics[width=\textwidth]{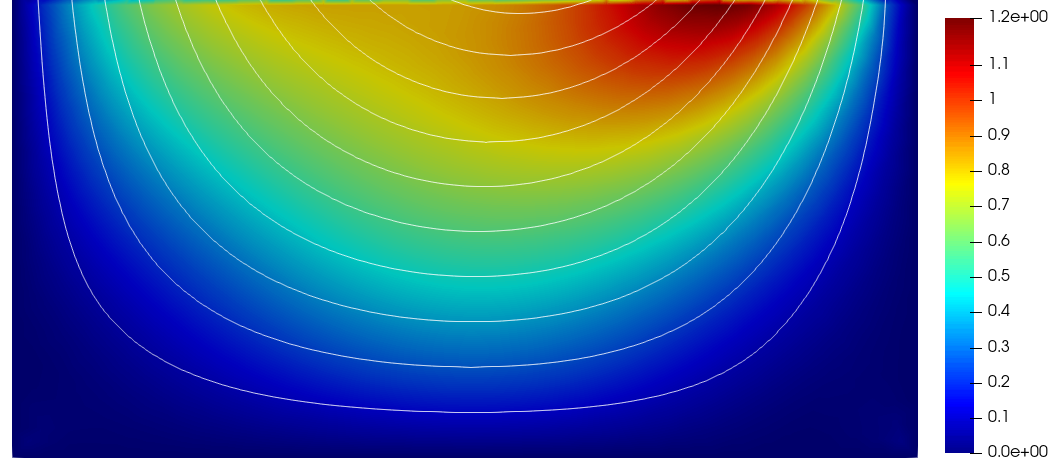}
         \caption{Iteration 25}
         \label{fig:reduced_u1_5005_25_cavity}

    \end{subfigure}
        
    \caption{Reduced--order solution for the velocity $u_1$. Values of the parameters $\bar U=5$, $\nu=0.05$ \reviewerA{and with $Re=100$}}
    \label{fig:reduced_u1_5005_cavity}
\end{figure}

\begin{figure}
    \centering
    \begin{subfigure}[b]{0.49\textwidth}
        \includegraphics[width=\textwidth]{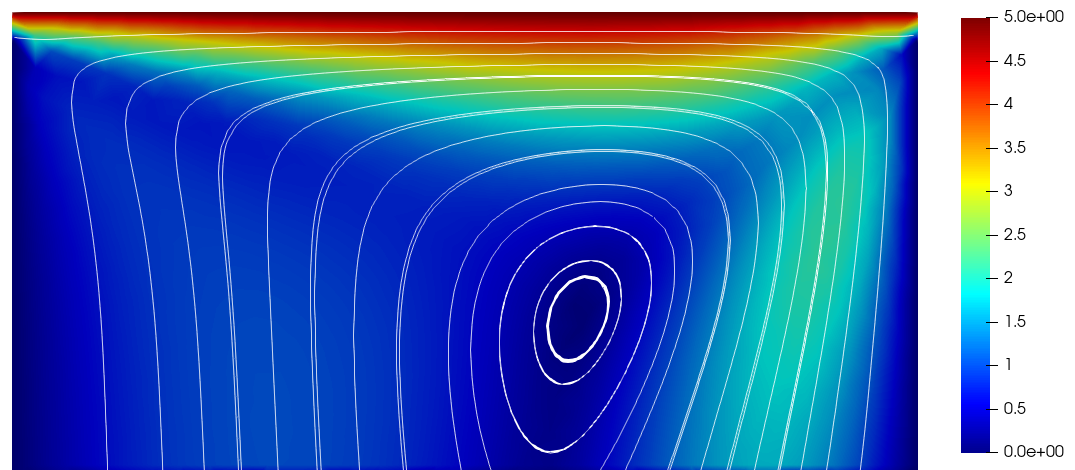}
        \caption{Iteration 0}
         \label{fig:reduced_u2_5005_0_cavity}
    \end{subfigure}
    \hfill
    \begin{subfigure}[b]{0.49\textwidth}
        \includegraphics[width=\textwidth]{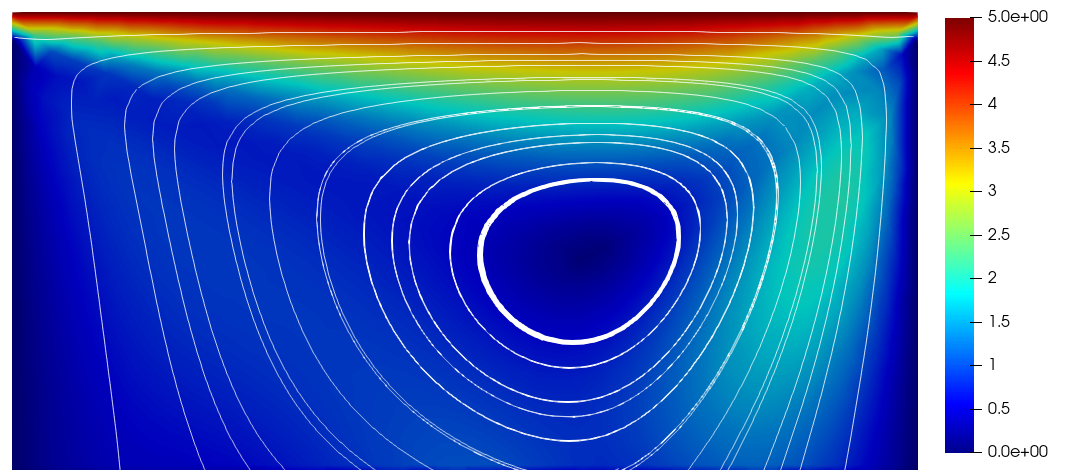}
         \caption{Iteration 5}
         \label{fig:reduced_u2_5005_5_cavity}

    \end{subfigure}
    
    \begin{subfigure}[b]{0.49\textwidth}
        \includegraphics[width=\textwidth]{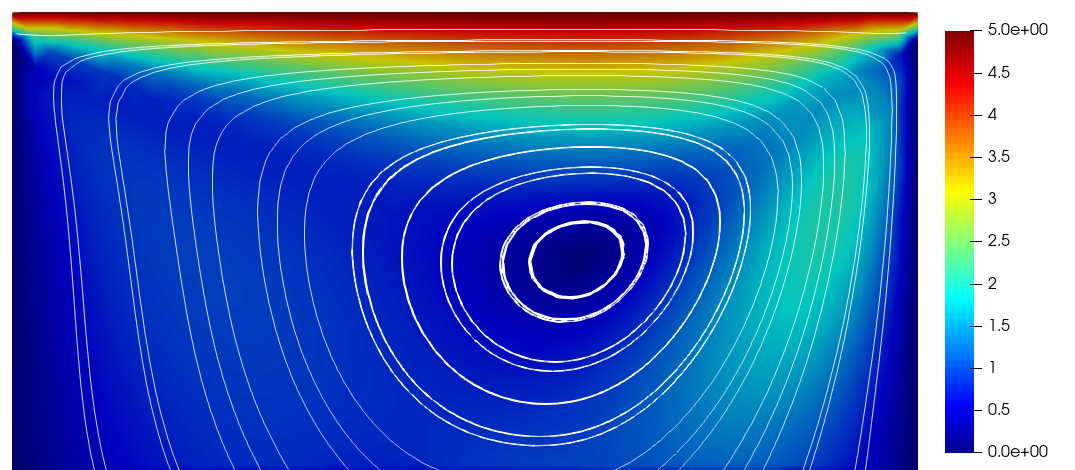}
        \caption{Iteration 10}
         \label{fig:reduced_u2_5005_10_cavity}
    \end{subfigure}
    \hfill
    \begin{subfigure}[b]{0.49\textwidth}
        \includegraphics[width=\textwidth]{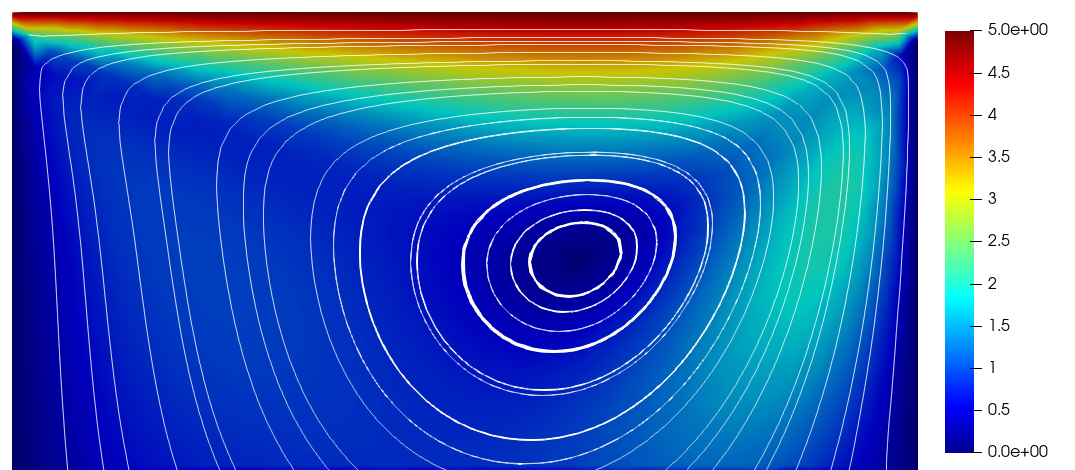}
         \caption{Iteration 25}
         \label{fig:reduced_u2_5005_25_cavity}

    \end{subfigure}
        
    \caption{Reduced--order solution for the velocity $u_2$. Values of the parameters $\bar U=5$, $\nu=0.05$ \reviewerA{and with $Re=100$}}
    \label{fig:reduced_u2_5005_cavity}
\end{figure}

In order to provide more visually representable results (the scale of the solution on the subdomains $\Omega_1$ and $\Omega_2$ has a few orders of the difference in the magnitude), we provide the graphs of the velocities $u_1$ and $u_2$ separately in Figures \ref{fig:reduced_u1_5005_cavity} and \ref{fig:reduced_u2_5005_cavity}.  \reviewerB{Additionally, in Table \ref{tab:fom_vs_rom_cavity} we provide a comparison between full--order and reduced--order models in terms of the relative errors between ROM solutions with respect to the corresponding FOM solutions. The considerations drawn in the previous section are valid also for this test case.}


\begin{table}[ht]
    \centering
    \begin{tabular}{|c|c|c|c|c|c|}
      \hline
      \multicolumn{2}{|c|}{\textbf{Parameter value}}  &\multicolumn{2}{|c|}{\textbf{Velocity relative error}} &\multicolumn{2}{|c|}{\textbf{Pressure relative error}}   \\ \hline  
      $\bar U$  & $\nu$ &  $\Omega_1$ & $\Omega_2$ & $\Omega_1$  & $\Omega_2$  \\  \hline  
      $1$ & $0.1$ & $0.020$ & $0.003$ & $0.014$ & $0.0007$ \\ 
       $5$ & $0.05$ & $0.040$ &  $0.005$ & $0.013$ & $0.002$ \\ \hline
    \end{tabular}
     \caption{\reviewerB{Relative errors between FOM and ROM solutions (in terms of $H^1$--norm for the velocity fields and $L^2$--norm for the pressure fields)}}
    \label{tab:fom_vs_rom_cavity}
\end{table}

\reviewerA{
\begin{remark}
 In both numerical cases presented above, it might seem that due to the fact that the non--homogeneous Dirichlet boundary condition is present only on the boundary of one of the subdomains only a few corrections are needed on this subdomain. On the other hand, this is true only for the velocity field, as it can be seen in the tables listing the errors (for instance in Table \ref{tab:truth_1_errors_bfs}). Indeed, the errors for the pressure on those subdomains are higher than on the other one. Regarding the cavity flow, our original idea was to split the domain vertically, but in that case, the convergence even at full--order level was much slower, hence, we opted for the horizontal split. 
\end{remark}
}
\reviewerB{
\begin{remark}[High Reynolds simulations]
 Also for this test case, the range of Reynolds number for which the DD solver converges is stricter than the one where the monolithic solver provides a solution. The reason is that the optimisation algorithms are very sensitive to the initial guess, and the authors suspect that some further stabilisation techniques should be used.
\end{remark}
}
\newpage
\section{Conclusions}
\label{conclusions}

In this work, we proposed a reduced--order model for the optimisation--based domain decomposition formulation of the parameter-dependent stationary incompressible Navier--Stokes equations. 

The original problem cast into the optimisation--based domain--decomposition framework leads to the optimal control problem aimed at minimising the coupling error at the interface; the problem, then, has been tackled using an iterative gradient--based optimisation algorithm, which allowed us to obtain a complete separation of the solvers on different subdomains.

On the reduced--order level, we have managed to build a model for which the generation of the reduced basis spaces is carried out separately in each subdomain and for each component of the problem solution. Furthermore, as the numerical results show, the reduction of the optimal--control problem can be observed not only in the dimensions of the different components of the problem, i.e., of the functional, the state and the adjoint equations but also in the number of the iterations of the optimisation algorithm. 

As it has been mentioned in the paper, the aforementioned techniques could be promising in the context of more complex time--dependent problems and, more importantly, multi--physics problems, where either pre-existing solvers can be used on each subcomponent or we do not have direct access to the codes. In particular, in future, we are planning to extend the methodology presented in this paper to \reviewerA{problems with several sub-domains}, to nonstationary fluid--dynamics problems and, eventually, to Fluid--Structure interaction problems. Moreover, this approach can be applied also to more complicated problems, where different types of numerical models are used in different subdomains.

\section*{Acknowledgements}
This work was supported by the European Union's Horizon 2020 research and innovation programme under the Marie Sklodowska-Curie Actions [grant agreement 872442]
(ARIA, Accurate Roms for Industrial Applications)  and by PRIN ``Numerical Analysis for Full and Reduced Order Methods for Partial Differential Equations'' (NA-FROM-PDEs) project. MN acknowledges the support of the Austrian Science Fund (FWF) project F65 ``Taming complexity in Partial Differential Systems'' and the Austrian Science Fund (FWF) project P 33477. DT has been funded by a SISSA Mathematical fellowship within Italian Excellence Departments initiative by Ministry of University and Research.
FB also thanks the project ``Reduced order modelling for numerical simulation of partial differential equations'' funded by Università Cattolica del Sacro Cuore.

\newpage
\bibliographystyle{abbrv} 
\bibliography{refs}

\end{document}